\documentclass[11pt,oneside,english]{amsart}
\usepackage{amstext}
\usepackage{amsthm,amssymb}
\usepackage[T1]{fontenc}
\usepackage[latin9]{inputenc}
\usepackage{geometry}
\usepackage{mathrsfs}
\usepackage{wasysym}
\usepackage{graphicx,epstopdf}
\usepackage[colorlinks,bookmarks,linkcolor=black,citecolor=black]{hyperref} 
\usepackage{enumitem}
\usepackage[english,algoruled,lined,noresetcount,norelsize]{algorithm2e}
\usepackage{bbm,bm}

\makeatletter
\numberwithin{equation}{section}
\numberwithin{figure}{section}

\newcommand{\By}[2]{\overset{\mbox{\tiny{#1}}}{#2}}
\newcommand{\ByRef}[2]{   \By{\eqref{#1}}{#2} }

\newcommand{\eqByRef}[1]{ \ByRef{#1}{=} }

\newcommand{\leByRef}[1]{ \ByRef{#1}{\le} }

\newtheorem{thm}{Theorem}
\newtheorem{lem}[thm]{Lemma}

\newtheorem{setting}[thm]{Setting}
\newtheorem{meta}[thm]{Meta-Conjecture}
\newtheorem{qu}[thm]{Question}

\newcounter{counterProbint}
\newcounter{counterDiet}
\newcounter{counterProbedge}
\newtheorem{claimProbint}{Claim}[counterProbint]
\newtheorem{claimDiet}{Claim}[counterDiet]
\newtheorem{claimProbedge}{Claim}[counterProbedge]

\newtheorem{defn}[thm]{Definition}

\theoremstyle{remark}

\newcommand{\oldqed}{}
\newcommand{\qedClaim}{\hfill\scalebox{.6}{$\Box$}}
\renewcommand{\epsilon}{\varepsilon}
\newenvironment{claimproof}[1][Proof]{
  \renewcommand{\oldqed}{\qedsymbol}
  \renewcommand{\qedsymbol}{\qedClaim}
  \begin{proof}[#1]
}{
  \end{proof}
  \renewcommand{\qedsymbol}{\oldqed}
} 
\newcommand{\JUSTIFY}[1]{\mbox{\tiny{(#1)}}\quad}
\newcommand{\ENHANCED}[1]{}
\newcommand{\EXP}{\mathbb{E}}
\newcommand{\PROB}{\mathbb{P}}
\newcommand{\NBH}{\mathsf{N}}
\newcommand{\LNBH}{\NBH^{-}}
\newcommand{\LEFTDEG}{\deg^{-}}
\newcommand{\CANDSET}{C}
\newcommand{\im}{\operatorname{im}}

\newcommand{\cA}{\mathcal{A}}
\newcommand{\cB}{\mathcal{B}}
\newcommand{\cE}{\mathcal{E}}
\newcommand{\cY}{\mathcal{Y}}
\newcommand{\cD}{\mathcal{D}}
\newcommand{\cF}{\mathcal{F}}

\newcommand{\DietEvent}{\mathsf{DietE}}
\newcommand{\CoverEvent}{\mathsf{CoverE}}
\newcommand{\CodietEvent}{\mathsf{CoDietE}}

\newcommand{\eps}{\varepsilon}
\renewcommand{\rho}{\varrho}
\renewcommand{\subset}{\subseteq}

\newcommand{\RandomEmbedding}{\emph{RandomEmbedding}}
\newcommand{\ModRandomEmbedding}{\emph{ModifiedRandomEmbedding}}
\newcommand{\ModmodRandomEmbedding}{\emph{MoreModifiedRandomEmbedding}}
\newcommand{\PackingProcess}{\emph{PackingProcess}}
\newcommand{\hist}{\mathscr{H}}
\newcommand{\histens}{\mathscr{L}}

\newcommand{\Var}{\mathrm{Var}}
\newcommand{\AlgMap}{\hookrightarrow}

\newcommand{\corrig}[1]{#1}

\begin{document}

\title{Packing degenerate graphs}

\author{Peter Allen}

\author{Julia Böttcher}

\address{Department of Mathematics, London School of Economics, Houghton Street,
London,\linebreak WC2A~2AE, UK}

\email{\{p.d.allen, j.boettcher\}@lse.ac.uk}

\author{Jan Hladký}
\address{Institute of Mathematics of the Czech Academy of Sciences. \v Zitn\'a 25, 110 00, Praha, Czechia. The Institute of Mathematics of the Czech Academy of Sciences is supported by RVO:67985840.}
\email{honzahladky@gmail.com}
\thanks{The research leading to these results has received funding from the
People Programme (Marie Curie Actions) of the European Union's Seventh
Framework Programme (FP7/2007-2013) under REA grant agreement number
628974.}

\author{Diana Piguet}

\address{Institute of Computer Science of the Czech Academy of Sciences, Pod Vodárenskou
v\v{e}ží 2, 182~07 Prague, Czech Republic. With institutional support
RVO:67985807.}
 \thanks{PA was partially supported by the EPSRC, grant number EP/P032125/1.}
 \thanks{JB was partially supported by the EPSRC, grant number EP/R00532X/1.}
\thanks{JH and DP were supported by the Czech Science Foundation, grant number
GJ16-07822Y}

\email{piguet@cs.cas.cz}
\begin{abstract}
Given $D$ and $\gamma>0$, whenever $c>0$ is sufficiently small and $n$ sufficiently large, if $\mathcal{G}$ is a family of $D$-degenerate
graphs of individual orders at most $n$, maximum degrees at most $\tfrac{cn}{\log n}$, and total number of edges at most $(1-\gamma)\binom{n}{2}$, then $\mathcal{G}$ packs into the complete graph $K_{n}$. Our
proof proceeds by analysing a natural random greedy packing algorithm. 

\corrig{This version of the manuscript corrects a small error that appeared in the published version [Adv Math, 354 (2019), 106739].}
\end{abstract}

\maketitle

\section{Introduction}

A \emph{packing} of a family $\mathcal{G}=\{G_1,\dots,G_k\}$ of graphs into a graph $H$ is a colouring of the edges of $H$ with the colours $0,1,\dots,k$ such that the edges of colour $i$ form an isomorphic copy of $G_i$ for each $1\le i\le k$. The packing is \emph{perfect} if no edges have colour $0$. We will often say an edge is \emph{covered} in a packing if it has colour at least $1$, and \emph{uncovered} if it has colour zero.

Packing problems have been studied in graph theory for several decades. Many classical theorems and conjectures of extremal graph theory can be written as packing problems. For example, Tur\'an's theorem can be read as the statement that if the $n$-vertex $G$ does not have too many edges (depending on $r$), then $G$ and $K_r$ pack into $K_n$. Putting extremal statements into this context often suggests interesting generalisations, such as asking for packings of more graphs. However packings in this context are usually very far from being perfect packings, with a large fraction of $E(H)$ uncovered. By contrast, in this paper we are interested in \emph{near-perfect packings}, that is, packings in which $o\big(e(H)\big)$ edges are uncovered.

The first problems asking for perfect packings in graphs actually predate modern graph theory: Pl\"ucker~\cite{Pluecker} in 1835 found perfect packings of $\tfrac13\binom{n}{2}$ copies of $K_3$ into $K_n$ for various values of $n$, and more generally, Steiner~\cite{Steiner} in 1853 asked the following question (phrased then in set-theoretic terms).
\begin{qu}\label{qu:design}
 Given $2\le k\le r$, for which values of $n$ does the complete $k$-uniform hypergraph $K^{(k)}_n$ have a perfect packing with copies of $K^{(k)}_r$?
\end{qu}
A packing of this form is called a \emph{combinatorial design}. There are some simple divisibility conditions on $n$ which are necessary for an affirmative answer. Recently and spectacularly, Keevash~\cite{Kee:ExistenceOfDesigns} proved that for sufficiently large $n$ these conditions are also sufficient. This result was reproved, using a more combinatorial method, by Glock, K\"uhn, Lo and Osthus~\cite{GKLO:Designs}, who were also able to extend the result to pack with arbitrary fixed hypergraphs in~\cite{GKLO:Fdesigns}. A related problem tracing back to Kirkman~\cite{Kirkman} in 1846 asks for packings with copies of the $n$-vertex $K^{(k)}_r$-factor (Kirkman posed specifically the case $k=2$, $r=3$, asking for $K_n$ to be packed with $\tfrac{n-1}{6}$ copies of the graph consisting of $\tfrac{n}{3}$ disjoint triangles). Such packings are called \emph{resolvable designs}, and although Ray-Chaudhuri and Wilson~\cite{RCW} solved Kirkman's problem (Kirkman's designs exist if and only if $n$ is congruent to $3$ modulo $6$), in general the problem is wide open.

The focus of this paper is in packings of large connected graphs. In 1963 Ringel~\cite{Ringel1963} conjectured that if $T$ is any $(n+1)$-vertex tree, then $2n+1$ copies of $T$ pack into $K_{2n+1}$, and in 1976 Gy\'arf\'as~\cite{Gyarfas1978} made the Tree Packing Conjecture, that if $T_i$ is an $i$-vertex tree for each $1\le i\le n$ then $\{T_1,\dots,T_n\}$ packs into $K_n$. Note that both conjectures ask for perfect packings. These problems are both unsolved, although there are many partial results. It is easy (in both cases) to verify that the conjecture holds when the trees are all stars, or all paths. In both cases, the conjectures were also settled for some specific families of trees (see a rather outdated survey by
Hobbs~\cite{Hobbs1981}), but until recently there existed no general results.

Intuitively, perfect packing results are hard precisely because every edge must be used. If the graphs $\mathcal{G}$ were embedded in order to $H$, on coming to the last graph of $\mathcal{G}$ we would need to find that a hole is left in $H$ of precisely the right shape to accommodate it; this clearly requires some foresight in the packing. If some edges will remain uncovered at the end, this difficulty decreases. Bollob\'as~\cite{Bollobas1983} was the first to utilise this, making the observation that one can pack the $2^{-1/2}n$ smallest trees of the Tree Packing Conjecture, and assuming the Erd\H{o}s--S\'os Conjecture\footnote{The Erd\H{o}s--S\'os Conjecture states that if an $n$-vertex graph has more than $\frac12(k-1)n$ edges then it contains each tree of order $k+1$. A proof (of a slightly weaker form of) the Erd\H{o}s--S\'os Conjecture was announced by Ajtai, Koml\'os, Simonovits and Szemer\'edi in the 1990s.} even the $\sqrt{3}n/2$ smallest trees. More recently Balogh and Palmer~\cite{Balogh2013} showed that for large $n$ the $\tfrac14n^{1/3}$ largest trees pack, provided their maximum degree is at most $2n^{2/3}$, and without degree restriction that the $\tfrac1{10}n^{1/4}$ largest trees pack in $K_{n+1}$ (i.e.\ using an extra vertex). These results do not give near-perfect packings --- a significant fraction of the complete graph is uncovered --- but until recently they were the only general results on the Tree Packing Conjecture allowing high-degree trees.

The first approximate result on the tree packing conjectures is due to B\"ottcher, Hladk\'y, Piguet and Taraz~\cite{Boettcherb}, who showed that one can pack into $K_n$ any family of trees whose maximum degree is at most $\Delta$, whose order is at most $(1-\delta)n$, and whose total number of edges is at most $(1-\delta)\binom{n}{2}$, provided that $n$ is sufficiently large given $\Delta$ and $\delta>0$. This provides approximate versions of both Ringel's Conjecture and the Tree Packing Conjecture for bounded degree graphs. A flurry of generalisations followed, beginning with Messuti, R\"odl and Schacht~\cite{MeRoSch:PackingMinorClosed}, who showed that one can replace trees with graphs from any nontrivial minor-closed family (but still requiring the other conditions), and then by Ferber, Lee and Mousset~\cite{FeLeMou:PackingSpanning} who showed that the restriction to at most $(1-\delta)n$ vertex graphs is unnecessary. Then, Kim, K\"uhn, Osthus and Tyomkyn~\cite{KiKuOsTy:Packing} proved a near-perfect packing result for families of graphs with bounded maximum degree which are otherwise unrestricted. At last, Joos, Kim, K\"uhn and Osthus~\cite{JoKiKuOs:Packing} obtained exact solutions of both Ringel's conjecture and the Tree Packing conjecture when all trees have degree bounded by a constant $\Delta$ and $n$ is sufficiently large compared to $\Delta$. This is an impressive and difficult result: what remains (which, unfortunately, is almost all cases) is to consider trees with some vertices of large degree.

Generalising in the direction of removing the restriction to bounded degree graphs, Ferber and Samotij~\cite{FeSa:PackingTrees} showed two near-perfect packing results for trees, one for spanning trees of maximum degree $O\big(n^{1/6}\log^{-6}n\big)$, and one for almost spanning trees of maximum degree $O\big(n/\log n\big)$. The latter result also follows in the particular case of Ringel's Conjecture from the work of Adamaszek, Allen, Grosu, Hladk\'y~\cite{AAGH:AlmostAllTrees}. The focus of~\cite{AAGH:AlmostAllTrees}  is the so-called Graceful Tree Conjecture but there is a well-known observation that this conjecture would imply Ringel's Conjecture, see~\cite[Section~1.1]{AAGH:AlmostAllTrees}.

To state our main result, we need to define the \emph{degeneracy} of a graph $G$. An ordering of $V(G)$ is $D$-degenerate if every vertex has at most $D$ neighbours preceding it, and $G$ is $D$-degenerate if $V(G)$ has a $D$-degenerate ordering. Every graph from a non-trivial minor-closed class has bounded degeneracy. In particular, trees are 1-degenerate, planar graphs are 5-degenerate. Of course, every bounded-degree graph has automatically bounded degeneracy.

Our main result then reads as follows.
\begin{thm}
\label{thm:MAINunbounded}For each $\gamma>0$ and each $D\in\mathbb{N}$
there exists $c>0$ and a number $n_{0}$ such that the following holds for each
integer $n>n_{0}$. Suppose that $\left(G_{t}\right)_{t\in [t^{*}]}$
is a family of $D$-degenerate graphs, each of which has at most $n$
vertices and maximum degree at most $\tfrac{cn}{\log n}$. Suppose further that the
total number of edges of $\left(G_{t}\right)_{t\in [t^{*}]}$ is at
most $(1-\gamma)\binom{n}{2}$. Then $\left(G_{t}\right)_{t\in [t^{*}]}$
packs into $K_{n}$. 
\end{thm} 
Theorem~\ref{thm:MAINunbounded} thus strengthens the main results about packings into complete graphs from ~\cite{Boettcherb,MeRoSch:PackingMinorClosed,FeLeMou:PackingSpanning,KiKuOsTy:Packing,FeSa:PackingTrees}.\footnote{Some of these papers deal also with packings into non-complete graphs, and most of these results are summarised below.} The main  features of the result are that guest graphs may be spanning, expanding, and have very high maximum degree.

Moving away from packing into complete graphs, there are several classical conjectures which ask for packing results similar to the above when $K_n$ is replaced by a graph of sufficiently high minimum degree, perhaps with additional constraints (such as regularity). Advances have recently been made on several of these, especially by the Birmingham Combinatorics group (see for example~\cite{CsKueLoOsTr:HamDec,barber2015fractional,glock2016decomposition}). In particular, we should observe that the near-perfect packing for bounded degree graphs~\cite{KiKuOsTy:Packing} mentioned above actually works in the setting of $\eps$-regular partitions, which turned out to be necessary for the perfect packing results of~\cite{JoKiKuOs:Packing}.

Finally, in line with the current trend in extremal combinatorics of asking for random analogues of classical extremal theorems, one can ask for packing results when $K_n$ is replaced by a typical binomial random graph $\mathbb G(n,p)$. This is actually the focus of the paper of Ferber and Samotij~\cite{FeSa:PackingTrees}, and they are able to prove near-perfect packing results even in $\mathbb G(n,p)$ when $p$ is not much above the threshold for connectivity. Our approach also proves near-perfect packing results (for the same family of graphs) in sufficiently quasirandom graphs of any positive constant edge density (see Theorem~\ref{thm:maintechGENERAL}), and hence in Erd\H{o}s--R\'enyi random graphs (see Theorem~\ref{thm:packingIntoErdosRenyi}). It might be possible to modify our approach to work in somewhat sparse random graphs as well, but certainly not sparse enough to compete with~\cite{FeSa:PackingTrees}. 

Although our current progress with actually proving exact packing conjectures is limited, at least we have not found counterexamples. The existing conjectures point in the following direction.

\begin{meta} Let $\mathcal{G}$ be any family of sparse graphs, and $H$ be an $n$-vertex dense graph. If there is no simple obstruction to packing $\mathcal{G}$ into $H$, then a packing exists.
\end{meta}

Some obvious examples of obstructions include the total number of edges in the family $\mathcal{G}$ being larger than $e(H)$, or any graph in $\mathcal{G}$ having more vertices than $H$. Certainly more subtle obstructions exist. For example it is possible that the total number of edges in graphs of $\mathcal{G}$ equals $e(H)$, but all graphs in $\mathcal{G}$ have only vertices of even degree, while some vertices of $H$ have odd degree, so that there is a parity obstruction to packing $\mathcal{G}$ into $H$, or that $\mathcal{G}$ contains two graphs with vertices of degree $n-1$ (or more generally too many vertices of very high degree). More such examples exist, see for example the discussions in~\cite{Boettcherb} (Section~9.1) and~\cite{JoKiKuOs:Packing} (after Theorem~1.7). The meta-conjecture can be read as claiming that there is nevertheless a finite list.

Note that without restriction the problem of packing a given $\mathcal{G}$ into a given $H$ is NP-complete (the survey~\cite{Yus:Sur} gives several NP-completeness results of which the one in~\cite{DorTar} is arguably the most convincing), so in particular we do not expect to find any finite list of simple obstructions to the general packing problem. It follows that `dense' in the meta-conjecture cannot simply mean large edge-density: one can artificially boost edge density without changing the outcome of this decision problem by taking the disjoint union with a very large clique and adding large connected graphs to $\mathcal{G}$ which perfectly pack the very large clique. However a typical random, or quasirandom, graph seems to be a reasonable candidate for `dense', as does a graph with high minimum degree (in this case, the minimum degree bound must depend on parameters of the graphs $\mathcal{G}$ such as chromatic number, otherwise a reduction similar to the edge-density reduction exists).

Finally, on the topic of what constitutes a `sparse graph', observe that bounded degeneracy is a fairly common and unrestrictive notion. One might ask whether degeneracy growing as a function of $n$ is reasonable (of course, in Theorem~\ref{thm:MAINunbounded} one can have a very slowly growing function). However, observe that we do not know the answer to Question~\ref{qu:design} when $r$ grows superlogarithmically, even for $k=2$, and it seems reasonable to believe that the answer will often be `no' even when the simple divisibility conditions are met. It is less clear that the maximum degree restriction of Theorem~\ref{thm:MAINunbounded} is necessary, and we expect that it can at least be relaxed. However, with no degree restriction at all Theorem~\ref{thm:MAINunbounded} becomes false, see Section~\ref{ssec:limits}.

\subsection*{Proof outline and organisation of the paper}
Our proof of Theorem~\ref{thm:MAINunbounded} amounts to the analysis of a quite natural randomised algorithm. We first describe a procedure which works if each graph in $\mathcal{G}$ has order at most $(1-\delta)n$. We take graphs in $\mathcal{G}$ in succession. For each $G$, we embed vertex by vertex into $K_n$ in a degeneracy order, at each time embedding to a vertex of $K_n$ chosen uniformly at random subject to the constraints that we do not re-use a vertex previously used in embedding $G$, or an edge used in embedding a previous graph. This procedure succeeds with high probability, and after each stage of embedding a graph, the unused edges in $K_n$ are quasirandom (in a sense we will later make precise).

To allow for spanning graphs, we modify this slightly. We adjust the degeneracy order so that the last $\delta n$ vertices are independent and all have the same degree; this can be done while at worst doubling the degeneracy of the order. Then for each graph we follow the above procedure to embed the first $(1-\delta)n$ vertices, and finally complete the embedding arbitrarily using a matching argument. We will see that this last step is with high probability always possible. The only slight subtlety is that we have to split $E(K_n)$ into a very dense main part, whose edges we use only for the embedding of the first $(1-\delta)n$ vertices, and a sparse reservoir which we use only for the completion; we do this randomly.

This paper is organised as follows.
In Section~\ref{sec:prelim} we introduce martingale concentration inequalities needed for the analysis of our algorithm. We also establish some basic properties of degenerate graphs. In Section~\ref{sec:reduction} we state our main technical result (Theorem~\ref{thm:maintechGENERAL}) and show how to deduce Theorem~\ref{thm:MAINunbounded} from it. In Section~\ref{sec:prooftech} we describe in detail our packing algorithm, \PackingProcess, and outline the main steps of its analysis. We also state our main lemmas and show how they imply Theorem~\ref{thm:maintechGENERAL}. In Sections~\ref{sec:diet},~\ref{sec:maintainingquasirandomness} and~\ref{sec:completion} we prove these lemmas. Finally in Section~\ref{sec:concluding} we give some concluding remarks.

\corrig{
\subsection{About the current version of the manuscript}
The current version of the manuscript, which we made available in April~2022, corrects a small error in Definition~\ref{def:cover} (version of June 2021) and an insufficiently small choice of constants (April 2022). These errors appear in the published version~[Advances in Mathematics, Volume 354, 106739]. Calculations had to be adjusted appropriately throughout the paper. The paper was not updated otherwise. In particular, the introduction and cited literature represent the state at the time of publication.
}

\section{Notation and preliminaries}
\label{sec:prelim}
\subsection{Notation} When we write $x=y\pm \alpha$, we mean $x\in[y-\alpha,y+\alpha]$. When we write $y\pm \alpha=z\pm\beta$, we mean $[y-\alpha,y+\alpha]\subset [z-\beta,z+\beta]$. Note that the latter convention is not symmetric, that is, $y\pm \alpha=z\pm\beta$ is not the same as $z\pm\beta=y\pm \alpha$.

The neighbourhood of a vertex $v$ in the graph $G$ is denoted $\NBH_{G}(v)$.
We write $\NBH_{G}(U)=\bigcap_{v\in U}\NBH_{G}(v)$ for the common neighbourhood of the set $U\subset V(G)$. 

The definition of degenerate graphs naturally suggests to label the vertices of a graph by integers. Suppose that the vertices of a graph $G$ are $V(G)=[\ell]$. Suppose that $i\in V(G)$. We write $\LNBH(i)=\NBH(i)\cap[i-1]$ and $\LEFTDEG(i)=|\LNBH(i)|$ for the \emph{left-neighbourhood} and the \emph{left-degree} of $i$. We make use of the natural order on $[\ell]$ also in other ways, like referring to sets of the form $[\ell_1]\subset V(G)$ and $\{\ell_2,\ell_2+1,\ldots,\ell\}\subset V(G)$ as \emph{initial vertices} and \emph{final vertices}, respectively.
The \emph{density} of a graph $H$ is the quantity ${e(H)}/{\binom{v(H)}{2}}$.

The graphs to be packed in Theorem~\ref{thm:MAINunbounded} are denoted $G_t$ because they are \textbf{g}uest graphs. By contrast, during our packing procedure, we shall work with \textbf{h}ost graphs $H_s$ which are obtained from the original $K_n$ by removing what was used previously.

\subsection{Probability}
\subsubsection{Probability basics}
All probability spaces considered in this paper are finite. The implicit
sigma-algebra underlying each such space is the sigma-algebra generated by all
singletons; in particular, the notion of measurability is trivial in this
setting. Recall that if $\Omega$ is finite probability space then a sequence of
partitions $\mathcal{F}_0$, $\mathcal{F}_1$,\dots,  $\mathcal{F}_n$ of~$\Omega$ is a \emph{filtration} if each partition $\mathcal{F}_i$ refines its predecessor $\mathcal{F}_{i-1}$.\footnote{Readers familiar with measure-theoretic probability will notice that the standard definition is a sequence of $\sigma$-algebras, namely those generated by our partitions; in the finite setting this is an unnecessary complication.} In this setting, a function $f:\Omega\rightarrow \mathbb{R}$ is called \emph{$\mathcal{F}_i$-measurable} if $f$ is constant on each part of $\mathcal{F}_i$.

Recall also that if $\Omega$ is a finite probability space and $f:\Omega\rightarrow\mathbb{R}$ is a function, then the \emph{conditional expectation} $\EXP(f|\mathcal{F}):\Omega\rightarrow\mathbb{R}$ and the \emph{conditional variance} $\Var(f|\mathcal{F}):\Omega\rightarrow\mathbb{R}$ of $f$ with respect to a given partition $\mathcal{F}$ of $\Omega$ are defined by
\begin{equation*}
\begin{aligned}
\EXP(f|\mathcal{F})(x)&=\EXP(f|X),\\
\Var(f|\mathcal{F})(x)&=\Var(f|X),
\end{aligned}
\qquad
\text{where $X\in\mathcal{F}$ is such that $X\ni x$.}
\end{equation*}

\subsubsection{Sequential dependence and concentration}\label{ssec:seqdependence}
In this section we introduce some convenient consequences of standard martingale inequalities. These are generally useful in the analysis of randomised processes, so we try to provide some brief background and motivation.

Suppose that we have a randomised algorithm which proceeds in $m$ rounds. We can
then denote by $\Omega:=\prod_{i=1}^{m}\Omega_i$ the probability space that
underlies an execution of the algorithm. Here $\Omega_i$ is the set of all
possible choices the algorithm may make in step $i$. It is important, however, that $\Omega$ as a probability space is not necessarily a product of probability spaces $\Omega_i$; in other words, the algorithms can (and typically will) make choices for the step $i$ depending on the choices it made in steps $1,\ldots,i-1$. By \emph{history up to time $t$} we mean a set of the form $\{\omega_1\}\times\cdots\times\{\omega_t\}\times \Omega_{t+1}\times\cdots\Omega_m$, where $\omega_i\in\Omega_i$. We shall use the symbol $\hist_{t}$ to denote any particular history of such a form. By a \emph{history ensemble up to time $t$} we mean any union of histories up to time $t$; we shall use the symbol $\histens$ to denote any one such. Observe that there are natural filtrations associated to such a probability space: given times $t_1<t_2<\dots$ we let $\mathcal{F}_{t_i}$ denote the partition of $\Omega$ into the histories up to time $t_i$. We introduce formally a probability space of this type, which we use for the key part of our argument, in Section~\ref{ssec:AlgorithmProbSpace}.

\medskip

We recall that if $Y_1,\dots,Y_n$ are a collection of independent random variables, whose ranges are not too large compared to $n$, 
%then $\sum_{i=1}^nY_i$ is close in distribution to a Gaussian random variable with mean $\mu=\sum_{i=1}^n\EXP(Y_i)$ and variance $\sum_{i=1}^n\Var(Y_i)$. In particular, 
we have Hoeffding's inequality for the tails of such sums:
\begin{equation}\label{eq:hoeffding}
 \PROB\Big(\sum_{i=1}^n\big(Y_i-\EXP(Y_i)\big)\ge \rho\Big)\le\exp\Big(-\frac{2\rho^2}{\sum_{i=1}^n(\max Y_i-\min Y_i)^2}\Big)\,,
\end{equation}
for each $\rho>0$. One should think of the squared range of $Y_i$ as a crude upper bound for $\Var(Y_i)$. There are various improvements, such as the Bernstein inequalities, which take into account the actual values $\Var(Y_i)$ in order to obtain stronger concentration results such as
\begin{equation}
\label{eq:bernstein}
 \PROB\Big(\sum_{i=1}^n\big(Y_i-\EXP(Y_i)\big)\ge \rho\Big)\le\exp\Big(-\frac{\rho^2}{2R\rho/3+2\sum_{i=1}^n\Var(Y_i)}\Big)\,,
\end{equation}
valid when $0\le Y_i\le R$ for each $i$. When the sum of variances is much larger than $R\rho$, this probability bound is optimal up to small order terms in the exponent; for most applications this means it cannot usefully be improved.

However when analysing randomised algorithms, usually one has to deal with a sum of random variables which are not independent, but rather are \emph{sequentially dependent}, meaning that they come in an order in which earlier outcomes affect the later random variables. A good example is the following procedure (a variant of which we use in this paper) for embedding a graph $G$ on vertex set $[n/2]$ into a graph $H$ on $n$ vertices. We simply embed vertices in order $1,\dots,n/2$, at each time $t$ embedding vertex $t$ uniformly at random to the set of all valid choices: that is, choices which give an embedding of $G[1,\dots,t]$. In order to show that this procedure is likely to succeed (which is true if~$G$ has small degeneracy and $H$ is sufficiently quasirandom) we will want to know how vertices are embedded over time to some subsets $S\subset V(H)$. In other words, we define (in this case, Bernoulli) random variables $Y_t$ to be~1 if $t$ is embedded to $S$ and $0$ otherwise, and we want to know how the partial sums of these random variables, which are certainly not independent but are sequentially dependent, behave. The point of this section is to observe that in fact more or less the same concentration bounds hold as for independent random variables, except that one has to replace the sum of expectations with a sum of observed expectations, that is, $\sum_{i=1}^n\EXP\big(Y_i|\hist_{i-1}\big)$, where $\hist_{i-1}$ denotes the history up to time $i-1$, and the sum of variances with a sum of observed variances, similarly defined.

In combinatorial applications, one is usually interested in showing that a sum of random variables (which might in general not be Bernoulli) is close to its expectation~$\mu$. It is not \emph{a priori} obvious that concentration bounds such as the above help: after all, the sum of observed expectations is itself a random variable and might not be concentrated near $\mu$ (it is easy to come up with examples in which it is not). We deal with this in what follows by defining a \emph{good event} $\cE$, within which the observed sum of expectations is $\mu\pm\nu$ for some (small) $\nu>0$. In applications $\cE$ will often be a combinatorial statement about the process, and hence we refer to $\nu$ as the \emph{combinatorial error}, to distinguish it from the \emph{probabilistic error} $\rho>0$, as in~\eqref{eq:hoeffding} and~\eqref{eq:bernstein}. It is important to note that $\cE$ is usually not determined before the random variables $Y_i$ (i.e.\ it may well not be $\mathcal{F}_i$-measurable for any member~$\mathcal{F}_i$ of the filtration), so we do not condition on $\cE$, rather we aim to estimate the probability that $\cE$ holds and yet $\sum_{i=1}^nY_i\neq\mu\pm(\nu+\rho)$.

In order to avoid mentioning any particular process, it is convenient to state
the following lemmas in terms of a finite probability space $\Omega$ with a
filtration $(\mathcal{F}_0,\mathcal{F}_1,\dots,\mathcal{F}_n)$. We should stress
that though in our applications we will always use the same probability space, which underlies our
packing process, we will consider different filtrations, always given by the histories up to increasing times, depending on the random variables we wish to sum. 

The following lemma, from~\cite{AAGH:AlmostAllTrees}, is a sequential dependence
version of Hoeffding's inequality. Note that the lemma as stated
in~\cite{AAGH:AlmostAllTrees} includes the condition $\PROB(\cE)>0$. However if
$\PROB(\cE)=0$ the lemma statement is trivially true, so we drop the condition
below.

\begin{lem}[Lemma~7, \cite{AAGH:AlmostAllTrees}]
\label{lem:concentration}Let $\Omega$ be a finite probability space,
and $(\mathcal{F}_0, \mathcal{F}_{1},\dots,\mathcal{F}_{n})$ be filtration. Suppose
that for each $1\le i\le n$ we have a nonnegative real number $a_{i}$,
an $\mathcal{F}_i$-measurable random variable $Y_{i}$ satisfying
$0\le Y_{i}\le a_{i}$, nonnegative real numbers~$\mu$ and~$\nu$,
and an event $\cE$. Suppose that almost surely, either $\cE$ does
not occur or $\sum_{i=1}^{n}\EXP\left(Y_{i}\big|\mathcal{F}_{i-1}\right)=\mu\pm\nu$.
Then for each $\varrho>0$ we have 
\[
\PROB\left(\mathcal{\cE}\text{ and }\sum_{i=1}^{n}Y_{i}\neq\mu\pm(\nu+\varrho)\right)\le2\exp\Big(-\frac{2\varrho^{2}}{\sum_{i=1}^{n}a_{i}^{2}}\Big)\,.
\]
Furthermore, if we weaken the assumption, requiring only that either $\cE$ does not occur or $\sum_{i=1}^n\EXP\left(Y_i\big|\mathcal{F}_{i-1}\right)\le\mu+\nu$, then for each $\rho>0$ we have
\[
\PROB\left(\mathcal{\cE}\text{ and }\sum_{i=1}^{n}Y_{i}>\mu+\nu+\varrho\right)\le\exp\Big(-\frac{2\varrho^{2}}{\sum_{i=1}^{n}a_{i}^{2}}\Big)\,.
\]
\end{lem}

We should note that the probability bound in this lemma is what one would
obtain from standard martingale inequalities for
$\PROB(\sum_{i=1}^{n}Y_{i}\neq\mu\pm(\nu+\varrho))$ if the condition $\sum_{i=1}^{n}\EXP\left(Y_{i}\big|\mathcal{F}_{i-1}\right)=\mu\pm\nu$ held almost surely. The r\^ole of $\cE$ is that we can allow this condition to fail outside of~$\cE$ but still obtain the same concentration within $\cE$; this is probabilistically fairly trivial but very useful. The same applies for the next lemma.

Lemma~\ref{lem:concentration} gives close to optimal (up to a constant factor in the exponential) results when the random variables $Y_i$ are relatively often close to $0$ and $a_i$; in other words, when $a_i^2$ is not much larger than the variance $\Var(Y_i)$. This will turn out to be the case for most of the random sums we need to estimate in this paper. However, when it is not the case, at the cost of a second moment calculation the following version of Freedman's inequality~\cite{Freedman} gives much stronger bounds, corresponding to a Bernstein inequality for independent random variables.

\begin{lem}[Freedman's inequality on a good event]
\label{lem:freedman}Let $\Omega$ be a finite probability space,
and $(\mathcal{F}_0, \mathcal{F}_{1},\dots,\mathcal{F}_{n})$ be a filtration. Suppose
that we have $R>0$, and for each $1\le i\le n$ we have an $\mathcal{F}_i$-measurable non-negative random variable $Y_{i}$, nonnegative real numbers~$\mu$,~$\nu$ and $\sigma$,
and an event $\cE$. Suppose that almost surely, either $\cE$ does
not occur or we have $\sum_{i=1}^{n}\EXP\left(Y_{i}\big|\mathcal{F}_{i-1}\right)=\mu\pm\nu$, and $\sum_{i=1}^n\Var\left(Y_i\big|\mathcal{F}_{i-1}\right)\le\sigma^2$, and $0\le Y_i\le R$ for each $1\le i\le n$.
Then for each $\varrho>0$ we have 
\[
\PROB\left(\mathcal{\cE}\text{ and }\sum_{i=1}^{n}Y_{i}\neq\mu\pm(\nu+\varrho)\right)\le2\exp\Big(-\frac{\varrho^{2}}{2\sigma^2+2R\rho}\Big)\,.
\]
Furthermore, if we assume only that either $\cE$ does
not occur or we have $\sum_{i=1}^{n}\EXP\left(Y_{i}\big|\mathcal{F}_{i-1}\right)\le\mu+\nu$, and $\sum_{i=1}^n\Var\left(Y_i\big|\mathcal{F}_{i-1}\right)\le\sigma^2$, and $0\le Y_i\le R$ for each $1\le i\le n$,
then for each $\varrho>0$ we have 
\[
\PROB\left(\mathcal{\cE}\text{ and }\sum_{i=1}^{n}Y_{i}>\mu+\nu+\varrho\right)\le\exp\Big(-\frac{\varrho^{2}}{2\sigma^2+2R\rho}\Big)\,.
\]
\end{lem}

  As with the Bernstein inequality, this result is essentially optimal when the sum of observed variances is much larger than $R\rho$. We would like to point out that since $\cE$ is often a combinatorial statement which is not tailored to the specific random variables $Y_i$ we are summing, when we use either lemma to estimate tail probabilities for several sums of random variables, we will often use the same event $\cE$ repeatedly; since it will appear only once in union bounds, both lemmas are useful for showing that a.a.s.\ a collection of many (rapidly growing with $n$) sums are simultaneously close to their expectations, even when the probability of $\cE$ only tends to one quite slowly with $n$.
  
  We deduce Lemma~\ref{lem:freedman} from Freedman's martingale inequality, which we now state.

\begin{thm}[Proposition~(2.1),~\cite{Freedman}]\label{thm:freedman}
Let $\Omega$ be a finite probability space,
and $(\mathcal{F}_0, \mathcal{F}_{1},\dots,\mathcal{F}_{n})$ be a filtration. Suppose
that for some $R>0$, for each $1\le i\le n$,
we have an $\mathcal{F}_i$-measurable random variable $Y_{i}$ that takes values in the range $-R\le Y_i\le R$, and we have $\EXP(Y_i|\mathcal{F}_{i-1})=0$ almost surely. Suppose that for some $\sigma$ we have $\sigma^2\ge\sum_{i=1}^n\Var(Y_i|\mathcal{F}_{i-1})$ almost surely. Then for each $\rho>0$, we have
\[\PROB\left(\sum_{i=1}^nY_i\ge\rho\right)\le \exp\Big(-\frac{\rho^2}{2\sigma^2+2R\rho}\Big)\,.\]
\end{thm}

We now deduce Lemma~\ref{lem:freedman}, using a similar approach as was used in~\cite{AAGH:AlmostAllTrees} to prove Lemma~\ref{lem:concentration}. 
\begin{proof}[Proof of Lemma~\ref{lem:freedman}]
 We show the required upper bound
 \begin{equation}\label{eq:freedman:up}
  \PROB\left(\mathcal{\cE}\text{ and }\sum_{i=1}^{n}Y_{i}>\mu+\nu+\varrho\right)\le\exp\Big(-\frac{\varrho^{2}}{2\sigma^2+2R\rho}\Big)\,,
 \end{equation}
 and the corresponding lower bound follows by symmetry, replacing each $Y_i$ with $R-Y_i$. This gives the desired two-sided result by the union bound.
 
 Observe that if $\PROB(\cE)=0$,~\eqref{eq:freedman:up} holds trivially. We may thus assume $\PROB(\cE)>0$. Now, given $Y_1,\dots,Y_n$, we define random variables $U_1,\dots,U_n$ as follows. We set $U_i=\max(Y_i,R)$ if $\PROB(\cE|\mathcal{F}_{i-1})>0$, and otherwise $U_i=0$. Observe that $U_i$ is constant on each part of $\mathcal{F}_{i}$ by definition. We claim that for each $1\le t\le n$ we have almost surely
 \begin{equation}\label{eq:freedman:Ui}
  \sum_{i=1}^t\EXP(U_i|\mathcal{F}_{i-1})\le\mu+\nu\quad\text{and}\quad\sum_{i=1}^t\Var(U_i|\mathcal{F}_{i-1})\le\sigma^2\,.
 \end{equation}
 Indeed, suppose that $t$ is minimal such that this statement fails, and let $F$ be a set in $\mathcal{F}_{t-1}$ with $\PROB(F)>0$ witnessing its failure. By minimality of $t$, at least one of $\EXP(U_t|F)$ and $\Var(U_t|F)$ is strictly positive. By definition of $U_t$ we have $\PROB(\cE|F)>0$. But since $\EXP(U_i|\mathcal{F}_{i-1})$ and $\Var(U_i|\mathcal{F}_{i-1})$ are nonnegative for each $i$, this shows that with probability at least $\PROB(F)\PROB(\cE|F)>0$, the event $\cE$ occurs and one of the assumptions $\sum_{i=1}^{n}\EXP\left(Y_{i}\big|\mathcal{F}_{i-1}\right)=\mu\pm\nu$ and $\sum_{i=1}^n\Var\left(Y_i\big|\mathcal{F}_{i-1}\right)\le\sigma^2$ fails. This is a contradiction, so we conclude~\eqref{eq:freedman:Ui} holds almost surely for each $t$. Furthermore, we have $0\le U_i\le R$ for each $1\le i\le n$.
 
 Next, define for each $1\le i\le n$ the random variable $W_i=U_i-\EXP(U_i|\mathcal{F}_{i-1})$. We have $-R\le W_i\le R$ for each $i$, by definition $W_i$ is $\mathcal{F}_i$-measurable, and by definition almost surely $\EXP(W_i|\mathcal{F}_{i-1})=0$ and $\Var(W_i|\mathcal{F}_{i-1})=\Var(U_i|\mathcal{F}_{i-1})$. Thus by Theorem~\ref{thm:freedman} we have
 \[\PROB\left(\sum_{i=1}^nW_i\ge\rho\right)\le \exp\Big(-\frac{\rho^2}{2\sigma^2+2R\rho}\Big)\,.\]
 Since almost surely we have $\sum_{i=1}^t\EXP(U_i|\mathcal{F}_{i-1})\le\mu+\nu$, we obtain
 \[\PROB\left(\sum_{i=1}^nU_i\ge\mu+\nu+\rho\right)\le \exp\Big(-\frac{\rho^2}{2\sigma^2+2R\rho}\Big)\,.\]
 Finally, if $\cE$ occurs then almost surely $Y_i=U_i$ for each $1\le i\le n$, giving the desired upper bound~\eqref{eq:freedman:up}.
\end{proof}

\medskip

Finally, let us note that we shall be using many statements of the form
\begin{equation}\label{eq:nocondwords}
\mbox{with probability at least $p$, provided event $\mathcal A$ we get event $\mathcal B$.}
\end{equation}
We emphasize that such statements are not statements about conditional probabilities. That is, the meaning of~\eqref{eq:nocondwords} is $\PROB(\mathcal A\setminus \mathcal B)\le 1-p$. A prototypical example is \emph{with probability at least $1-o(1)$, if a given randomized algorithm does not fail, then it produces an output with certain desired properties}.

\subsection{Simple properties of degenerate graphs}

We need to bound $\sum_{x\in V(G)}\deg(x)^{2}$ for
degenerate graphs $G$. In several applications of Lemma \ref{lem:concentration}
the numbers $a_{i}$ will be upper bounded by the degrees of vertices
in $G$, where $G$ is one of the graphs to be packed, so that $\sum_{x\in V(G)}\deg(x)^{2}$ is an upper bound
for the sum $\sum_{i}a_{i}^{2}$ appearing in Lemma \ref{lem:concentration}.
\begin{lem}
\label{lem:squarebound}Let $G$ be an $n$-vertex graph with degeneracy
$D$ and maximum degree $\Delta$. Then we have
\[
\sum_{x\in V(G)}\deg(x)^{2}\le2Dn\Delta\,.
\]
\end{lem}
\begin{proof}
We have \[
\sum_{x\in V(G)}\deg(x)^{2}\le \sum_{x\in V(G)}\deg(x)\cdot \Delta=2 e(G)\cdot \Delta \le 2 D n\cdot \Delta\,.
\]
\end{proof}

We also need to show that degenerate graphs contain large independent sets all of whose vertices have the same degree.
\begin{lem}\label{lem:degindept}
 Let $G$ be a $D$-degenerate $n$-vertex graph. Then there exists an integer $0\le d\le 2D$ and a set $I\subseteq V(G)$ with $|I|\ge (2D+1)^{-3}n$ which is independent, and all of whose vertices have the same degree $d$ in $G$.
\end{lem}
\begin{proof}
 We first claim that at least $(2D+1)^{-1}n$ vertices of $G$ have degree at most~$2D$. Indeed, if this were false then there would be more than $2Dn/(2D+1)$ vertices of $G$ all of whose degrees are at least $2D+1$, so that we obtain $e(G)>Dn$, which contradicts the $D$-degeneracy of $G$. Let $0\le d\le 2D$ be chosen to maximise the number of vertices in $G$ of degree $d$, and let $S$ be the set of vertices in $G$ with degree $d$. We thus have $|S|\ge (2D+1)^{-2}n$. Now let $I$ be a maximal independent subset of $S$. Each vertex of $I$ has at most $d\le 2D$ neighbours in $S$, so that $\big|I\cup \bigcup_{i\in I}\NBH(i)\big|\le (2D+1)|I|$. By maximality $I\cup \bigcup_{i\in I}\NBH(i)$ covers $S$, hence $|I|\ge(2D+1)^{-1}|S|\ge(2D+1)^{-3}n$, as desired.
\end{proof}

\section{Reducing the main theorem}
\label{sec:reduction}

We deduce Theorem \ref{thm:MAINunbounded} from the following technical
result.
\begin{thm}
\label{thm:maintech}For each $\gamma>0$ and each $D\in\mathbb{N}$
there exists $c>0$ and a number $n_{0}$ such that the following holds for each
integer $n>n_{0}$. Suppose that $s^{*}\le 2n$ and that for each $s\in [s^{*}]$
the graph $G_{s}$ is a graph on vertex set $[n]$, with maximum degree
at most $\tfrac{cn}{\log n}$, such that $\LEFTDEG(x)\le D$ for each $x\in V(G_{s})$
and such that the last $(D+1)^{-3}n$ vertices of $[n]$ form an independent
set in $G_{s}$, and all have the same degree $d_s$ in $G_s$. Suppose further that the total number of edges of
$\left(G_{s}\right)_{s\in [s^{*}]}$ is at most $(1-3\gamma)\binom{n}{2}$.
Then $\left(G_{s}\right)_{s\in [s^{*}]}$ packs into $K_{n}$. 
\end{thm}

Actually, we prove Theorem~\ref{thm:maintech} in a slightly more general form using the concept of quasirandomness which is crucial for our approach. This concept was introduced by several authors independently in the 1980s (of which the paper~\cite{Chung1989} is the most comprehensive) and captures a property that the edges of graph are distributed evenly among its vertices.  We give a definition tailored for our needs which is somewhat stronger than the usual definition of quasirandom graphs.
  \begin{defn}[quasirandom]
  	\label{def:quasirandomgraph}Suppose that $H$ is a graph with $n$ vertices and with density~$p$. We say that such graph $H$
  	 is \emph{$(\alpha,L)$-quasirandom} if for
  	every set $S\subset V(H)$ of at most $L$ vertices we have
  	$|\NBH_{H}(S)|=(1\pm\alpha)p^{|S|}n$.
  \end{defn}
\ENHANCED{We can now give our full main technical result. We note that the `furthermore' part of Theorem~\ref{thm:maintechGENERAL}, was added in~2018, is not needed for the proof of Theorem~\ref{thm:MAINunbounded}. However, establishing these additional properties will be useful in ????. This furthermore part uses the so-called diet condition, which we define in Definition~\ref{def:dietcondition} below.}
\begin{thm}[Main technical result]
	\label{thm:maintechGENERAL}For each $\gamma>0$ and each $D\ENHANCED{,L}\in\mathbb{N}$
	there exist numbers $n_{0}\in\mathbb{N}$ and $c,\xi>0$ such that the following holds for each
	$n>n_{0}$. Suppose that~$\widehat H$ is an $(\xi,2D+3)$-quasirandom graph with $n$ vertices and density $p>0$. \ENHANCED{Suppose that we are given sets $U_1,\ldots,U_L\subset V(\widehat H)$.} Suppose that $s^{*}\le 2n$ and that for each $s\in [s^{*}]$
	the graph $G_{s}$ is a graph on vertex set $[n]$, with maximum degree
	at most $\tfrac{cn}{\log n}$, such that $\LEFTDEG(x)\le D$ for each $x\in V(G_{s})$
	and such that the last $(D+1)^{-3}n$ vertices of $[n]$ form an independent
	set in $G_{s}$, and all have the same degree $d_s$ in $G_s$. Suppose further that the total number of edges of
	$\left(G_{s}\right)_{s\in [s^{*}]}$ is at most $(p-3\gamma)\binom{n}{2}$.
	Then $\left(G_{s}\right)_{s\in [s^{*}]}$ packs into~$\widehat H$. 
	
	\ENHANCED{Furthermore, let $H^\dagger$ by the graph formed by the edges not used in the packing. Then $H^\dagger$ is $(\gamma,2D+3)$-quasirandom. For $s\in [s^*]$, let $\mathrm{img}(G_s)\subset V(\widehat H)$ be the image of $G_s$ under the said packing. Then for each $s\in[s^*]$ and each $\ell\in[L]$, the pair $(H^\dagger,\mathrm{img}(G_s))$ satisfies the $(\gamma,2D+3)$-diet condition relative to the set $U_\ell$.}
\end{thm}
Theorem~\ref{thm:maintechGENERAL} indeed generalizes Theorem~\ref{thm:maintech} because it can be easily checked that for any fixed $D\in\mathbb{N}$ and $\alpha>0$, the graph $K_n$ is $(\alpha,2D+3)$-quasirandom for $n$ sufficiently large. The reason why we give the proof in this greater generality is that it is clear that the only feature of $K_n$ we actually use is its quasirandomness. We show that Theorem~\ref{thm:maintech} implies Theorem~\ref{thm:MAINunbounded}. Note that starting with Theorem~\ref{thm:maintechGENERAL} the same deduction would yield a version of Theorem~\ref{thm:MAINunbounded} for quasirandom host graphs. We state such a version for dense Erd\H{o}s--R\'enyi random graphs $\mathbb G(n,p)$, an $n$-vertex graph, where each pair of vertices forms an edge independently with probability~$p$. Those graphs are well-known to have asymptotically almost surely error in quasirandomness (even in our Definition~\ref{def:quasirandomgraph}) tending to zero.
\begin{thm}\label{thm:packingIntoErdosRenyi}
	For each $p,\gamma>0$ and each $D\in\mathbb{N}$
	there exists $c>0$ such that the following holds asymptotically almost surely, as $n\rightarrow\infty$. Suppose that $\left(G_{t}\right)_{t\in [t^{*}]}$
	is a family of $D$-degenerate graphs, each of which has at most $n$
	vertices and maximum degree at most $\tfrac{cn}{\log n}$. Suppose further that the
	total number of edges of $\left(G_{t}\right)_{t\in [t^{*}]}$ is at
	most $(p-\gamma)\binom{n}{2}$. Then $\left(G_{t}\right)_{t\in [t^{*}]}$
	packs into $\mathbb G(n,p)$. 
\end{thm}

\begin{proof}[Proof of Theorem~\ref{thm:MAINunbounded}]
  To deduce Theorem \ref{thm:MAINunbounded} from
  Theorem~\ref{thm:maintech}, observe that given an integer $D$ and
  graphs~$\mathcal{G}=(G_{t})_{t\in[t^*]}$ to pack, we may assume without loss of generality that
  none of the graphs in $\mathcal{G}$ has isolated vertices, since such vertices can be erased and then easily packed in the last step.
  
  We now successively modify the family $\mathcal{G}$ as follows. If there are two graphs $G,G'\in\mathcal{G}$ with $v(G),v(G')\le n/2$, we replace $G$ and $G'$ with the disjoint union $G\cup G'$. We repeat this until no further such pairs exist, giving $\mathcal{G}'$.
  
  Observe that the maximum degree and the degeneracy of the graphs in $\mathcal{G}$ is the same as in $\mathcal{G}'$. Furthermore a packing of $\mathcal{G}'$ is also a packing of $\mathcal{G}$. Finally, there is at most one graph in $\mathcal{G}'$ with less than $n/2$ vertices. Hence all but at most one graph has at least $n/4$ edges. We conclude that the total number $s^*$ of graphs in $\mathcal{G}'$ satisfies $(s^*-1)n/4\le (1-\gamma)\binom{n}{2}$, and hence $s^*\le 2n$. Finally, we let the graphs $(G'_s)_{s=1}^{s^*}$ be obtained from the graphs $\mathcal{G}'$ by adding if necessary isolated vertices to each in order to obtain $n$-vertex graphs.
  
  Now, for each
  $G'_{s}$ we choose an order on~$V(G'_{s})$ as follows. First, we pick an
  order witnessing $D$-degeneracy of $G'_{s}$. Next, we pick an integer $0\le d_s\le 2D$ and an independent $I_s$
  set of $(2D+1)^{-3}n$ vertices each of which has degree $d_s$ in
  $G'_{s}$ and change the order by moving these vertices to the
  end. Such an integer $d_s$ and independent set exist by Lemma~\ref{lem:degindept}. The result is an ordering of $V(G'_s)$ with degeneracy at most $2D$, as required for Theorem \ref{thm:maintech} with input $2D$ and $\gamma/3$. Then Theorem~\ref{thm:maintech} returns the desired packing.
\end{proof}

\section{Proof of Theorem~\ref{thm:maintechGENERAL}}
\label{sec:prooftech}

For the proof of Theorem~\ref{thm:maintechGENERAL}, we need some algorithms and definitions. We give these now along with a sketch of the proof.

We prove Theorem~\ref{thm:maintechGENERAL}
by analysing a randomised algorithm, which we call \PackingProcess{}, that packs the
guest graphs $G_s$ into~$\widehat H$. We prove that this algorithm succeeds with
high probability. In this algorithm we assume that the last~$\delta n$ vertices of each
graph~$G_s$ form an independent set, where $\delta<(D+1)^{-3}$  is to be chosen later.

\PackingProcess{} begins by splitting the edges of the input graph $\widehat H$ into a \emph{bulk} $H_0$ and a \emph{reservoir} $H^*_0$ by independently selecting edges into the latter with probability chosen such that $e(H^*_0)\approx\gamma\binom{n}{2}$. As a result, the graphs $H_0$ and $H_0^*$ are with high probability quasirandom.

Now \PackingProcess{} proceeds in~$s^*$ stages. In each stage~$s$, it runs a randomised embedding algorithm, called
\RandomEmbedding{} and explained below, to embed the first $n-\delta n$ vertices of~$G_s$ into the bulk~$H_{s-1}$. Then in the \emph{completion phase} the last $\delta n$ vertices of $G_s$ are embedded into the reservoir $H^*_{s-1}$. Since there are exactly~$\delta n$ vertices of $G_s$ left to embed and exactly~$\delta n$ vertices of $V(\widehat{H})$ unused so far in this stage, we want to find a bijection between these. Since all neighbours of each yet unembedded vertex are already embedded, this completion amounts to choosing a system of distinct representatives. The completion phase does not use randomness: the system of disjoint representatives is obtained using Hall's theorem. Now $H_s$ and $H^*_s$ are defined simply by removing the edges used in this embedding.

Both \RandomEmbedding{} and the completion phase may \emph{fail} at any stage $s$; this means that it is not possible to embed a certain part of $G_s$. In that case \PackingProcess{} fails, too. If \PackingProcess{} does not fail then it always produces a valid packing of $(G_s)$ into $H$. So, we need to show that \PackingProcess{} (see Algorithm~\ref{alg:pack}) succeeds with positive probability.

\begin{algorithm}[ht]
    \caption{\PackingProcess{}}\label{alg:pack}
    \SetKwInOut{Input}{Input}
    \Input{graphs $G_1,\dots,G_{s^*}$, with $G_s$ on vertex set $[n]$
      such that  the last $\delta n$ vertices of~$G_s$ form an independent set; a graph $\widehat H$ on $n$ vertices}
    choose $H^*_0$ by picking edges of $\widehat H$ independently with probability $\gamma\binom{n}{2}/e(\widehat{H})$ \;
    let $H_0=\widehat H-H^*_0$ \;
    \For{$s=1$ \KwTo $s^*$}{
      run \RandomEmbedding($G_s$,$H_{s-1}$) to get
      an embedding $\phi_s$ of~$G_s[{\scriptstyle [n-\delta n]}]$ into~$H_{s-1}$\;   
      let $H_{s}$ be the graph obtained from $H_{s-1}$ by removing the
      edges of $\phi_{s}\big(G_s[{\scriptstyle [n-\delta n]}]\big)$\;
      choose an extension $\phi^*_s$ of $\phi_s$ embedding all of $G_s$ and embedding the edges of $G_s-G_s[{\scriptstyle [n-\delta n]}]$ into $H^*_{s-1}$ \;
      let $H^*_s$ be the graph obtained from $H^*_{s-1}$ by removing the edges of $\phi^*_s\big(G_s-G_s[{\scriptstyle [n-\delta n]}]\big)$ \;
    }    
  \end{algorithm}

\medskip

For describing our randomised embedding algorithm
\RandomEmbedding{} we need the following definitions. We shall use the symbol $\AlgMap$ to denote embeddings produced by \RandomEmbedding{}. We write $G\AlgMap H$ to indicate that the graph $G$ is to be embedded into $H$. Also, if $t\in V(G)$, $v\in V(H)$ and $A\subset V(H)$ then $t\AlgMap v$ means that $t$ is embedded on $v$, and $t\AlgMap A$ means that $t$ is embedded on a vertex of $A$. 

\begin{defn}[partial embedding, candidate set]
  Let~$G$ be a graph with vertex set $[v(G)]$, and~$H$ be a graph with
  $v(H)\ge v(G)$.
  Further, assume $\psi_{j}\colon[j]\rightarrow V(H)$ is a \emph{partial embedding}
  of $G$ into~$H$ for $j\in[v(G)]$, that is, $\psi_j$ is a graph embedding
  of $G\big[[j]\big]$ into~$H$. Finally, let $t\in[v(G)]$ be such that $\LNBH_G(t)\subset[j]$.
  Then the \emph{candidate set of $t$} (with respect to~$\psi_j$) is 
  \[\CANDSET_{G\AlgMap H}^{j}(t)=\NBH_{H}\Big(\psi_{j}\big(\LNBH_{G}(t)\big)\Big)\,.\]
  When $j=t-1$, we call $\CANDSET_{G\AlgMap H}^{j}(t)$ the \emph{final candidate set of $t$}.
\end{defn}

\RandomEmbedding{} (see Algorithm~\ref{alg:embed}) randomly embeds a guest graph~$G$ into a host
graph~$H$. The algorithm is simple: we
iteratively embed the first $(1-\delta)n$ vertices of~$G$ randomly to one
of the vertices of their candidate set which was not used for embedding
another vertex already.

  \begin{algorithm}[ht]
    \caption{\RandomEmbedding{}}\label{alg:embed}
    \SetKwInOut{Input}{Input}
    \Input{graphs~$G$ and~$H$, with $V(G)=[v(G)]$ and $v(H)=n$}
    $\psi_0:=\emptyset$\;
    $t^*:=(1-\delta)n$\;
    \For{$t=1$ \KwTo $t^*$}{
      \lIf{$\CANDSET_{G\AlgMap H}^{t-1}(t)\setminus\im(\psi_{t-1})=\emptyset$}{
        halt with failure}
      choose $v\in\CANDSET_{G\AlgMap
        H}^{t-1}(t)\setminus\im(\psi_{t-1})$ uniformly at random\;
      $\psi_{t}:=\psi_{t-1}\cup\{t\AlgMap v\}$\;
    }
    \Return $\psi_{t^*}$
  \end{algorithm}

  To show that \PackingProcess{} does not fail at any stage, we shall show
  that the host graph $H_{s}$ constructed in \PackingProcess{} in embedding
  stage~$s$ is quasirandom in the sense of Definition~\ref{def:quasirandomgraph}.
  In fact, in order to analyse the completion phase of \PackingProcess{} we need quasirandomness of the pair $(H_s,H^*_0)$, where $H^*_0$ is the initial reservoir. We now define this \emph{coquasirandomness} of a pair of graphs. Recall that quasirandomness of one graph means that common neighbourhoods are always about the size one would expect in a random graph of a similar density. Coquasirandomness of two graphs means that the intersection of a common neighbourhood in the first graph and another in the second graph has about the size one would expect in two independent random graphs of the respective densities.
  
  \begin{defn}[coquasirandom]
    \label{def:conj-quasirandom}For $\alpha>0$ and $L\in\mathbb N$, we say that a pair of graphs $(F, F^*)$, both on the same vertex set~$V$ of order~$n$
    and with densities $p$ and $p^*$, respectively, is \emph{$(\alpha,L)$-coquasirandom} if for
    every set $S\subset V$ of at most $L$ vertices and every subset $R\subseteq S$ we have
    \[|\NBH_{F}(R)\cap \NBH_{F^*}(S\setminus R)|
    =(1\pm\alpha)p^{|R|}(p^*)^{|S\setminus R|}n\,.\]
  \end{defn}
  
  With this we can state the setting of our main lemmas and fix various constants which we will use in the remainder of the paper.
  
  \begin{setting}\label{set:graphs}
Let $D,n\in\mathbb N$ and $\gamma>0$ be given. We define
\begin{equation}\label{eq:defconsts}\begin{split}
 \eta&=\frac{\gamma^D}{200D}\,,\quad  \delta=\frac{\gamma^{10D}\eta}{10^{6}D^{4}}\,,\quad C=40D\exp\big(1000D\delta^{-2}\gamma^{-2D-10}\big)\,,\\
  \alpha_x&=\frac{\delta}{10^8CD}\exp\Big(\frac{10^8C^2D^3\delta^{-1}\gamma^{-4D-6}(x-2n)}{n}\Big) \qquad\text{for each $x\in\mathbb{R}$},\\
  \eps&= \alpha_0\delta^{\corrig{4}}\gamma^{10D}/1000CD\,,\quad c= D^{-4}\eps^{4}/100\,\quad\text{and}\quad \xi=\alpha_0/100\,.
\end{split}\end{equation}
Let
$G_{1},G_{2},\ldots,G_{s^{*}}$ (for some $s^*\le 2n$) be graphs on $[n]$, such that for each $s$ and $x\in V(G_s)$ we have $\LEFTDEG_{G_s}(x)\le D$, such that $\Delta(G_s)\le cn/\log n$, and such that the final $\delta n$ vertices of $G_s$ all have degree $d_s$ and form an independent set.

Let $H_0$ and $H_0^*$ be two edge-disjoint graphs on the same vertex set of order $n$ such that $(H_0,H_0^*)$ is
$(\frac14\alpha_{0},2D+3)$-coquasirandom, and
 $\sum_{s\in[s^{*}]}e(G_{s})\le e(H_0)-\gamma n^{2}$.
  \end{setting}

  Note that in~\eqref{eq:defconsts} we give numbers $\alpha_x$ which we call `constant' even though $n$ appears in their definition. Observe that $\alpha_x$ is strictly increasing in $x$. We will be interested only in values $0\le x\le 2n$ (though it is technically convenient to have the definition for all $x\in\mathbb{R}$), and it is easy to check that neither $\alpha_0$ nor $\alpha_{2n}$ depends on $n$.   
  
\smallskip  
  
The main lemmas for the analysis of \PackingProcess{} are now the following.
  Lemma~\ref{lem:firstquasirandomness} states that $(H_0,H^*_0)$ is coquasirandom with high probability. 
Lemma~\ref{lem:aggregate}  states that with high probability $(H_s,H^*_0)$ continues to be coquasirandom for each stage $s$. To prove this lemma will be the main work of this paper. Lemma~\ref{lem:dietsimplified} states that, provided that $H_{s}$ has the quasirandomness provided by Lemma~\ref{lem:aggregate}, the \RandomEmbedding\ of $G_{s+1}$ into $H_{s}$ is very likely to succeed.
Lemma~\ref{lem:degHstar} states that with high probability very few edges of $H^*_0$ are removed at each vertex to form $H^*_s$. This then implies that $(H_s,H^*_s)$ is also likely to be coquasirandom (though with a much worse error parameter). Finally, in Lemma~\ref{lem:completion}, using the coquasirandomness of $(H_s,H^*_s)$, we argue that at each stage it is very likely that the completion phase is possible.
  
We start with the lemma concerning the coquasirandomness of the initial bulk and reservoir.
  
\begin{lem}
\label{lem:firstquasirandomness}
For each $D\in\mathbb{N}$ and each $\gamma>0$, and for each $n$ sufficiently large, let us suppose that the constants $\alpha_0$ and $\xi$ are as in Setting~\ref{set:graphs}.

Suppose that $\widehat H$ is a $(\xi,2D+3)$-quasirandom graph of order $n$ and density $p\ge 3\gamma$. Let $H_0^*$ be a random subgraph of $\widehat H$ in which each edge of $\widehat H$ is kept with probability $q=\gamma/p$. Let $H_0$ be the complement of $H_0^*$ in $\widehat H$. Then with probability at least $1-n^{-6}$, we have that $e(H_0^*)= (1\pm\alpha_0)\gamma \binom{n}{2}$ and the pair $(H_0,H^*_0)$ is $\big(\tfrac14\alpha_{0},2D+3\big)$-coquasirandom.
\end{lem}

The next lemma states that coquasirandomness of $(H_s,H^*_0)$ is preserved.

\begin{lem}
  \label{lem:aggregate}For each $D\in\mathbb{N}$ and each $\gamma>0$, and for each $n$ sufficiently large, the following holds with probability
  at least $1-n^{-5}$. Suppose that the constants and 
  $G_{1},G_{2},\ldots,G_{s^{*}}$ and the graph $H_0\cup H_0^*=H$ are as in Setting~\ref{set:graphs}. 
  When \PackingProcess{} is run, for each $s\in[s^*]$ either \PackingProcess{} fails before completing stage $s$, or  the pair $(H_{s},H^*_0)$ is $(\alpha_{s},2D+3)$-coquasirandom.
\end{lem}

The next lemma estimates the probability that a single execution of \RandomEmbedding{} succeeds. 

\begin{lem}
\label{lem:dietsimplified}For each $D$, each $\gamma>0$, and any sufficiently large~$n$, let $\delta, \eta, \alpha_0,\alpha_{2n}, \varepsilon$ and $c$ be as in Setting~\ref{set:graphs}. Given any $\alpha_0\le\alpha\le\alpha_{2n}$, let~$G$ be a graph on vertex set~$[n]$ with  maximum degree at most
$cn/\log n$ such that $\deg^-(x)\le D$ for each $x\in V(G)$, and let~$H$ be any $(\alpha,2D+3)$-quasirandom $n$-vertex
graph with at least~$\gamma\binom{n}{2}$ edges. When \RandomEmbedding{} 
is run then it fails with probability at most $2n^{-9}$.
\end{lem}

Our final two main lemmas concern the completion phase of \PackingProcess{}.
The first states that the completion phase is likely to delete very few edges at any vertex of~$H^*_0$.

\begin{lem}\label{lem:degHstar}
 Given $D\in\mathbb{N}$ and $\gamma>0$, let $n$ be sufficiently large. Suppose that the constants and
  $G_{1},G_{2},\ldots,G_{s^{*}}$ and $H$ are as in Setting~\ref{set:graphs}. 
  When \PackingProcess{} is run, with probability at least $1-n^{-50}$ one of the following three events occurs. First, \PackingProcess\ fails. Second, there is some $s\in[s^*]$ such that $(H_s,H^*_0)$ is not $(\alpha_s,2D+3)$-coquasirandom. Third, for each $s\in[s^*]$ and $v\in V(H^*_s)$ we have $\deg_{H^*_0}(v)-\deg_{H^*_s}(v)\le 50\gamma^{-D}D\delta n$, and $(H_s,H^*_s)$ is $(\eta,2D+3)$-coquasirandom.
\end{lem}
We will show in the proof of Theorem~\ref{thm:maintechGENERAL} that the first two events are unlikely, so that the likely event is the last.

Our last lemma states that with high probability, at any stage $s$, provided $(H_{s-1},H^*_{s-1})$ is sufficiently coquasirandom, running \RandomEmbedding{} to partially embed~$G_s$ into $H_{s-1}$ is likely to give a partial embedding which can be completed to an embedding of $G_s$ using $H^*_s$.

\begin{lem}\label{lem:completion}
 For each $D\in\mathbb{N}$ and each $\gamma>0$, and for each $n$ sufficiently large, let the constants be as in Setting~\ref{set:graphs}. Suppose that
  $G$ is a graph on $[n]$, such that we have $\LEFTDEG(x)\le D$ for each $x\in V(G)$, we have $\Delta(G)\le cn/\log n$, and such that the final $\delta n$ vertices of $G$ form an independent set, and all have degree $d$. Suppose $(H,H^*)$ are a pair of $(\eta,2D+3)$-coquasirandom graphs on $n$ vertices, and $H$ is $(\alpha_{s^*},2D+3)$-quasirandom, with $e(H)=p\binom{n}{2}$ and $e(H^*)=(1\pm\eta)\gamma\binom{n}{2}$, where $p\ge\gamma$. When \RandomEmbedding\ is run to embed $G[{\scriptstyle [n-\delta n]}]$ into $H$, with probability at least $1-5n^{-9}$ it returns a partial embedding $\phi$ which can be extended to an embedding $\phi^*$ of~$G$ into $H\cup H^*$, with all the edges using a vertex in $\{n-\delta n+1,\dots,n\}$ mapped to~$H^*$.
\end{lem}

  Let us briefly explain why we cannot simply perform the whole embedding in the quasirandom $\widehat{H}$, but have to split it into a bulk and a reservoir. In order to analyse \RandomEmbedding, we require that the bulk is very quasirandom, but \RandomEmbedding\ is very well-behaved and preserves this good quasirandomness. In contrast, we are not able to show that the completion stage, where we choose a system of distinct representatives for the remaining vertices, is so well-behaved. If we used the bulk for this embedding the errors would rapidly become unacceptably large. However, to show that choosing such a system of distinct representatives is possible, we do not need much quasirandomness. Thus the reservoir $H^* _s$ does rapidly lose its quasirandomness (compared to $H_s$), but it is sufficient for the completion.

\smallskip

We now argue that our main lemmas imply Theorem~\ref{thm:maintechGENERAL}.
\begin{proof}[Proof of Theorem~\ref{thm:maintechGENERAL}]
We can assume that $p>3\gamma$ as the statement is vacuous otherwise.

Suppose that we run \PackingProcess{} on the input graphs $G_1,\ldots,G_{s^*}$. For the course of the analysis of this run, we shall first ignore possible failures during the completion phase. That is, if any failure during the completion phase occurs, we ignore it and continue embedding using \RandomEmbedding{} into the bulk. Clearly, this does not change behaviour of future rounds of \RandomEmbedding{} or the evolution of the bulk.

As we said earlier, we need to argue that  with positive probability \PackingProcess{} does not fail. Rather than proving this directly, we introduce additional quasirandomness conditions, and prove that with positive probability, all these conditions are satisfied up to any given stage, and that if we have the said quasirandomness conditions up to that stage, then \RandomEmbedding{} will proceed successfully through the next stage. (Of course, it could happen that \PackingProcess{} succeeds in the overall embedding even though some of our quasirandomness conditions failed during the course of the packing; we shall pessimistically view such an execution of \PackingProcess{} as unsuccessful.) More precisely, it is clear that \PackingProcess{} does not fail (in the \RandomEmbedding{} stage) unless at least one of the following
exceptional events occurs:
\begin{enumerate}[label=(\roman*)]
	\item\label{f:1} $(H_0,H^*_0)$ is not $(\frac14\alpha_0,2D+3)$-coquasirandom.
	\item\label{f:2} \RandomEmbedding{} proceeded through stages $s=1,\ldots,r$ (for some $r\in [s^*-1]$) without failure,  the pairs $(H_s,H_0^*)$ are $(\alpha_s,2D+3)$-coquasirandom for $s< r$, and $(H_r,H_0^*)$ is not an $(\alpha_r,2D+3)$-coquasirandom pair.
	\item\label{f:3} \RandomEmbedding{} proceeded through stages $s=1,\ldots,r$ (for some $r\in \{0,\ldots,s^*-1\}$) without failure, the graphs $H_s$ are $(\alpha_s,2D+3)$-quasirandom for $s\le r$. Then, in stage $r+1$, \RandomEmbedding{} fails.
\end{enumerate}
Lemma~\ref{lem:firstquasirandomness} gives an upper bound on the probability of the event in~\ref{f:1}. Lemma~\ref{lem:aggregate} gives an upper bound on the probability of all the events in~\ref{f:2}. For each fixed $r\in \{0,\ldots,s^*-1\}$, the event in~\ref{f:3} can be bounded using Lemma~\ref{lem:dietsimplified}. Thus, the probability that \PackingProcess{} fails in the \RandomEmbedding{} part is at most
$n^{-6}+n^{-5}+s^*\cdot 2n^{-9}$.

Let us now analyse the completion phases of \PackingProcess{}. If \PackingProcess{} fails in one of the completion phases then one of the following events occurs:
\begin{enumerate}[label=(\roman*),resume]
	\item One of the events described under~\ref{f:1}-\ref{f:3}.
	\item\label{f:4} None of~\ref{f:1}-\ref{f:3} occurs. \RandomEmbedding{} and the completion phase proceed successfully through the first $r$ stages (for some $r\in\{1,\ldots,s^*-1\}$. For $s\in[r]$ all the pairs $(H_s,H_0^*)$ are $(\alpha_s,2D+3)$-coquasirandom. However, there is a stage $s\in[r]$ where $(H_s,H^*_s)$ is not $(\eta,2D+3)$-coquasirandom.
	\item\label{f:5} None of~\ref{f:1}-\ref{f:3} occurs. \RandomEmbedding{} and the completion phase proceeds successfully through the first $r$ stages (for some $r\in\{0,\ldots,s^*-1\}$, and throughout all the pairs $(H_s,H_0^*)$ and $(H_s,H^*_s)$ are $(\alpha_s,2D+3)$-coquasirandom and $(\eta,2D+3)$-coquasirandom, respectively. In stage $r+1$, \RandomEmbedding{} successfully embeds but the completion phase fails.
\end{enumerate}
Lemma~\ref{lem:degHstar} bounds the probability of the event in~\ref{f:4} by $n^{-50}$. Finally, Lemma~\ref{lem:completion} bounds the probability of events in~\ref{f:5} for each given $r$ by $5n^{-9}$. Thus, the total probability of failure due to~\ref{f:4} or~\ref{f:5} is at most $n^{-50}+s^*\cdot 5n^{-9}$.

We conclude that \PackingProcess{} packs the graphs $G_1,\ldots,G_{s^*}$ into $\widehat H$ with positive probability.
\end{proof}

\subsection{The probability space for \RandomEmbedding}\label{ssec:AlgorithmProbSpace}
Algorithm~\ref{alg:embed} gives a sound definition of a randomised algorithm
which either provides an embedding of $G{\scriptstyle [n-\delta n]}$ into $H$ or
fails, and the probability of any output can be in principle computed. To handle
the analysis of \RandomEmbedding, which is the most demanding part of this
paper, it is useful to properly set up a probability space as indicated at the
beginning of Section~\ref{ssec:seqdependence}. Given $G$ and $H$ as in
Algorithm~\ref{alg:embed} (recall that $V(G)=[n]$), let $\Omega^{G\AlgMap
  H}:=(V(H)\cup\{\frownie\})^{n-\delta n}$. We now need to define the
probability measure on $\Omega^{G\AlgMap H}$. Let
$\bm{\omega}=(\omega_1,\ldots,\omega_{n-\delta n})\in \Omega^{G\AlgMap H}$ be
given. Suppose first that $\bm{\omega}$ consists only of vertices of $V(H)$.
Then we define $\PROB^{G\AlgMap H}(\bm{\omega})$ as the probability that
\RandomEmbedding{} succeeds embedding $G{\scriptstyle [n-\delta n]}$ into $H$,
and maps each vertex $t\in [n-\delta n]$ of $G$ on vertex $\omega_t$. Suppose
next that $\bm{\omega}$ contains some $\frownie$'s, and that these form a
terminal segment of $\bm{\omega}$, say starting from position $t_0$. Then we
define $\PROB^{G\AlgMap H}(\bm{\omega})$ as the probability that
\RandomEmbedding{} succeeds in the first $t_0-1$ steps, and for each
$t\in[t_0-1]$ it maps vertex $t$ on $\omega_t$, and then in step $t$ it halts
with failure. Last, suppose that $\bm{\omega}$ contains some $\frownie$'s but
these do not form a terminal segment of $\bm{\omega}$. We then define
$\PROB^{G\AlgMap H}(\bm{\omega}):=0$. It is clear that $\PROB^{G\AlgMap
  H}(\bm{\omega})$ is a probability measure on $\Omega^{G\AlgMap H}$ which
corresponds to possible runs of \RandomEmbedding{}.

We shall use the concept of histories and history ensembles, as introduced in Section~\ref{ssec:seqdependence}, in connection with $\Omega^{G\AlgMap H}$.

\subsection{Organisation of the technical part of the paper}
It thus remains to prove all the main lemmas from this section. 
Lemmas~\ref{lem:firstquasirandomness} and~\ref{lem:aggregate} are proven in Section~\ref{sec:maintainingquasirandomness}.
Lemma~\ref{lem:dietsimplified} is stated here in a simplified form. In actuality, we prove a stronger statement (of which Lemma~\ref{lem:dietsimplified} is a straightforward consequence) in Lemma~\ref{lem:diet}. This stronger form is also needed for proving Lemma~\ref{lem:aggregate}, and its proof spans the entire Section~\ref{sec:diet}.
Lemmas~\ref{lem:degHstar} and~\ref{lem:completion} are proven in Section~\ref{sec:completion}.

\section{Staying on a diet}\label{sec:diet}

In this section we consider the running of \RandomEmbedding\ to embed one degenerate graph $G$ into a quasirandom graph $H$. The results of this section will always be used to analyse one stage $s$, when we take $G=G_s$ and $H=H_{s-1}$. We also analyse how \RandomEmbedding{} behaves with respect to the graph $H^*=H^*_{s-1}$. We analyse carefully how fast common neighbourhoods of vertices in $H$ are eaten up by \RandomEmbedding, and how often individual vertices of $H$ appear in candidate sets. To make this precise, we introduce the following two definitions.

The diet condition states that during the running of \RandomEmbedding, for
each $t\in[n-\delta n]$, the fraction of each set $\NBH_{H}(S)$ which is
covered by $\im(\psi_{t})$ is roughly as expected, that is, roughly
proportional to $|\im(\psi_{t})|/n$. As with (co)quasirandomness, we also require a codiet condition, considering the intersection of some vertex neighbourhoods in $H$ and $H^*$.

\begin{defn}[diet condition, codiet condition]
  \label{def:dietcondition} 
  Let~$H$ be a graph with $n$ vertices and $p\binom{n}{2}$ edges, and let
  $X\subseteq V(H)$ be any vertex set. We say that the pair $(H,X)$
  satisfies the \emph{$(\beta,L)$-diet condition} if for every set
  $S\subset V(H)$ of at most~$L$ vertices we have $|\NBH_{H}(S)\setminus
  X|=(1\pm\beta)p^{|S|}(n-|X|)$.
  
    Let~$H, H^*$ be two graphs with vertex set~$V$ of order~$n$ and $p\binom{n}{2}$ and $p^*\binom{n}{2}$ edges, respectively, and let
  $X\subseteq V$ be any vertex set. We say that the triple $(H, H^*,X)$
  satisfies the \emph{$(\beta,L)$-codiet condition} if for every set
  $S\subset V$ of at most~$L$ vertices and for every subset $R\subseteq S$ we have \[\Big|\big(\NBH_{H}(R)\cap \NBH_{H^*}(S\setminus R)\big)\setminus
  X\Big|=(1\pm\beta)p^{|R|}(p^*)^{|S\setminus R|}(n-|X|)\;.\]
\end{defn}
Observe that the $(\beta,L)$-diet condition holding for $(H,\emptyset)$ is simply the statement that~$H$ is $(\beta,L)$-quasirandom, and similarly for the codiet condition.

\ENHANCED{
We now give a definition of the diet condition relative to a given set of vertices. This definition is used in the `furthermore'	
\begin{defn}[diet condition relative to a set]
	\label{def:dietconditionrelative} 
	Let~$H$ be a graph on $n$ vertices and $p\binom{n}{2}$ edges, and let
	$X\subseteq V(H)$ be any vertex set. We say that the pair $(H,X)$
	satisfies the \emph{$(\beta,L)$-diet condition} if for every set
	$S\subset V(H)$ of at most~$L$ vertices we have $|\NBH_{H}(S)\setminus
	X|=(1\pm\beta)p^{|S|}(n-|X|)$.
\end{defn}
}

The cover condition, defined below, roughly states that for each~$v$ in the
host graph~$H$ during the embedding of~$G$ into~$H$ 
by \RandomEmbedding{}, the right fraction of vertices~$x$ of~$G$ have~$v$ in their final candidate set.
For making precise what we mean by `the right fraction' some care is
needed. Firstly, how likely it is that~$v$ is in the final candidate set
of~$x$ depends on the number neighbours of~$x$ preceding~$x$. Therefore we
will partition $V(G)$ according to this number of previous neighbours.  For
technical reasons we actually further want to control this fraction in intervals of $V(G)$ of length~$\eps n$, where $n$ is the order of $H$. Hence we define for a given $\varepsilon>0$ the set \[X_{i,d}:=\{x\in
V(G)\colon i\le x<i+\varepsilon n,|\LNBH(x)|=d\}\,.\] 
When $G$ is given with a $D$-degenerate ordering it is enough to consider $d\in\{0,1,\ldots,D\}$. \corrig{That is, 
\begin{equation}\label{eq:SumXid}
\sum_{d=0}^D |X_{i,d}|=\epsilon
n\;.
\end{equation}
}
So if $H$
is quasirandom and has $p\binom{n}{2}$ edges, then for an arbitrary
$v\in V(H)$, we would expect that about a $p^{d}$-fraction of
vertices $x$ in each~$X_{i,d}$ have $v$ in their final candidate sets (let us remind that the candidate set may include also vertices used by the embedding). \corrig{However, this expectation turns out not to be quite true. If a vertex $y$ of $G$ has been embedded to $v$, and a vertex $x$ has a left-neighbour $z$ which is adjacent to $y$, then $x$ is much more likely to have $v$ in its candidate set, because we get `for free' that $z$ is embedded to a neighbour of $v$. This is the only reason why the above intuition can fail: in particular the expectation does hold true if $v$ is not in the image of $\psi$, and it turns out that this is all we need.}

\begin{defn}[cover condition]
\label{def:cover}
Suppose that $G$ and $H$ are two graphs such that $H$ has order $n$, the vertex set of $G$ is $[n]$, and $H$ has density $p$.  Suppose that numbers $\beta,\epsilon>0$ and $i\in [n-\epsilon n]$ are given. \corrig{Suppose that $\psi$ is a partial embedding of $G$ into $H$ which embeds at least the first $i+\eps n-1$ vertices of $G$. We say that $\psi$ satisfies the \emph{$(\varepsilon,\beta,i)$-cover condition} if for each $v\in V(H)$ such that $v\not\in\im\psi_{\restriction [i+\eps n-1]}$}, and for each $d\in\mathbb N$,
 we have
  \[
  \big|\big\{ x\in X_{i,d}:v\in\NBH_{H}\big(\psi(\LNBH(x))\big)\big\}\big|=(1\pm\beta)p^{d}|X_{i,d}|\pm\varepsilon^{2}n\,.
  \]
Note that a corresponding condition for $d=0$ is trivial, even with zero error parameters.
\end{defn}

We use Definitions~\ref{def:conj-quasirandom}, \ref{def:dietcondition} and~\ref{def:cover}, to define key events $\DietEvent(\cdot;\cdot)$, $\CoverEvent(\cdot;\cdot)$, $\CodietEvent(\cdot)$ on $\Omega^{G\AlgMap H}$.
\begin{defn}
 Suppose that $D$, $\delta$ and $\epsilon$ are as in Setting~\ref{set:graphs}.
 Suppose that $\lambda>0$. Suppose that we have graphs $G$ and $H$ as in
 Algorithm~\ref{alg:embed}. Suppose that we run \RandomEmbedding{} to partially
 embed $G$ into $H$. Let $(\psi_{i})_{i\in[t_*]}$ be the partial embeddings of
 $G\big[[i]\big]$ into $H$, where $t_*=n-\delta n$ if \RandomEmbedding{} succeeded, and otherwise
 $t_*+1$ is the step in which \RandomEmbedding{} halted with failure.
\begin{itemize}
	\item For each $t\in[n-\delta n]$, let $\DietEvent(\lambda;t)\subset\Omega^{G\AlgMap H}$ correspond to executions of \RandomEmbedding{} for which $t_*\ge t$ and the pair $(H,\im\psi_{t})$
	satisfies the $(\lambda,2D+3)$-diet condition.
	\item For each $t\in[n-\delta n]$, let $\CoverEvent(\lambda;t)\subset\Omega^{G\AlgMap H}$ correspond to executions of \RandomEmbedding{} for which $t_*\ge t+\eps n$ and the embedding $\psi_{t^*}$ of~$G$ into~$H$
	satisfies the $(\varepsilon,\lambda,t)$-cover condition.
	\item Suppose further that we have a graph $H^*$ with $V(H)=V(H^*)$. For each $t\in[n-\delta n]$, let $\CodietEvent(t)\subset\Omega^{G\AlgMap H}$ correspond to executions of \RandomEmbedding{} for which $t_*\ge t$ and the triple $(H,H^*,\im\psi_t)$ satisfies the $(2\eta,2D+3)$-codiet condition.
\end{itemize}
\end{defn}
Note that the events $\DietEvent(\cdot;t)$ and $\CodietEvent(t)$ are determined by histories (as defined in Sections~\ref{ssec:seqdependence} and~\ref{ssec:AlgorithmProbSpace}) up to time $t$. That is, for any $\lambda>0$ and any history $\hist_t$, we have that $\DietEvent(\lambda;t)$ either contains $\hist_t$ or is disjoint from $\hist_t$. We have similar the same property for $\CodietEvent(t)$. The event $\CoverEvent(\cdot;t)$ is somewhat different since its definition involves the set $X_{t,d}$ which looks $\epsilon n-1$ many steps forward in time. So, for any history $\hist_{t+\epsilon n-1}$, we have that $\CoverEvent(\lambda;t)$ either contains $\hist_{t+\epsilon n-1}$ or is disjoint from $\hist_{t+\epsilon n-1}$.

\medskip

The following lemma is the crucial accurate analysis of \RandomEmbedding\ which we need in order to show that \RandomEmbedding\ is likely to succeed and in order to derive further properties of the final embedding.

\begin{lem}[Diet-and-cover lemma]
  \label{lem:diet}For each $D\in\mathbb N$, each $\gamma>0$, and any sufficiently large~$n$, let $\delta, \eta, \alpha_0,\alpha_{2n}, \varepsilon$ and $c,C$ be as in Setting~\ref{set:graphs}. Let $\alpha\in[\alpha_0,\alpha_{2n}]$ be arbitrary. Let~$G$ be a graph on vertex set~$[n]$ with  maximum degree at most
  $cn/\log n$ such that $\deg^-(x)\le D$ for each $x\in V(G)$, and let~$H$ be any $(\alpha,2D+3)$-quasirandom $n$-vertex
  graph with at least~$\gamma\binom{n}{2}$ edges. Suppose in addition that~$H^*$ is a graph on $V(H)$ such that $(H,H^*)$ is $(\eta,2D+3)$-coquasirandom. Then we have
\begin{equation}\label{eq:Th}
\PROB^{G\AlgMap H}\left(
\bigcap_{t\in[n-\delta n]} \DietEvent(C\alpha;t)\cap \bigcap_{t\in [n+1-\epsilon n]} \CoverEvent(C\alpha;t) \cap \bigcap_{t\in[n-\delta n]} \CodietEvent(t)
\right)\ge 1-2n^{-9}\;.
\end{equation}
\end{lem}

This lemma immediately implies Lemma~\ref{lem:dietsimplified}.
\begin{proof}[Proof of Lemma~\ref{lem:dietsimplified}]
 Recall that \RandomEmbedding\ fails if and only if $\CANDSET_{G\AlgMap H}^{t-1}(t)\setminus\im(\psi_{t-1})=\emptyset$ for some $t$, and $\DietEvent(C\alpha;t-1)$ in particular gives a formula lower bounding the size of $\CANDSET_{G\AlgMap H}^{t-1}(t)\setminus\im(\psi_{t-1})$ which is greater than $0$. Since the likely event of Lemma~\ref{lem:diet} is contained in $\DietEvent(C\alpha;t-1)$ for each $t\ge 2$, and the same lower bound is trivially implied by $(\alpha,2D+3)$-quasirandomness of $H$ for $t=1$ (since $\im\psi_0=\emptyset$), we conclude that within the likely event of Lemma~\ref{lem:diet}, \RandomEmbedding\ does not fail.
\end{proof}

The main difficulty is to establish that the cover and diet conditions hold. We will see that the codiet condition is an easy byproduct. The reason for the difficulty is that the error terms in the cover and diet conditions for small times $t$ feed back into the calculations which will establish the cover and diet conditions for larger times $t$, and we have to ensure that this feedback loop does not allow the errors to spiral out of control. To that end, we define a new sequence of error terms, which we need only in the proof of Lemma~\ref{lem:diet}. 
The following constants
$\{\beta_t:t\in\mathbb{R}\}$ are a carefully chosen increasing
sequence (depending on $\alpha$) such that $\beta_{0}=\alpha$ and such that 
$\beta_{n}/\beta_{0}$ is bounded by a constant which does not depend on $\alpha$ (though it does depend on $D$, $\gamma$ and $\delta$). Given~$D$ and $\alpha,\delta,\gamma>0$, we define
\begin{equation}\label{eq:defnbeta}
  \beta_{t}:=2\alpha\exp\big(\tfrac{1000D\delta^{-2}\gamma^{-2D-10}t}{n}\big)\;.
\end{equation}
We will mainly take $t$ integer in the range $[0,n]$, but it is convenient to allow $t$ to be any real number. In particular, for each $t\ge 0$, we have
\begin{equation}\label{eq:betabound}\begin{split}
 &\tfrac{1}{n}\int_{i=0}^t 1000D\delta^{-2}\gamma^{-2D-10}\beta_i \,\mathrm{d}i\\
 \le &2\alpha\int_{i=-\infty}^t \frac{1000D\delta^{-2}\gamma^{-2D-10}}{n}\exp\big(\tfrac{1000D\delta^{-2}\gamma^{-2D-10}i}{n}\big) \,\mathrm{d}i=\beta_t\,.
\end{split}\end{equation}
Suppose that we have Setting~\ref{set:graphs}, and suppose that $\alpha\ge \alpha_{0}$ is given. Then for each $t\ge 0$ we have
\begin{equation}\label{eq:ttt}
\beta_t \gamma^{2D+3}\delta\ge \beta_0 \gamma^{2D+3}>\epsilon\;.
\end{equation}

\medskip

We split the proof of Lemma~\ref{lem:diet} into two parts. The cover lemma (Lemma~\ref{lem:cover}) states that if the $(\beta_t,2D+3)$-diet condition holds for $(H,\im\psi_i)$
for each $i\in[t-1]$, then it is very
unlikely that the $(\eps,20D\beta_t,t)$-cover condition fails for $\psi_{t+\eps n-2}$. Note that the time $t+\eps n-2$ is the first time at which the $(\eps,20D\beta_t,t)$-cover condition is guaranteed to be determined, since at this time all left-neighbours of all vertices $t,t+1,\dots,t+\eps n-1$ have certainly been embedded.

\begin{lem}[Cover lemma]
  \label{lem:cover}
  For each $D$, each $\gamma>0$ and sufficiently large $n$, let $\alpha_0,\alpha_{2n}, \varepsilon, \delta $ and~$c$ be as in Setting~\ref{set:graphs}. Suppose
  that $\alpha_0\le\alpha\le\alpha_{2n}$ and $G$ is a graph on vertex set~$[n]$, with $\LEFTDEG(x)\le D$ for each $x\in[n]$, with maximum degree at most
  $cn/\log n$, 
  and suppose that $H$ is an $n$-vertex graph of density at least~$\gamma$.
  Let~$\beta_t$ for $0\le t\le n$ be defined as
  in~\eqref{eq:defnbeta} and assume that $\beta_n\le\frac1{10}$. Let~$t$
  with $1\le t\le n-\delta n-\eps n+1$ be fixed.

Then we have
$$
\PROB^{G\AlgMap H}\left(
\bigcap_{i=1}^{t-1}\DietEvent(\beta_t;i)
\setminus
\CoverEvent(20D\beta_t;t)
\right)
\le n^{-10}
\;.
$$
\end{lem}

Let us consider Setting~\ref{set:graphs}. Suppose that for some $0\le t\le n-\delta n-\eps n$, \RandomEmbedding\ runs up to time $t$ and the $(\beta_t,2D+3)$-diet condition holds for $(H,\im\psi_t)$. Let $p:=e(H)/\binom{n}{2}$ and suppose that $p\ge \gamma$. Then for each $t+1\le j\le t+\eps n$, and each set $S\subset V(H)$ of at most $2D+3$ vertices, we have
\begin{align*}
|\NBH_{H}(S)\setminus \im\psi_j|&\ge |\NBH_{H}(S)\setminus \im\psi_t|-\eps n\\
\JUSTIFY{diet for $(H,\im\psi_t)$}
&\ge (1-\beta_t)p^{|S|}(n-|\im\psi_t|)-\eps n\\
\JUSTIFY{$\eps<\beta_t\gamma^{2D+3}\delta$ by~\eqref{eq:ttt}}&\ge (1-2\beta_t)p^{|S|}(n-|\im\psi_t|)\;.
\end{align*}
Hence, the $(2\beta_t,2D+3)$-diet condition holds deterministically for $(H,\im\psi_j)$. In particular \RandomEmbedding{} cannot fail before time $t+\varepsilon n$. 

The diet lemma (Lemma~\ref{lem:dietcon}) states that when the $(\beta_i,2D+3)$-diet condition holds for
$(H,\im\psi_i)$ for each $i\in[t-1]$, and the $(\eps,20D\beta_i,i)$-cover condition holds for
$\psi_{i+\eps n-2}$ for each $i\in[t+1-\eps n]$, then it is unlikely that the
$(\beta_t,2D+3)$-diet condition fails for $(H,\im\psi_t)$. We also obtain the desired codiet condition.

\begin{lem}[Diet lemma]
\label{lem:dietcon}For each $D$, each  $\gamma>0$, and any sufficiently large~$n$, let $\alpha_0,\alpha_{2n}, \varepsilon, \delta$ and $\eta$ be as in Setting~\ref{set:graphs}.  For any $t\le(1-\delta)n$, and $\alpha_0\le\alpha\le\alpha_{2n}$ the following holds. Suppose that
$G$ is a graph on $[n]$ such that $\LEFTDEG(x)\le D$ for each $x\in[n]$, and $H$ is an ($\alpha, 2D+3$)-quasirandom graph with $n$
vertices with $p\binom{n}{2}$ edges, with $p\ge \gamma$. Suppose furthermore that $H^*$ is a graph on $V(H)$ and $\hat p \binom{n}{2}$ edges with $\hat p\ge (1-\eta)\gamma$, such that  $(H,H^*)$ satisfies the $(\eta,2D+3)$-coquasirandomness condition. Let~$\{\beta_\tau:\tau\in[0,n]\}$ be defined as
in~\eqref{eq:defnbeta} and assume that $\beta_n\le\frac1{10}$. Let $t$ with $1\le t\le n-\delta n$ be fixed.

Then we have
$$
\PROB^{G\AlgMap H}\left(
\bigcap_{j=1}^{t-1}\DietEvent(\beta_j;j)
\cap
\bigcap_{j=1}^{t+1-\epsilon n}\CoverEvent(20D\beta_j;j)
\setminus
(\DietEvent(\beta_t;t)\cap \CodietEvent(t))
\right)
\le n^{-10}
\;.
$$
\addtocounter{counterDiet}{\value{thm}}
\end{lem}

Since the graphs $G$ and $H$ are fixed in Lemmas~\ref{lem:diet}, \ref{lem:cover}, and~\ref{lem:dietcon},
in this section we drop the subscript in the notation $\CANDSET_{G\AlgMap H}^{j}(x)$
and write simply $\CANDSET^{j}(x)$. Likewise, we write $\PROB$ instead of $\PROB^{G\AlgMap H}$. Last, we write 
$(\psi_{i})_{i\in t_*}$ for partial embeddings of~$G$ into~$H$; here $t_*$ is the time at which \RandomEmbedding{} halts. Of course, $t_*$ and $(\psi_{i})_{i\in t_*}$ depend on a particular realization $\bm{\omega}\in\Omega^{G\AlgMap H}$ of the run of \RandomEmbedding{}.

We now show that Lemmas~\ref{lem:cover} and~\ref{lem:dietcon}, whose proofs are deferred to later in this section, imply Lemma~\ref{lem:diet}.

\begin{proof}[Proof of Lemma~\ref{lem:diet}]
	Suppose that we are given $D$ and $\gamma$. Now, given $\alpha>0$, we define $\beta_t$ for each $0\le t\le n$ as in~\eqref{eq:defnbeta}. For $t=0,\ldots,n-\delta n$, define 
	\begin{equation}\label{eq:kafe}
	\cA_t:=\bigcap_{j=1}^t\DietEvent(\beta_j;j)\cap\bigcap_{j=\epsilon n}^t\CoverEvent(20D\beta_{t-\epsilon n+1};j-\epsilon n+1)\cap\bigcap_{j=1}^t\CodietEvent(j)\;.
	\end{equation}
	Our strategy is first to show that $\PROB(\cA_{t-1}\setminus\cA_{t})$ is tiny for each $t$. Since $\PROB(\cA_0)=1$, this will imply that $\PROB(\cA_{n-\delta n})$ is very close to~1. Last, we shall show that $\cA_{n-\delta n}$ is a subset of the event in~\eqref{eq:Th}.
	
	Indeed, suppose that the event $\cA_{t-1}$ holds.  This in particular means that the $(\beta_j,2D+3)$-diet condition holds for $(H,\im\psi_j)$ for each $1\le j<t$, and the $(\eps,20D\beta_{j-\eps n+1},j-\eps n+1)$-cover condition holds for~$\psi_j$ for each $\eps n-1\le j< t$.
	
	%Firstly, we show that \RandomEmbedding{} cannot fail at time~$t$. Then we use Lemma~\ref{lem:cover} to show that with high probability the $(\eps,20D\beta_{t-\eps n+2},t-\eps n+2)$-cover condition holds for $\psi_t$. 
	Because the $(\beta_{t-1},2D+3)$-diet condition holds for ($H$, $\psi_{t-1}$), picking $S=\psi_{t-1}(\LNBH(t))$, we have $\big|\CANDSET^{t-1}(t)\setminus\im\psi_{t-1}\big|=\big|\NBH_H(S)\setminus\im\psi_{t-1}\big|>0$. It follows that \RandomEmbedding{} cannot fail at time~$t$. 
	
	Firstly, let us focus on the term $\CoverEvent(20D\beta_{t-\epsilon n+1};t-\epsilon n+1)$ in~\eqref{eq:kafe}. This term does not exist when $t<\epsilon n$, so let us assume the contrary. Lemma~\ref{lem:cover} then tells us that
	$$
	\PROB\left(
	\bigcap_{i=1}^{t-\epsilon n}\DietEvent(\beta_{t-\epsilon n+1};i)
	\setminus
	\CoverEvent(20D\beta_{t-\epsilon n+1};t-\epsilon n+1)
	\right)
	\le n^{-10}\;.$$
In particular, 
\begin{equation}\label{eq:beansI}
\PROB(\cA_{t-1}\setminus 	\CoverEvent(20D\beta_{t-\epsilon n+1};t-\epsilon n+1))\le n^{-10}\;.
\end{equation}
	
	Secondly, we  use Lemma~\ref{lem:dietcon} to show that with high probability neither the diet condition nor the codiet condition fails at time~$t$. Indeed, Lemma~\ref{lem:dietcon} tells us that
	\[\PROB\left(
	\bigcap_{j=1}^{t-1}\DietEvent(\beta_j;j)
	\cap
	\bigcap_{j=\epsilon n}^{t-1}\CoverEvent(20D\beta_{j+1-\epsilon n};j+1-\epsilon n)
	\setminus
	(\DietEvent(\beta_t;t)\cap \CodietEvent(t))
	\right)
	\le n^{-10}
	\]
In particular, 
\begin{equation}\label{eq:beansII}
\PROB(\cA_{t-1}\setminus 
(\DietEvent(\beta_t;t)\cap \CodietEvent(t))
)
\le n^{-10}\;.
\end{equation}
		
Summing up~\eqref{eq:beansI} and~\eqref{eq:beansII}, we conclude that $\PROB(\cA_{t-1}\setminus\cA_{t})\le 2n^{-10}$. Taking a union bound over the at most $n$ choices of~$t$, we see that with probability at least $1-2n^{-9}$ the good event from the statement of Lemma~\ref{lem:diet} holds, i.e.,  that \RandomEmbedding{} does not fail, and by the choice of~$C$ and by~\eqref{eq:defnbeta}, for each $1\le t\le (1-\delta)n$ the pair $(H,\im\psi_t)$ satisfies the $(C\alpha,2D+3)$-diet condition and the triple $(H,H^*,\im\psi_t)$ satisfies the $(2\eta,2D+3)$-codiet condition, and for each $1\le t\le n+1-\eps n$ the embedding $\psi_{(1-\delta)n}$ satisfies the $(\eps,C\alpha,t)$-cover condition, as desired.
\end{proof}

We now prove the cover lemma.

\begin{proof}[Proof of Lemma~\ref{lem:cover}]
  Let $e(G)=p\binom{n}{2}\ge\gamma\binom{n}{2}$. 
  Let $\cD$ be the event
  that the $(\beta_{t},2D+3)$-diet condition holds for each
  $(H,\im\psi_{i})$ with $1\le i\le t-1$, $\cD:=\bigcap_{i=1}^{t-1}\DietEvent(\beta_t;i)$.  We also fix $1\le d\le
  D$. Define $\cB_{v,d}$ as the event that~$\cD$
  holds, and that~$v$ and~$d$ witness the failure of the
  $(\varepsilon,20D\beta_{t},t)$-cover condition for $\psi_{t+\eps n-2}$. \corrig{More formally,
  \begin{multline*}
  \cB_{v,d}:=\cD\cap\Big\{\bm{\omega}\in\Omega^{G\AlgMap H}:\quad\text{$v\not\in\im\psi_{t+\eps n-1}$ and}\\
  \big|\big\{ x\in X_{t,d}:v\in\NBH_{H}\big(\psi_{t+\eps n-2}(\LNBH(x))\big)\big\}\big|
  \neq (1\pm 20D\beta_{t})p^{d}|X_{t,d}|\pm\varepsilon^{2}n
  \Big\}
  \,.
  \end{multline*}}
  Our aim is to show that 
\begin{equation}\label{eq:PrBvd}
\PROB\left(\cB_{v,d}\right)\le n^{-12}/D\;.
\end{equation} 
  
  A union bound over the
  choices of~$v$ and~$d$ then gives the lemma.

\medskip

  Our strategy for proving~\eqref{eq:PrBvd} is as follows. Ideally, we would like to assert that for
  each $x\in X_{t,d}$ the probability of $v\in\CANDSET^{x-1}(x)$ is roughly
  $p^{d}$ and apply Lemma~\ref{lem:concentration} to bound the probability
  of the bad event~$\cB_{v,d}$. To this end, we consider a dynamical version of candidate sets, where we track changes in the set potentially suitable to accommodate $x$ as we gradually embed more and more left-neighbors of $x$. More precisely, for each $i\le x-1$, let $\CANDSET^{i,\mathsf{dyn}}(x):=\NBH_{H}\left(\psi_{x-1}\big([i]\cap\LNBH_{G}(x)\big)\right)$. At time $i=0$, we have
  $v\in\CANDSET^{i,\mathsf{dyn}}(x)$, and as $i$ increases, the set $\CANDSET^{i,\mathsf{dyn}}(x)$ shrinks exactly at times $y\in\LNBH(x)$ when left-neighbors of $x$ are embedded.
  
  Unfortunately we are not able to carry out this ideal strategy, because when we apply Lemma~\ref{lem:concentration} what we need to calculate is not the probability of $v\in\CANDSET^{x-1}(x)$, but this probability in the conditioned space given by the history up to some earlier time. Because the sets $\LNBH(x)$ interleave each other, this conditional probability will generally not be close to $p^d$ and we were not able to find a good way to estimate it. Hence we refine this strategy by rewriting the event $\{v\in\CANDSET^{x-1}(x)\}$ as
  \begin{equation}\label{eq:rme}
   \bigcap_{k=1}^d\{y_1,y_2,\ldots,y_k \AlgMap \NBH_H(v)\}\;,
  \end{equation}where $y_1,\ldots,y_d$ are the neighbours of $x$, ordered from left to right. The event $\{y_1,y_2,\ldots,y_d \AlgMap \NBH_H(v)\}$, of course, equals the entire intersection~\eqref{eq:rme}. However, this more complicated way of expressing~\eqref{eq:rme} suggests to introduce, for each $k$, a sequence of random variables that count the events of the form $\{y_1,y_2,\ldots,y_k \AlgMap \NBH_H(v)\}$, ordered by $y_k$. Intuitively, conditioning on $\{y_1,y_2,\ldots,y_k \AlgMap \NBH_H(v)\}$ holding (which is determined by the history up to the time at which we embed $y_k$) we should expect that the probability that $\{y_1,y_2,\ldots,y_{k+1} \AlgMap \NBH_H(v)\}$ holds is about $p$. We will be able to demonstrate this is true, even if we condition on a typical history up to the time immediately before embedding $y_{k+1}$, and this allows us to use Lemma~\ref{lem:concentration}.

  More formally, given $1\le k\le d$ and $y\in V(G)$, we define random variables
  $Y_{k,1},\dots,Y_{k,t+\eps n-2}$ as follows.  Let $Y_{k,y}$ be the number of
  vertices $x\in X_{t,d}$ such that $y$ is the $k$-th leftmost vertex of $\LNBH(x)$
  and the first $k$ vertices of $\LNBH(x)$ are all embedded to
  $\NBH_{H}(v)$.
  Further, for each $0\le k\le d$, we let $\cY_{k}$ be the event that \corrig{ either $v\in\im\psi_{t+\eps n-1}$ or }
  $(1\pm10\beta_{t})^{k}p^{k}|X_{t,d}|\pm k\varepsilon^{2}n/d$ vertices
  $x\in X_{t,d}$ have all of the first~$k$ vertices of $\LNBH(x)$ embedded
  to $\NBH_{H}(v)$. 
  Observe that the
  event $\cY_{k}$ is precisely the statement that
  \begin{equation}\label{eq:eventcalY_k}
  \corrig{\text{either $v\in\im\psi_{t+\eps n-1}$ or }}\quad\sum_{y=1}^{t+\eps n-2}Y_{k,y}=(1\pm10\beta_{t})^{k}p^{k}|X_{t,d}|\pm
  k\varepsilon^{2}n/d\,.
  \end{equation}
  Our bad event then satisfies
  \[\cB_{v,d} \subset  \cD \setminus \cY_d\,,\]
  because $(1\pm10\beta_{t})^{d}=1\pm 20D\beta_{t}$.
  In order to bound the probability of~$\cB_{v,d}$ we cover~$\cB_{v,d}$ with $d$ events, each of whose probabilities we can bound with
  Lemma~\ref{lem:concentration}.
  For this purpose we define the event 
  \[\cE_k = \cY_{k-1}\cap \cD\]
  for each $1\le k\le d$. Note that $\cE_1=\cD$ since $\cY_0$ holds trivially with probability one. We thus have
  \[\cB_{v,d} \subset  \cD \setminus \cY_d \subset \bigcup_{1\le k\le d} \big( \cE_k \setminus \cY_k \big)\,.\]
  Our aim then is to show that for each $1\le k\le d$ we have
  \begin{equation}\label{eq:cover:goal}
    \PROB(\cE_k \setminus \cY_k)\le n^{-12}/(d\cdot D)\,.
  \end{equation}
  Note that this and a union bound over the~$d$ choices of~$k$ gives~\eqref{eq:PrBvd}.

  To establish~\eqref{eq:cover:goal} we would like to apply
  Lemma~\ref{lem:concentration}. Hence we need to argue that either $\cE_k$ fails, or we can estimate
  $\sum_{y=1}^{t+\eps n-2}\EXP\left(Y_{k,y}|\hist_{y-1}\right)$, where $\hist_{y-1}$ is the history of
  embedding decisions taken in \RandomEmbedding{} up to and including the
  embedding of vertex~$y-1$.
  To this end, for $y\in[t+\eps n-2]$ let $Z_{k,y}$ be the number of vertices $x\in X_{t,d}$ such that $y$ is the
  $k$-th leftmost vertex of $\LNBH(x)$ and the first $k-1$ vertices of $\LNBH(x)$ are
  embedded to $\NBH_H(v)$. Then the quantity
  $Z_{k,y}$ is determined by $\hist_{y-1}$ and 
  \begin{equation}\label{eq:cover:EXP}
    \EXP\left(Y_{k,y}|\hist_{y-1}\right) = Z_{k,y} \cdot \PROB\big( y \AlgMap \NBH_H(v) | \hist_{y-1}\big)\,.
  \end{equation}
  Observe further that
  \begin{equation}\label{eq:cover:ZY}
    \sum_{y=1}^{t+\eps n-2}Z_{k,y}=\sum_{y=1}^{t+\eps n-2}Y_{k-1,y}\,,
  \end{equation}
  because both sums count the number of vertices $x\in X_{t,d}$ such that
  the first $k-1$ vertices of $\LNBH(x)$ are embedded to $\NBH_H(v)$, in the
  first sum grouped by their $k$-th left neighbour, and in the second sum
  by their $(k-1)$-st left neighbour.

  Assume now that~$y\in V(G)$ is fixed and that~$\hist_{y-1}$ is such that
  $\hist_{y-1}\cap\cE_k\neq\emptyset$, and
  let us
  bound $\PROB\big( y \AlgMap\NBH_H(v) | \hist_{y-1}\big)$. \corrig{Observe that if $v\in\im\psi_{y-1}$, then we are by definition in the event $\cY_k$ and hence not contributing to the probability of~\eqref{eq:cover:goal}. Thus we can assume in what follows that 
  \begin{equation}\label{eq:vout}
  v\not\in\im\psi_{y-1}\;.
  \end{equation}} 
  Since~$\hist_{y-1}\cap\cE_k\neq\emptyset$ and
  $\cD\supseteq\cE_k$, by definition of~$\cD$ the $(\beta_{t},2D+3)$-diet condition holds for
  $(H,\im\psi_{y-\eps n})$, where we have to subtract $\eps n$ in the index
  of $\psi_{y-\eps n}$ because~$y$ could be as large as~$t+\eps n-2$ (and we only
  know that the diet condition holds up to time $t-1$).
  This implies that for each set~$S$ of vertices in~$H$ with $|S|\le 2D+3$
  we have
  \begin{equation*}\begin{split}
  \big|\NBH_H(S)\setminus\im\psi_{y-1}\big|
  &=(1\pm\beta_t)p^{|S|}(n-y+\eps n)\pm\eps n \\
  &=(1\pm\beta_t)p^{|S|}(n-y+1)\pm2\eps n
  =(1\pm2\beta_t)p^{|S|}(n-y+1)\,,
  \end{split}\end{equation*}
  where the last inequality follows from $\gamma\le p$ and
  $\eps\le\alpha\gamma^{2D+3}\le\frac12\beta_t\gamma^{2D+3}$. We conclude that
  the $(2\beta_{t},2D+3)$-diet condition holds for
  $(H,\im\psi_{y-1})$.
  Since $\LEFTDEG(y)\le D$ it follows that
  \begin{align*}
    \big|\CANDSET^{y-1}(y)\setminus\im\psi_{y-1}\big|&=(1\pm2\beta_{t})p^{\LEFTDEG(y)}(n-y+1)\quad\text{and}\\
    \big|\NBH_{H}(v)\cap\CANDSET^{y-1}(y)\setminus\im\psi_{y-1}\big|&=(1\pm2\beta_{t})p^{1+\LEFTDEG(y)}(n-y+1)\,.
  \end{align*}
  \corrig{Here we used the diet condition twice, once with the set of vertices $\psi_{y-1}\big(\LNBH_G(y)\big)$ and once with the set $\psi_{y-1}\big(\LNBH_G(y)\big)\cup\{v\}$. The latter set is indeed one larger than the former since $\psi_{y-1}\big(\LNBH_G(y)\big)$ is by definition contained in the image of $\psi_{y-1}$ and $v$ is not by~\eqref{eq:vout}.}
  Therefore we have
  \begin{equation*}
    \PROB\big( y \AlgMap\NBH_H(v) | \hist_{y-1}\big)
    =\frac{\big|\NBH_{H}(v)\cap\CANDSET^{y-1}(y)\setminus\im\psi_{y-1}\big|}{\big|\CANDSET^{y-1}(y)\setminus\im\psi_{y-1}\big|}
    =(1\pm 10\beta_{t})p\,.
  \end{equation*}
  We conclude from~\eqref{eq:cover:EXP} that
  \begin{equation}\label{eq:SiM}
    \sum_{y=1}^{t+\eps n-2}\EXP(Y_{k,y}|\hist_{y-1})=(1\pm10\beta_t)p\sum_{y=1}^{t+\eps n-2}Z_{k,y}\,,
  \end{equation}
  unless~$\cE_k$ fails.
  Further, unless~$\cE_k$ fails, we have
  \[\sum_{y=1}^{t+\eps n-2}Z_{k,y}
    \eqByRef{eq:cover:ZY}
  \sum_{y=1}^{t+\eps n-2}Y_{k-1,y}
  \eqByRef{eq:eventcalY_k}
  (1\pm10\beta_{t})^{k-1}p^{k-1}|X_{t,d}|\pm
  (k-1)\varepsilon^{2}n/d\,.\]
  Plugging this in~\eqref{eq:SiM}, we get that $\cE_k$ fails or we have
  \[\sum_{y=1}^{t+\eps n-2}\EXP\left(Y_{k,y}|\hist_{y-1}\right)= (1\pm10\beta_t)^kp^k
  |X_{t,d}|\pm (k-1)\varepsilon^{2}n/d\,.
  \]
  Since $0\le Y_{k,y}\le\deg(y)$ for each $y$,
  we can thus apply Lemma~\ref{lem:concentration} with the event
  $\cE=\cE_k$, with $\mu\pm\nu=(1\pm10\beta_{t})^{k}p^{k}|X_{t,d}|\pm
  (k-1)\varepsilon^{2}n/d$, and with $\varrho=\eps^2 n/d$ to conclude that
  \[ \PROB\left(\cE_k \text{ and not } \cY_k\right) = 
  \PROB\left(\cE_k \text{ and } \sum_{y=1}^{t+\eps n-2}
  Y_{k,y}\neq\mu\pm(\nu+\rho)\right)\le 2\exp\left(-\frac{2\rho^2}{\sum_{y=1}^{t+\eps n-2}\deg(y)^2}\right)\,.\]
  By Lemma~\ref{lem:squarebound} applied to~$G$, and because $\Delta(G)\le cn/\log n$, we have
  \begin{equation*}
    \frac{2\rho^2}{\sum_{y=1}^{t+\eps n-2}\deg(y)^2} = \frac{2\eps^4n^2}{d^2\sum_{y=1}^{t+\eps n-2}\deg(y)^2}
    \ge \frac{\varepsilon^{4}\log n}{d^{2}Dc}\,,
  \end{equation*}
  and hence, because $c\le D^{-4}\eps^4/100$ and $d\le D$, we obtain~\eqref{eq:cover:goal} as desired.
\end{proof}

Finally, we prove the diet lemma.

\begin{proof}[Proof of Lemma~\ref{lem:dietcon}] First observe that if $\psi_{t-1}$ satisfies the $(\beta_{t-1},2D+3)$-diet condition, \RandomEmbedding{} cannot fail at time $t$, so $\psi_t$ exists. We first state a claim that if the diet condition holds up to time $t-\eps n$, then for any given large set $T\subseteq V(H)$, with high probability either the cover condition fails at some time before $t-\eps n$, or $\psi_{t-1}$ embeds about the expected fraction of each interval of $\eps n$ vertices to $T$.
\begin{claimDiet}\label{cl:diet:prob}
 For every $1\le j\le t-\eps n+1$, and for every $T\subseteq V(H)\setminus\im\psi_j$ with $|T|\ge\tfrac12\gamma^{2D+3}\delta n$, if the $(\beta_j,2D+3)$-diet condition holds for $(H,\im\psi_j)$, then with probability at least $1-n^{-2D-19}$, one of the following occurs.
 \begin{enumerate}[label=(\alph*)]
  \item\label{e:part1} $\psi_{t}$ does not have the $(\eps,20D\beta_j,j)$-cover condition, or
  \item\label{e:part2} $\big|\{x\::\:j\le x<j+\eps n, \psi_{t-1}(x)\in T\}\big|=(1\pm40D\beta_j)\tfrac{|T|\eps n}{n-j}$.
 \end{enumerate}
\end{claimDiet}

We defer the proof of this claim until later, and move on to state a second claim, which we will deduce from Claim~\ref{cl:diet:prob}. Let $\ell=\lfloor\tfrac{t}{\eps n}\rfloor$. We claim that either we witness a failure of the diet or cover conditions before time $t$, or the set $\NBH_H(R)\cap\NBH_{H^*}(S\setminus R)\setminus\im\psi_{\ell\eps n}$ has about the expected size for each $R\subseteq S\subseteq V(H)$ with $|S|\le 2D+3$.

\begin{claimDiet}\label{cl:diet:errors}
 With probability at least $1-n^{-10}$, one of the following holds.
 \begin{enumerate}[label=(\alph*)]
  \item The $(\beta_j,2D+3)$-diet condition fails for $(H,\im\psi_j)$ for some $1\le j\le t-1$, or
  \item the $(\eps,20D\beta_j,j)$-cover condition fails for $\psi_{t-1}$ for some $1\le j\le t+1-\eps n$, or
  \item for every $R\subset S\subseteq V(H)$ with $|S|\le 2D+3$, we have
   \begin{equation}\label{eq:diet:NHS}\begin{split}
    &\big|\NBH_H(R)\cap\NBH_{H^*}(S\setminus R)\setminus\im\psi_{\ell\eps n}\big|=\\
    &\big|\NBH_H(R)\cap\NBH_{H^*}(S\setminus R)\big|\prod_{k=0}^{\ell-1}\Big(1-\big(1\pm40D\beta_{k\eps n}\big)\tfrac{\eps n}{n-k\eps n}\Big)\,.
   \end{split}\end{equation}
 \end{enumerate}
\end{claimDiet}

Before proving these claims, we show that Claim~\ref{cl:diet:errors} implies the lemma. We want to show that~\eqref{eq:diet:NHS} holding implies that we do not have witnesses for a failure of the diet condition nor the codiet condition at time $t$. Indeed, taking logs, we have
\begin{align*}
 &\log \big|\NBH_H(R)\cap\NBH_{H^*}(S\setminus R)\setminus \im\psi_{\ell\eps n}\big|\\
 &=\log \big|\NBH_H(R)\cap\NBH_{H^*}(S\setminus R)\big|+\sum_{k=0}^{\ell-1}\log \Big(1-(1\pm 40D\beta_{k\varepsilon n})\tfrac{\varepsilon n}{n-k\varepsilon n}\Big)\\
 &=\log \big|\NBH_H(R)\cap\NBH_{H^*}(S\setminus R)\big|+\sum_{k=0}^{\ell-1}\Big(\log\tfrac{n-(k+1)\eps n}{n-k\eps n}+\log \big(1\pm \tfrac{40D\beta_{k\eps n}\eps n}{n-(k+1)\eps n}\big)\Big)\\
 &=\log \big|\NBH_H(R)\cap\NBH_{H^*}(S\setminus R)\big|+\log \big(1-\ell\eps\big) \pm 2\sum_{k=0}^{\ell-1}\tfrac{40D\beta_{k\eps n}\eps }{1-(k+1)\eps}\,,
\end{align*}
where the final equality holds since $1-(k+1)\eps\ge\delta$, and hence by choice of $\eps$ the quantity $\tfrac{40D\beta_{k\eps n}\eps}{1-(k+1)\eps}$ is close to $0$. Since at most $\eps n$ vertices are removed from $\NBH_H(R)\cap\NBH_{H^*}(S\setminus R)\setminus \im\psi_{\ell\eps n}$ to obtain $\NBH_H(R)\cap\NBH_{H^*}(S\setminus R)\setminus\im\psi_t$, we conclude
\begin{equation}\label{eq:diet:nearly}\begin{split}
 &\big|\NBH_H(R)\cap\NBH_{H^*}(S\setminus R)\setminus\im\psi_t\big|\\
 &=|\NBH_H(R)\cap\NBH_{H^*}(S\setminus R)|\cdot\frac{n-t\pm\eps n}{n}\cdot \exp\Big(\pm 80D\delta^{-1}\eps\sum_{k=0}^{\ell-1}\beta_{k\eps n}\Big)\pm\eps n\;.
\end{split}\end{equation}

We first consider the case $R=S$, when $\NBH_H(R)\cap\NBH_{H^*}(S\setminus R)=\NBH_H(S)$, and deduce that $S$ does not witness a failure of the $(\beta_t,2D+3)$-diet condition for $(H,\im\psi_t)$. Indeed, from~\eqref{eq:diet:nearly} we have
\begin{align*}
 \big|\NBH_H(S)\setminus\im\psi_t\big|&=
 |\NBH_H(S)|\cdot\frac{n-t\pm\eps n}{n}\cdot \exp\Big(\pm 80D\delta^{-1}\eps\sum_{k=0}^{\ell-1}\beta_{k\eps n}\Big)\pm\eps n\\
 &=(1\pm \alpha)p^{|S|}(n-t\pm \eps n)\Big(1\pm 200D\delta^{-1}\eps\sum_{k=0}^{\ell-1}\beta_{k\eps n}\Big)\Big (1\pm \frac{2\eps n}{p^{|S|}(n-t)}\Big)
\end{align*}
where the second equality uses the fact that $H$ is ($\alpha, 2D+3$)-quasirandom. We thus have
\begin{align*}
 \big|\NBH_H(S)\setminus\im\psi_t\big|&=(1\pm\alpha)p^{|S|}(n-t)\Big(1\pm 200D\delta^{-1}\eps\sum_{k=0}^{\ell-1}\beta_{k\eps n}\Big)(1\pm 4\eps\delta^{-1}\gamma^{-|S|})\\
 &\!\!\eqByRef{eq:betabound}(1\pm \alpha)p^{|S|}(n-t)(1\pm\beta_t/4)(1\pm 4\varepsilon \delta^{-1}\gamma^{-|S|})\\
 &=(1\pm\beta_t)p^{|S|}(n-t)\,.
\end{align*}

Now, we let $R$ be any subset of $S$ and aim to establish the codiet condition. Again from~\eqref{eq:diet:nearly}, we have
\begin{align*}
 \big|\NBH_H(R)\cap & \NBH_{H^*}(S\setminus R)\setminus\im\psi_t\big| \\
 &=|\NBH_H(R)\cap\NBH_{H^*}(S\setminus R)|\cdot\frac{n-t\pm\eps n}{n}\cdot \exp\Big(\pm 80D\delta^{-1}\eps\sum_{k=0}^{\ell-1}\beta_{k\eps n}\Big)\pm\eps n\\
&=(1\pm \eta)p^{|R|}\hat p^{|S\setminus R|}(n-t\pm \eps n)\Big(1\pm 200D\delta^{-1}\eps\sum_{k=0}^{\ell-1}\beta_{k\eps n}\Big)\Big (1\pm \frac{2\eps n}{p^{|S|}(n-t)}\Big)
\end{align*}
since $(H,H^*)$ is ($\eta, 2D+3$)-coquasirandom. Therefore
\begin{align*}
 \big|\NBH_H(R)\cap & \NBH_{H^*}(S\setminus R)\setminus\im\psi_t\big| \\
 &=(1\pm \eta)p^{|R|}\hat p^{|S\setminus R|}(n-t)\Big(1\pm 200D\delta^{-1}\eps\sum_{k=0}^{\ell-1}\beta_{k\eps n}\Big)(1\pm 4\eps\delta^{-1}\gamma^{-|S|})\\
 &\!\!\eqByRef{eq:betabound}(1\pm \eta)(1\pm \beta_k)(1\pm 4\varepsilon \delta^{-1}\gamma^{-|S|})p^{|R|}\hat p^{|S\setminus R|}(n-t)\\
 &=(1\pm 2\eta)p^{|R|}\hat p^{|S\setminus R|}(n-t)\,.
\end{align*}

This concludes the proof of the lemma, modulo the proofs of Claim~\ref{cl:diet:prob} and Claim~\ref{cl:diet:errors}, which we now provide.

\begin{claimproof}[Proof of Claim~\ref{cl:diet:prob}]
Let $j$ and $T$ be as in the statement. Fix $0\le d\le D$. We want to show how to make use of the $(\varepsilon,20D\beta_{j},j)$-cover condition for \corrig{$\psi_{t}$} (which we have when Part~\ref{e:part1} fails) to deduce that the assertion of Part~\ref{e:part2} holds with high probability. That is, we consider the number
of vertices in $X_{j,d}$ embedded to $T$. In order to apply Lemma~\ref{lem:concentration}, we want to estimate the sum over $x\in X_{j,d}$
of the probability that $x$ is embedded to $T$, conditioning on
$\psi_{x-1}$, that is, we need to estimate the number 
\begin{equation}\label{eq:divadlo1}
\frac{\big|T\cap\CANDSET^{x-1}(x)\setminus\im\psi_{x-1}\big|}{\big|\CANDSET^{x-1}(x)\setminus\im\psi_{x-1}\big|}\;.
\end{equation}
By the diet condition, we have $\big|\CANDSET^{x-1}(x)\setminus\im\psi_{j}\big|=(1\pm\beta_{j})p^{d}(n-j)$.
Since $j< t\le(1-\delta)n$, since $x\le j+\varepsilon n$, since
$p\ge\gamma,$ and by choice of $\varepsilon$, we have 
\begin{equation}\label{eq:divadlode}
\big|\CANDSET^{x-1}(x)\setminus\im\psi_{x-1}\big|=(1\pm2\beta_{j})p^{d}(n-j)\;,
\end{equation} thus providing a bound on the denumerator in~\eqref{eq:divadlo1}. (Note that this bound on the denumerator does not depend on the choice of $x\in X_{j,d}$.)
Now $x$ is embedded uniformly at random into $\CANDSET^{x-1}(x)\setminus\im\psi_{x-1}$,
so it remains to determine the sum of the numerators in~\eqref{eq:divadlo1},
\begin{align}
\nonumber
 \sum_{x\in X_{j,d}}\big|T\cap\CANDSET^{x-1}(x)\setminus\im\psi_{x-1}\big|&=\sum_{x\in X_{j,d}}\big|T\cap\CANDSET^{x-1}(x)\setminus\im\psi_{j}\big|\pm\varepsilon|X_{j,d}|n\\
\label{eq:divadlonu1}
 &=\sum_{x\in X_{j,d}}\big|T\cap\CANDSET^{x-1}(x)\big|\pm\varepsilon^{2}n^{2}\,,
\end{align}
where the first equality uses $j\le x<j+\varepsilon n$, and the second
the fact that $T\subseteq V(H)\setminus\im\psi_{j}$ and that $|X_{j,d}|\le\varepsilon n$.

\corrig{Consider a vertex $v\in T$. If $v\not\in\im\psi_{j+\eps n-1}$, then the $(\eps,20D\beta_j,j)$-cover condition tells us that $v$ contributes $(1\pm20D\beta_{j})p^{d}|X_{j,d}|\pm\varepsilon^{2}n$ to the summation $\sum_{x\in X_{j,d}}\big|T\cap\CANDSET^{x-1}(x)\big|$. Since $T\subseteq V(H)\setminus\im\psi_j$, there are at most $\eps n$ vertices $v\in T$ such that $v\in\im\psi_{j+\eps n-1}$, and these contribute between $0$ and $\eps n$ to the summation; in particular they contribute $(1\pm20D\beta_{j})p^{d}|X_{j,d}|\pm\varepsilon^{2}n\pm\eps n$ to the summation. Putting this together, we have
\[\sum_{x\in X_{j,d}}\big|T\cap\CANDSET^{x-1}(x)\big|=|T|(1\pm20D\beta_{j})p^{d}|X_{j,d}|\pm\varepsilon^{2}|T|n\pm\eps^2n^2\,.\]
Putting this into~\eqref{eq:divadlonu1} we have
\begin{equation}\label{eq:divadlonufi}
\sum_{x\in X_{j,d}}\big|T\cap\CANDSET^{x-1}(x)\setminus\im\psi_{x-1}\big|=(1\pm20D\beta_{j})p^{d}|T||X_{j,d}|\pm3\varepsilon^{2}n^{2}\;.
\end{equation}}
%But now if the $(\varepsilon,20D\beta_{j},j)$-cover condition holds
%for $\psi_{j}$, then summing over $v\in T$ we obtain
%\[
%\sum_{x\in X_{j,d}}\big|T\cap\CANDSET^{x-1}(x)\big|=|T|(1\pm20D\beta_{j})p^{d}|X_{j,d}|\pm\varepsilon^{2}|T|n\,,
%\]
%which, combined with~\eqref{eq:divadlonu1}, gives
%\begin{equation}\label{eq:divadlonufi}
%\sum_{x\in X_{j,d}}\big|T\cap\CANDSET^{x-1}(x)\setminus\im\psi_{x-1}\big|=(1\pm20D\beta_{j})p^{d}|T||X_{j,d}|\pm2\varepsilon^{2}n^{2}\;.
%\end{equation}

We can thus apply Lemma~\ref{lem:concentration}, setting $\cE$ to
be the event that the $(\varepsilon,20 D\beta_{j},j)$-cover condition
holds for $\psi_{j}$. The random variables whose sum we are estimating are the Bernoulli random variables indicating whether each $x\in X_{j,d}$ is embedded to $T$, so the sum of squares of their ranges is at most $\varepsilon n$. Combining~\eqref{eq:divadlode} and~\eqref{eq:divadlonufi}, the expected number of vertices of $X_{j,d}$
embedded to $T$ is
\corrig{\[
\frac{(1\pm20D\beta_{j})p^{d}|T||X_{j,d}|\pm3\varepsilon^{2}n^{2}}{(1\pm2\beta_{j})p^{d}(n-j)}=(1\pm30D\beta_{j})\frac{|T||X_{j,d}|}{n-j}\pm4\varepsilon^{2}\gamma^{-d}\delta^{-1}n\,,
\]}
%\[
%\frac{(1\pm20D\beta_{j})p^{d}|T||X_{j,d}|\pm2\varepsilon^{2}n^{2}}{(1\pm2\beta_{j})p^{d}(n-j)}=(1\pm30D\beta_{j})\frac{|T||X_{j,d}|}{n-j}\pm4\varepsilon^{2}\gamma^{-d}\delta^{-1}n\,,
%\]
where we use $n-j\ge\delta n$ and $p\ge\gamma$. The probability
that the $(\varepsilon,20D\beta_{j},j)$-cover condition holds for $\psi_{j}$
and the outcome differs from this by more than $\varepsilon^{2}n$
is at most $2\exp(-2\eps^3n)\le n^{-2D-20}$, so taking the union bound over the $D+1$
choices of $d$ and summing, we conclude that with probability at
most $n^{-2D-19}$ the $(\varepsilon,20 D\beta_{j},j)$-cover condition
holds for $\psi_{j}$ and the number of vertices $x$ with $j\le x<j+\varepsilon n$
embedded to~$T$ is not equal to
\[
(1\pm30D\beta_{j})\frac{|T|\varepsilon n}{n-j}\pm4(D+1)\varepsilon^{2}\gamma^{-D}\delta^{-1}n\pm(D+1)\varepsilon^{2}n=(1\pm40D\beta_{j})\frac{|T|\varepsilon n}{n-j}\,,
\]
where the final equality uses our lower bound on $|T|$ and the choice
of $\varepsilon$. This is what we wanted to show.
\end{claimproof}

\begin{claimproof}[Proof of Claim~\ref{cl:diet:errors}]
Given a set $S\subseteq V(H)$ with $|S|\le 2D+3$ and a subset $R\subseteq S$, for each integer $0\le k< \ell$, we set $T_k=\NBH_H(R)\cap \NBH_{H^*}(S\setminus R)\setminus\im\psi_{k\varepsilon n}$. Observe that as $(H,H^*)$ is  $(\eta, 2D+3)$-coquasirandom, we have
\[|T_0|\ge (1-\eta)p^{|R|}\hat p^{|S\setminus R|}n\ge (1-\eta)^{2D+4}\gamma^{2D+3}n\,.\]

For each $0\le k<\ell$, suppose that
\begin{align*}
|T_k|&\ge (1-\eta)^{2D+4}(1-80D\beta_n\delta^{-1}\eps)^k\gamma^{2D+3}(n-k\epsilon n)\\
&\ge \tfrac{9}{10}(1-80D\beta_n\delta^{-1}\eps)^{1/\eps}\gamma^{2D+3}\delta n\ge \tfrac{9}{10}\exp\big(-200D\beta_n\delta^{-1}\big)\gamma^{2D+3}\delta n\\
&>\frac 12\gamma^{2D+3}\delta n\,,
\end{align*}
where the final line follows since $200D\beta_n\delta^{-1}\le 400CD\alpha\delta^{-1}<1/100$ by choice of $\alpha$. We can thus apply Claim~\ref{cl:diet:prob}   with $T=T_{k}$ and obtain that with probability at least $1-n^{-2D-19}$ either we have a failure of the diet or the cover condition is witnessed before time~$k$, or we have 
\[|T_{k+1}|= |T_k|\Big(1-\big(1\pm 40D\beta_{k\eps n}\big)\frac{\eps n}{n-k\eps n}\Big)\;.\] 
Observe that then 
\begin{align*}|T_{k+1}|&\ge |T_k|\left(1 -\frac{\eps n}{n-k\eps n}-40D\beta_{n}\delta^{-1}\eps \right)\\
	&> (1-\eta)^{2D+4}(1-80D\beta_n\delta^{-1}\eps)^{k+1}\gamma^{2D+3}\big(n-(k+1)\epsilon n\big)\;,
\end{align*}
providing the assumption for using of Claim~\ref{cl:diet:prob} in step $k+1$.

Repeating this process for each $0\le k\le \ell-1$ we get that with probability at least $1-\varepsilon^{-1} n^{-2D-19}$ either a failure of the diet or cover condition is witnessed before time $\ell \varepsilon n$, or we have
\[\big|T_{\ell}\big|=\Big|\NBH_H(R)\cap \NBH_{H^*}(S\setminus R)\Big|\prod_{k=0}^{\ell-1}\Big(1-(1\pm 40D\beta_{k\varepsilon n})\tfrac{\varepsilon n}{n-k\varepsilon n}\Big)\,.\]

Taking a union bound over the at most $(2D+3)n^{2D+3}$ choices of $S$ and the at most~$2^{2D+3}$ choices of $R\subseteq S$, we see that with probability at least $1-n^{-10}$ either a failure of the diet or cover condition is witnessed before time $t$, or the above equation holds for all $|S|\le 2D+3$ and $R\subseteq S$.
\end{claimproof} 
\end{proof}

\section{Maintaining quasirandomness}
\label{sec:maintainingquasirandomness}

In this section we provide the proofs of Lemma~\ref{lem:firstquasirandomness} and Lemma~\ref{lem:aggregate}.

\subsection{Initial coquasirandomness}

We begin with the easy proof of Lemma~\ref{lem:firstquasirandomness}, which states that splitting the edges of a quasirandom graph randomly gives a coquasirandom pair with high probability.

\begin{proof}[Proof of Lemma~\ref{lem:firstquasirandomness}]
Using~\eqref{eq:hoeffding} we see that the densities $p_0$ and $p_0^*$ of $H_0$ and $H_0^*$ satisfy
\begin{equation}\label{eq:p0p0}
p_0=(1\pm \tfrac{\alpha_0}{1000D})(p-\gamma)
\quad\mbox{and}\quad
p^*_0=(1\pm\tfrac{\alpha_0}{1000D})\gamma
\end{equation}
with probability at least $1-n^{-10}$, giving the first part of Lemma~\ref{lem:firstquasirandomness}.

Now, let $R\subset S\subset V(\widehat H)$ be two sets of size at most $2D+3$. By quasirandomness of $\widehat H$ we have $|\NBH_{\widehat H}(S)|=(1\pm \xi)p^{|S|}n$. Observe that each vertex of $\NBH_{\widehat H}(S)$ appears with probability $q^{|R|}(1-q)^{|S\setminus R|}$ in $\NBH_{H_0^*}(R)\cap \NBH_{H_0}(S\setminus R)$. Hence, 
$$\EXP\left(\left|\NBH_{H_0^*}(R)\cap \NBH_{H_0}(S\setminus R)\right|\right)
=
q^{|R|}(1-q)^{|S\setminus R|}(1\pm \xi)p^{|S|}n\;.
$$
Observe also that for distinct vertices in $\NBH_{\widehat H}(S)$ the events whether these appear in $\NBH_{H_0^*}(R)\cap \NBH_{H_0}(S\setminus R)$ are independent. Using again~\eqref{eq:hoeffding}, with probability at least $1-n^{-2D-10}$ we have that
\begin{equation}\label{eq:takemetotheunion}
\left|\NBH_{H_0^*}(R)\cap \NBH_{H_0}(S\setminus R)\right|
=
q^{|R|}(1-q)^{|S\setminus R|}(1\pm 2\xi)p^{|S|}n
\;.
\end{equation}

Taking the union bound we conclude that~\eqref{eq:takemetotheunion} holds for all $S\subset V(\widehat H)$ with $|S|\le 2D+3$ and $R\subset S$ with probability at least $1-n^{-6}$.

Now, assume that~\eqref{eq:p0p0} holds. Then the right-hand side of~\eqref{eq:takemetotheunion} can be rewritten as
\begin{align*}
(1\pm 2\xi)\gamma^{|R|}(p-\gamma)^{|S\setminus R|}n
&=
(1\pm 2\xi)\left(\tfrac{p_0^*}{1\pm \tfrac{\xi_0}{1000D}}\right)^{|R|}\left(\tfrac{p_0}{1\pm \tfrac{\alpha_0}{1000D}}\right)^{|S\setminus R|}n\\
&=(1\pm 2\xi)(1\pm \tfrac{\alpha_0}{100})(p_0^*)^{|R|}p_0^{|S\setminus R|}=\big(1\pm \tfrac1{10}\alpha_0\big)(p_0^*)^{|R|}p_0^{|S\setminus R|}\;.
\end{align*}
We conclude that $(H_0^*,H_0)$ is $\big(\tfrac1{10}\alpha_0,2D+3)$-coquasirandom with probability at least $1-n^{-5}$.
\end{proof}

\subsection{Maintaining coquasirandomness}
In this subsection we prove Lemma~\ref{lem:aggregate}. We need to show that, provided coquasirandomness is maintained up to stage $s-1$ and \RandomEmbedding\ does not fail, it is likely that coquasirandomness holds after stage $s$, when $G_s$ is embedded into $H_{s-1}$ and we obtain $H_s$. Let us briefly sketch the idea (for convenience focusing only on quasirandomness of $H_s$). We fix a set $R\subset V(\widehat{H})$ with $|R|\le 2D+3$, and consider the running of \PackingProcess{} up to
stage~$s$. We want to show that it is very unlikely that $R$ witnesses the failure of $H_s$ to be quasirandom, since then the union bound over choices of $R$ tells us that it is likely that $H_s$ is quasirandom. In other words, we want to know that $\big|N_{H_s}(R)\big|$ is very likely close to the expected size. We write
\[\big|N_{H_s}(R)\big|=\big|N_{H_0}(R)\big|-Y_1-\dots-Y_s\,,\]
where $Y_i=\big|N_{H_i-1}(R)\big|-\big|N_{H_{i}}(R)\big|$ is the change at step $i$, and apply Lemma~\ref{lem:freedman} to show that the sum $Y_1+\dots+Y_s$ is very likely to be close to its expectation. So proving Lemma~\ref{lem:aggregate} boils down to estimating accurately $\EXP(Y_i|H_{i-1})$ and finding a reasonable upper bound for $\EXP(Y_i^2|H_{i-1})$. The latter turns out to be relatively straightforward and is done in Lemma~\ref{lem:starsecond}. We now outline the route to the former estimation.

Observe that $Y_i$ is equal to the number of stars in $H_{i-1}$ whose leaves are the vertices in $R$ and at least one of whose edges is used in embedding $G_i$ to $H_{i-1}$. By linearity of expectation, $\EXP(Y_i|H_{i-1})$ is equal to the sum, over stars in $H_{i-1}$ whose leaves are $R$, of the probability that at least one edge in the star is used in embedding $G_i$. We will see that this probability is about the same for any given  star $S$, and the problem is to calculate it. To do this we need to consider the running of \RandomEmbedding.

We begin in Lemma~\ref{lem:probint} by estimating the chance that a given vertex, or one of a given pair of vertices, is used in a short time interval in \RandomEmbedding. From this we deduce in Lemma~\ref{lem:probvtx} the probability that a given vertex, or one of a pair, is used in any given time interval. This helps us to establish, in Lemma~\ref{lem:probedge}, that any given edge of $H_{i-1}$ is about equally likely to be used in the embedding of $G_i$. Finally, in Lemma~\ref{lem:probstar} we show that the chance of two or more edges in $S$ being used in the embedding of $G_i$ is tiny, from which it follows that the chance of one or more is about $|R|$ times the probability of any given edge being used.

All of these estimations depend upon $H_{i-1}$ being sufficiently quasirandom, and the errors depend upon the quasirandomness $\alpha_{i-1}$. Because the errors add up over time, it is important that the $\alpha_{s}$ increase
quite fast with~$s$. Here it is very important that the
dependence of the error term in Lemma~\ref{lem:diet} is linear in
the input~$\alpha$ and not much worse: otherwise it would not be
possible to choose any sequence $\alpha_{s}$ such that the error
remains bounded by $\alpha_{s}$ at each stage $s$.

As the main work
is to estimate the probability that, for a given $H_{s-1}$ and $G_s$, and $R$ and~$v$, \RandomEmbedding\ uses an edge of the star with centre $v$ and leaves $R$ when embedding $G_s$ into $H_{s-1}$, for most of this section we will consider fixed graphs $G$ and $H$. We now embark upon this probability estimation.

First, for given $u,v\in V(H)$, we estimate the probability that \RandomEmbedding{} embeds a vertex to $\{u,v\}$
in the short interval of time $[t,t+\varepsilon n)$, conditioning on not
having done so before time $t$, and the probability that \RandomEmbedding{}
embeds a vertex to $v$ in the interval of time $[t,t+\varepsilon n),$
conditioning on not having done so before time $t$. In both cases, we need to assume that the history $\hist_{t-1}$ of embedding up to time $t-1$ is typical (in a sense which we now make precise).
\begin{lem}\label{lem:probint}
Given~$D\in \mathbb{N}$ and $\gamma>0$, let $\delta,\alpha_0,\alpha_{2n},C,\eps$ be as in Setting~\ref{set:graphs}. The following holds for any $\alpha_0\le\alpha\le\alpha_{2n}$ and all sufficiently large $n$. Suppose
that $G$ is a graph on $[n]$ such that $\LEFTDEG(x)\le D$ for each $x\in V(G)$, and~$H$ is an $(\alpha,2D+3)$-quasirandom graph
with~$n$ vertices and $p\binom{n}{2}$ edges, with $p\ge\gamma$. Suppose that $u$ and $v$ are two distinct vertices of $H$. When \RandomEmbedding{} is run to embed~$G[{\scriptstyle [n-\delta n]}]$ into~$H$, for any $1\le t\le n+1-(\delta+\eps) n$ we have the following two statements.
\begin{enumerate}[label=(\alph*)]
\item\label{e:probintv} Suppose the history $\hist_{t-1}$ up to and including embedding $t-1$ is such that $v\not\in\im\psi_{t-1}$, the $(C\alpha,2D+3)$-diet condition holds for $(H,\im\psi_{t-1})$, and 
\begin{equation}\label{eq:breakfast}
\PROB^{G\AlgMap H}\left(
\CoverEvent(C\alpha;t)|\hist_{t-1}\right)
\le n^{-3}\,.
\end{equation}
Then we have
\begin{equation*}
\PROB^{G\AlgMap H}\big(v\in\im\psi_{t+\varepsilon n-1}\big|\hist_{t-1}\big)=(1\pm 10C\alpha)\tfrac{\eps n}{n-t}\,.
\end{equation*}
\item\label{e:probintuv}
Suppose the history $\hist_{t-1}$ up to and including embedding $t-1$ is such that $u,v\not\in\im\psi_{t-1}$, the $(C\alpha,2D+3)$-diet condition holds for $(H,\im\psi_{t-1})$, and 
\[
\PROB^{G\AlgMap H}\left(
\CoverEvent(C\alpha;t)|\hist_{t-1}\right)
\le n^{-3}\,.\]
Then we have
\begin{equation*}
\PROB^{G\AlgMap H}\Big(\big|\{u,v\}\cap\im\psi_{t+\varepsilon n-1}\big|\ge1\Big|\hist_{t-1}\Big)=(1\pm 10C\alpha)\tfrac{2\eps n}{n-t}\,.
\end{equation*}
\end{enumerate}
\addtocounter{counterProbint}{\value{thm}}
\end{lem}
Before proving Lemma~\ref{lem:probint}, we first sketch its proof. For
Lemma~\ref{lem:probint}\ref{e:probintv}, the idea is that either the cover
condition fails, or $v$ is in candidate sets of roughly $p^d|X_{t,d}|$ vertices
$x$ of $X_{t,d}$ (for each $d$). Because the diet condition holds at time $t-1$,
each of these vertices $x$ is embedded uniformly at random to a set of roughly
$p^d(n-t)$ vertices. One would like to say that it follows that the probability
that $x$ is embedded to $v$ is thus about $1/(p^d(n-t))$ and the desired result
follows by summing these probabilities. Unfortunately this is not true: the
probability that $x$ is embedded to $v$ also depends on the probability that no
previous vertex was embedded to $v$. In order to get around this, we define the following  \ModRandomEmbedding, which generates a sequence of embeddings with an identical distribution to \RandomEmbedding, but which in addition generates a sequence of \emph{reported} vertices. The modification we make is simple: at each time $1\le t'\le n-\delta n$, \RandomEmbedding\ chooses a vertex of $\CANDSET_{G\AlgMap H}^{t'-1}(t')\setminus\im\psi_{t'-1}$. In \ModRandomEmbedding, we instead choose a vertex $w$ of $\CANDSET_{G\AlgMap H}^{t'-1}(t')\setminus(\im\psi_{t'-1}\setminus\{v\})$, and report this vertex. If the reported vertex $w$ is not in $\im\psi_{t'-1}$, we set $\psi_{t'}=\psi_{t'-1}\cup\{t'\AlgMap w\}$, as in \RandomEmbedding. If the reported vertex is in $\im\psi_{t'-1}$ (which happens only if $w=v$) we choose $w'$ uniformly at random in $\CANDSET_{G\AlgMap H}^{t'-1}(t')\setminus\im\psi_{t'-1}$, and set $\psi_{t'}=\psi_{t'-1}\cup\{t'\AlgMap w'\}$. We will see that it is easy to calculate the expected number of times $v$ is reported, and also easy to show that the contribution due to $v$ being reported multiple times is tiny. The point is that the probability of \RandomEmbedding\ using $v$ is the same as the probability that \ModRandomEmbedding\ reports $v$ at least once, which we can thus calculate.

Lemma~\ref{lem:probint}\ref{e:probintuv} is established similarly, using a slightly different version of \ModRandomEmbedding.

\begin{proof}[Proof of Lemma~\ref{lem:probint}\ref{e:probintv}]
Instead of \RandomEmbedding, we consider \ModRandomEmbedding\ as defined above, which creates the same embedding distribution.
For each $i$, let $r(i)$ be the vertex reported by  \ModRandomEmbedding\ at time $i$. 
 We shall use the following two auxiliary claims.
 
 Define $E$ as the random variable counting the times when $v$ is reported by \ModRandomEmbedding\ in the interval $t\le x<t+\eps n$,
\[E=\left|\:\big\{x\in [t,t+\eps n) \::\: r(x-1)=v\big\}\:\right|\;.\]
The probability that \RandomEmbedding\ uses $v$ in the interval $t\le x<t+\eps n$, conditioning on $\hist_{t-1}$, is equal to the probability that  \ModRandomEmbedding\ reports $v$ at least once in that interval, which probability is by definition at most $\EXP\left(E\:|\:\hist_{t-1}\right)$ and at least
\[\EXP\left(E\:|\:\hist_{t-1}\right)-\sum_{k=2}^{\eps n}\PROB\big(\text{$v$ is reported at least $k$ times in the interval $[t,t+\eps n)$}\big|\hist_{t-1}\big)\,.\]
Our first claim estimates $\EXP\left(E\:|\:\hist_{t-1}\right)$.
 \begin{claimProbint}\label{cl:probint:E}
 	We have that
 	\[\EXP\left(E\:|\:\hist_{t-1}\right)=(1\pm 4C\alpha)\frac{\eps n}{n-t}\pm \corrig{8(D+1)\eps^2\gamma^{-2D}\delta^{-2}}\,.\] %4(D+1)\eps^2\gamma^{-D}\delta^{-1}\;.\]
 \end{claimProbint}

Our second claim is that the sum in the expression above is small.
 
 \begin{claimProbint}\label{cl:probint:sum}
 	We have that
 	\[\sum_{k=2}^{\eps n}\PROB\Big(\big|\big\{x\in [t,t+\eps n)\,:\; r(x-1)=v\big\}\big|\ge k\Big)\le 8\eps^2\gamma^{-2D}\delta^{-2}\;.\]
 \end{claimProbint}
By choice of $\eps$, we have \corrig{$32(D+1)\eps^2\gamma^{-2D}\delta^{-2}<C\alpha\eps\delta$}.
 %$16(D+1)\eps^2\gamma^{-2D}\delta^{-2}<C\alpha\eps$
 Thus the two claims  give Lemma~\ref{lem:probint}\ref{e:probintv}. We now prove the auxiliary Claims~\ref{cl:probint:E} and~\ref{cl:probint:sum}.

\begin{claimproof}[Proof of Claim~\ref{cl:probint:E}]
	Note that since the $(C\alpha,2D+3)$-diet condition holds for $(H,\im\psi_{t-1})$, for each $t\le x<t+\eps n$, setting $S=\psi_{x-1}(\LNBH(x))$,  we have\footnote{\label{footnote:biggererror}We remark that in~\eqref{eq:probint:setsize}, the calculations are included with an error ``$\pm 2$'' and for this proof ``$\pm 1$'' would have sufficed. We reuse~\eqref{eq:probint:setsize} in the proof of Lemma~\ref{lem:probint}\ref{e:probintuv} where the bigger error is needed.} 
	\begin{equation}\label{eq:probint:setsize}\begin{split}
	\big|\CANDSET^{x-1}(x)\setminus\im\psi_{x-1}\big|\pm 2&=\big|\NBH_H(S)\setminus \im\psi_{t-1}\big|\pm \eps n\pm 2\\
	&=(1\pm C\alpha)p^{|\LNBH(x)|}(n-t)\pm\eps n\pm 2\\
	&=(1\pm 2C\alpha)p^{|\LNBH(x)|}(n-t)\,.
	\end{split}\end{equation}
	
    \corrig{Before trying to obtain an accurate estimate for $\EXP\left(E\:|\:\hist_{t-1}\right)$, we give an easy (and not very sharp) upper bound on $\PROB^{G\AlgMap H}\big(v\in\im\psi_{t+\varepsilon n-1}\big|\hist_{t-1}\big)$. When we embed any one vertex $x$ with $t\le x\le t+\eps n-1$, we have a probability at most $|\CANDSET^{x-1}(x)\setminus\im\psi_{x-1}\big|^{-1}$ of embedding $x$ to $v$ (in fact, the probability is equal either to this figure or to zero). Using the lower bound~\eqref{eq:probint:setsize} and summing over the $\eps n$ choices of $x$, we see
\begin{equation}\label{eq:lunch}
\PROB^{G\AlgMap H}\big(v\in\im\psi_{t+\varepsilon n-1}\big|\hist_{t-1}\big)\le\frac{\eps n}{\tfrac12p^D\delta n}\le 2\eps \gamma^{-D}\delta^{-1}\,.
\end{equation}
    We now try to estimate the desired expectation.} 	By linearity of expectation, we have
	\begin{align}
	\begin{split}\label{eq:obed}
	\EXP\Big(E\:|\:\hist_{t-1}\Big)&=\sum_{x=t}^{t+\eps n-1}\PROB\big(v\text{ is reported at time }x\big|\hist_{t-1}\big)\\
	&=\sum_{x=t}^{t+\eps n-1}\EXP\left(\frac{\mathbbm{1}\{v\in\CANDSET^{x-1}(x)\}}{|\CANDSET^{x-1}(x)\setminus(\im\psi_{x-1}\setminus\{v\})|}\:\Big|\hist_{t-1}\right)\\
	&=\sum_{x=t}^{t+\eps n-1}\EXP\left(\frac{\mathbbm{1}\{v\in\CANDSET^{x-1}(x)\}}{|\CANDSET^{x-1}(x)\setminus \im\psi_{x-1}|\pm 1}\:\Big|\hist_{t-1}\right)\;.
	\end{split}
	\end{align}
	Using~\eqref{eq:probint:setsize}, we get
	\[\EXP\left(E\:|\:\hist_{t-1}\right)=\sum_{x=t}^{t+\eps n-1}\frac{\PROB\big(v\in\CANDSET^{x-1}(x)\big|\hist_{t-1}\big)}{(1\pm 2C\alpha)p^{|\LNBH(x)|}(n-t)}\,.\]
	Splitting this sum up according to $|\LNBH(x)|$, and again using linearity of expectation, we have
\begin{equation}\label{eq:dinner}
\EXP\left(E\:|\:\hist_{t-1}\right)=\sum_{d=0}^D\frac{\EXP\left(|\{x\in X_{t,d}:v\in\CANDSET^{x-1}(x)\}|\big|\hist_{t-1}\right)}{(1\pm 2C\alpha)p^{d}(n-t)}\,.
\end{equation}
\corrig{Now fix $0\le d\le D$. If the $(\eps,C\alpha,t)$-cover condition holds, and if $v\not\in\im\psi_{t+\eps n-1}$, we have $|\{x\in X_{t,d}:v\in\CANDSET^{x-1}(x)\}|=(1\pm C\alpha)p^d|X_{t,d}|\pm\eps^2 n$. If the $(\eps,C\alpha,t)$-cover condition fails, or if $v\in\im\psi_{t+\eps n-1}$ (which occur with total probability at most $n^{-3}+2\eps \gamma^{-D}\delta^{-1}$, see~\eqref{eq:breakfast} and~\eqref{eq:lunch}), we have $0\le |\{x\in X_{t,d}:v\in\CANDSET^{x-1}(x)\}|\le \eps n$. In particular in such a case we have
\[|\{x\in X_{t,d}:v\in\CANDSET^{x-1}(x)\}|=(1\pm C\alpha)p^d|X_{t,d}|\pm\eps^2 n\pm\eps n\,.\]
Putting these together, we have
\begin{align*}
		\EXP\left(|\{x\in X_{t,d}:v\in\CANDSET^{x-1}(x)\}|\;\big|\;\hist_{t-1}\right)&=(1\pm C\alpha)p^d|X_{t,d}|\pm\eps^2 n\pm \big(n^{-3}+2\eps \gamma^{-D}\delta^{-1}\big)\cdot\eps n\\
		&=(1\pm C\alpha)p^d|X_{t,d}|\pm 4\gamma^{-D}\delta^{-1}\eps^2n \,.
	\end{align*}
	Substituting this in~\eqref{eq:dinner}, we have
	\[\EXP\left(E\:|\:\hist_{t-1}\right)=\sum_{d=0}^D\frac{(1\pm C\alpha)p^d|X_{t,d}|\pm 4\gamma^{-D}\delta^{-1}\eps^2n}{(1\pm 2C\alpha)p^{d}(n-t)}
	\eqByRef{eq:SumXid}(1\pm 4C\alpha)\tfrac{\eps n}{n-t}\pm 8(D+1)\eps^2\gamma^{-2D}\delta^{-2}\,,\]
	where the last equality uses $p\ge\gamma$ and $n-t\ge\delta n$.}
%
%
%
%
%
%	Now for each $0\le d\le D$, since the $(\eps,C\alpha,t)$-cover condition holds with probability at least $1-n^{-3}$ conditioning on $\hist_{t-1}$, we have
%	\begin{align*}
%		\EXP\left(|\{x\in X_{t,d}:v\in\CANDSET^{x-1}(x)\}|\big|\hist_{t-1}\right)&=(1-n^{-3})\big((1\pm C\alpha)p^d|X_{t,d}|\pm\eps^2 n\big)\pm n^{-3}\cdot\eps n\\
%		&=(1\pm C\alpha)p^d|X_{t,d}|\pm 2\eps^2n\,.
%	\end{align*}
%	Substituting this in, we have
%	\[\EXP\left(E\:|\:\hist_{t-1}\right)=\sum_{d=0}^D\frac{(1\pm C\alpha)p^d|X_{t,d}|\pm 2\eps^2n}{(1\pm 2C\alpha)p^{d}(n-t)}=(1\pm 4C\alpha)\tfrac{\eps n}{n-t}\pm 4(D+1)\eps^2\gamma^{-D}\delta^{-1}\,,\]
%	where the last equality uses $p\ge\gamma$ and $n-t\ge\delta n$.
\end{claimproof} 
 
\begin{claimproof}[Proof of Claim~\ref{cl:probint:sum}]
 Since the $(C\alpha,2D+3)$-diet condition holds for $(H,\im\psi_{t-1})$, since $p\ge\gamma$, and since $n-t\ge\delta n$, for each $x\in[t,t+\eps n)$, when we embed $x$ we report a uniform random vertex from a set of size at least $\tfrac12\gamma^D\delta n$. The probability of reporting $v$ when we embed $x$ is thus at most $2\gamma^{-D}\delta^{-1}n^{-1}$, conditioning on $\hist_{t-1}$ and any embedding of the vertices $[t,x)$. Since the conditional probabilities multiply, the probability that at each of a given $k$-set of vertices in $[t,t+\eps n)$ we report $v$ is at most $2^k\gamma^{-kD}\delta^{-k}n^{-k}$. Taking the union bound over choices of $k$-sets, we have
 \begin{align*}
 	&\sum_{k=2}^{\eps n}\PROB\big(\text{$v$ is reported at least $k$ times in the interval $[t,t+\eps n)$}\big|\hist_{t-1}\big)\\
 	\le&\sum_{k=2}^{\eps n}\binom{\eps n}{k}2^k\gamma^{-kD}\delta^{-k}n^{-k}\le\sum_{k=2}^{\eps n}\big(2\eps\gamma^{-D}\delta^{-1}\big)^k\le\tfrac{4\eps^2\gamma^{-2D}\delta^{-2}}{1-2\eps\gamma^{-D}\delta^{-1}}\le 8\eps^2\gamma^{-2D}\delta^{-2}\,,
 \end{align*}
 where we use the bound $\binom{\epsilon n}{k}\le (\epsilon n)^k$ and sum the resulting geometric series. 
\end{claimproof}
\end{proof}

The proof of Lemma~\ref{lem:probint}\ref{e:probintuv} is similar, and we only focus on the differences.
\begin{proof}[Proof of Lemma~\ref{lem:probint}\ref{e:probintuv}]
  We define  \ModmodRandomEmbedding\, this time reporting a uniform random vertex of $\CANDSET_{G\AlgMap H}^{t-1}(t)\setminus(\im\psi_{t-1}\setminus\{u,v\})$ at each time step $t$, and either embedding $t$ to it (if it is not in $\im\psi_{t-1}$) or otherwise picking as before a uniform random vertex of   $\CANDSET_{G\AlgMap H}^{t-1}(t)\setminus\im\psi_{t-1}$ to embed $t$ to. As before, the embedding distribution generated by this procedure is the same as for \RandomEmbedding.
  We let $E'$ be the number of times $u$ or $v$ are reported in the interval $t\le x<t+\eps n$. Again, the probability that \RandomEmbedding\ uses either $u$ or $v$ is equal to the probability that \ModmodRandomEmbedding\ reports $u$ or $v$ at least once, which by definition is
 $$\EXP\left(E'\:|\:\hist_{t-1}\right)-\sum_{k=2}^{\eps n}\PROB\big(\text{$u$ or $v$ is reported at least $k$ times in the interval $[t,t+\eps n)$}\big|\hist_{t-1}\big)\,.$$
By linearity of expectation, $\EXP\left(E'|\hist_{t-1}\right)$ is equal to the expected number of times $u$ is reported plus the expected number of times $v$ is reported. We now argue that these latter quantities  are \corrig{$(1\pm 4C\alpha)\tfrac{\eps n}{n-t}\pm 8(D+1)\eps^2\gamma^{-2D}\delta^{-2}$.}
%$(1\pm 4C\alpha)\tfrac{\eps n}{n-t}\pm 4(D+1)\eps^2\gamma^{-D}\delta^{-1}$. 
  This follows from calculations in Claim~\ref{cl:probint:E}, with a small change which we now describe. Note that  Claim~\ref{cl:probint:E} deals with \ModRandomEmbedding, where reported vertices are taken from $\CANDSET^{x-1}(x)\setminus(\im\psi_{x-1}\setminus\{v\})$ and not from $\CANDSET^{x-1}(x)\setminus(\im\psi_{x-1}\setminus\{u, v\})$. This is corrected if we rewrite~\eqref{eq:obed} as 
  $$\EXP\left(E'\:|\:\hist_{t-1}\right)=\sum_{x=t}^{t+\eps n-1} \EXP\left(\frac{\mathbbm{1}\{u\in\CANDSET^{x-1}(x)\}+\mathbbm{1}\{v\in\CANDSET^{x-1}(x)\}}{|\CANDSET^{x-1}(x)\setminus \im\psi_{x-1}|\pm 2}\:\Big|\hist_{t-1}\right)\;.$$
  Then the rest of the calculations in Claim~\ref{cl:probint:E} applies (see Footnote~\ref{footnote:biggererror})
  We thus have
 \corrig{$$\EXP\left(E'\:|\:\hist_{t-1}\right)=(1\pm 4C\alpha)\tfrac{2\eps n}{n-t}\pm 16(D+1)\eps^2\gamma^{-2D}\delta^{-2}\,.$$}
 %$$\EXP\left(E'\:|\:\hist_{t-1}\right)=(1\pm 4C\alpha)\tfrac{2\eps n}{n-t}\pm 8(D+1)\eps^2\gamma^{-D}\delta^{-1}\,.$$
 
 Again, it remains to show that the effect of reporting $u$ or $v$ multiple times is small. This time the probability at any step $x$ that one of $u$ and $v$ is reported, conditioning on the history up to time $x-1$, is at most $4\gamma^{-2D}\delta^{-2}n^{-1}$, and by the same calculation as above we conclude that the summation is bounded above by \corrig{$32\eps^2\gamma^{-2D}\delta^{-2}$},
 %$16\eps^2\gamma^{-D}\delta^{-1}$}
 which as before gives~Lemma~\ref{lem:probint}\ref{e:probintuv}.
\end{proof}

We now use Lemma~\ref{lem:probint} to estimate the probability of embedding a vertex to $v$, or to $\{u,v\}$, in the interval $(t_0,t_1]$ (which may be of any length). This time, we do not condition on one typical embedding history up to time $t_0$, but rather on a history ensemble up to time $t_0$ which is not very unlikely. This allows us to drop the typicality restriction, simply because only very few histories can be atypical.

 \begin{lem}\label{lem:probvtx} Given~$D\in \mathbb{N}$ and $\gamma>0$, let $\delta,\alpha_0,\alpha_{2n},C,\eps$ be as in Setting~\ref{set:graphs}. Then the following holds for any $\alpha_0\le\alpha\le\alpha_{2n}$ and all sufficiently large $n$. Suppose
that $G$ is a graph on $[n]$ such that $\LEFTDEG(x)\le D$ for each $x\in V(G)$, and $H$ is an $(\alpha,2D+3)$-quasirandom graph
with $n$ vertices and $p\binom{n}{2}$ edges, with $p\ge\gamma$.
  Let $0\le t_0<t_1\le n-\delta n$. Let $\histens$ be a history ensemble of \RandomEmbedding{} up to time $t_0$, and suppose that $\PROB(\histens)\ge n^{-4}$. Then the following hold for any distinct vertices $u,v\in V(H)$.
  \begin{enumerate}[label=(\alph*)]
   \item\label{e:f1} If $v\not\in\im\psi_{t_0}$ then we have
   \[\PROB^{G\AlgMap H} (v\not\in\im\psi_{t_1}|\histens)=(1\pm 100C\alpha\delta^{-1})\tfrac{n-1-t_1}{n-t_0}\,.\]
   \item\label{e:f2} If $u,v\not\in\im\psi_{t_0}$ then we have
   \[\PROB^{G\AlgMap H}(u,v\not\in\im\psi_{t_1}|\histens)=(1\pm 100C\alpha\delta^{-1})\big(\tfrac{n-1-t_1}{n-t_0}\big)^2\,.\]
  \end{enumerate}
 \end{lem}
 \begin{proof}
 We write $\PROB$ for $\PROB^{G\AlgMap H}$.
We shall first address part~\ref{e:f1}. We divide the interval $(t_0,t_1]$ into $k:=\lceil(t_1-t_0)/\eps n\rceil$ intervals, all but the last of length $\eps n$.
  Let $\histens_0:=\histens$. Let, for each $1\le i<k$, the set $\histens_i$ be the embedding histories up to time $t_0+i\eps n$ of \RandomEmbedding{} which extend histories in $\histens_{i-1}$ and are such that $v\not\in\psi_{t_0+i\eps n}$. Let $\histens_k$ be the embedding histories up to time $t_1$ extending those in $\histens_{k-1}$ such that $v\not\in\psi_{t_1}$. Thus we have
\[\PROB(v\not\in\im\psi_{t_1}|\histens)=\PROB(\histens_k)/\PROB(\histens_0)\,.\]
Finally, for each $1\le i\le k$, let the set $\histens'_{i-1}$ consist of all histories in $\histens_{i-1}$ such that the $(C\alpha,2D+3)$-diet condition holds for $(H,\im\psi_{t_0+(i-1)\eps n})$ and the probability that the $(\eps,C\alpha,t_0+1+(i-1)\eps n)$-cover condition fails, conditioned on $\psi_{t_0+(i-1)\eps n}$, is at most $n^{-3}$. In other words, $\histens'_i$ is the subset of $\histens_i$ consisting of typical histories, satisfying the conditions of Lemma~\ref{lem:probint}. 
  
  We now determine $\PROB(\histens_k)$ in terms of $\PROB(\histens_0)$, and in particular we show inductively that $\PROB(\histens_i)>n^{-5}$ for each $i$. Observe that for any time $t$, the probability (not conditioned on any embedding) that either the $(C\alpha,2D+3)$-diet condition fails for $(H,\im\psi_i)$ for some $i\le t$ or that the $(\eps,C\alpha,t+1)$-cover condition has probability greater than $n^{-3}$ of failing, is at most $2n^{-6}$ by Lemma~\ref{lem:diet}. In other words, for each $i$ we have $\PROB(\histens_i\setminus\histens'_i)\le 2n^{-6}$. Thus by Lemma~\ref{lem:probint}\ref{e:probintv} we have
  \begin{align*}
   \PROB(\histens_i)&=\big(1-(1\pm 10C\alpha)\tfrac{\eps n}{n-t_0-(i-1)\eps n}\big)\PROB(\histens'_{i-1})\pm 2n^{-6}\\
   &=\big(1-(1\pm 10C\alpha)\tfrac{\eps n}{n-t_0-(i-1)\eps n}\big)\big(\PROB(\histens_{i-1})\pm 2n^{-6}\big)\pm 2n^{-6}\\
   &=\big(1-(1\pm 20C\alpha)\tfrac{\eps n}{n-t_0-(i-1)\eps n}\big)\PROB(\histens_{i-1})\,,
  \end{align*}
  where the final equality uses the lower bound $\PROB(\histens_{i-1})\ge n^{-5}$. Similarly, we have $\PROB(\histens_k)=\big(1\pm(1+20C\alpha)\tfrac{\eps n}{n-t_1}\big)\PROB(\histens_{k-1})$.
  
  Putting these observations together, we can compute $\PROB(\histens_k)$:
  \[\PROB(\histens_k)=\big(1\pm(1+20C\alpha)\tfrac{\eps n}{n-t_1}\big)\PROB(\histens_0)\prod_{i=1}^{k-1}\Big(1-(1\pm 20C\alpha)\tfrac{\eps n}{n-t_0-(i-1)\eps n}\Big)\,.\]
 Observe that the approximation $\log(1+x)=x\pm x^2$ is valid for all sufficiently small $x$. In particular, since $n-t_0-(i-1)\eps n\ge n-t_1\ge \delta n$ and by choice of $\eps$, for each $i$ we have
\[\log \Big(1-(1\pm 20C\alpha)\tfrac{\eps n}{n-t_0-(i-1)\eps n}\Big)=-(1\pm 30C\alpha)\tfrac{\eps n}{n-t_0-(i-1)\eps n}\,.\]
Thus we obtain
 \begin{align}
\nonumber  \log\PROB(\histens_k)&=\log\PROB(\histens_0)\pm (1+30C\alpha)\tfrac{\eps n}{n-t_1}-\sum_{i=1}^{k-1}(1\pm 30C\alpha)\tfrac{\eps n}{n-t_0-(i-1)\eps n}\\
\nonumber  &=\log\PROB(\histens_0)\pm 2\delta^{-1}\eps-(1\pm 40C\alpha)\int_{x=0}^{(k-1)\eps n}\tfrac{1}{n-t_0-x}\, \mathrm{d}x\\
\nonumber  &=\log\PROB(\histens_0)\pm 2\delta^{-1}\eps-(1\pm 50C\alpha)\big(\log(n-t_0)-\log(n-1-t_1)\big)\\
 \label{eq:histensnotsmall} &=\log\PROB(\histens_0)+\log\tfrac{n-1-t_1}{n-t_0}\pm 2\delta^{-1}\eps\pm 50C\alpha\log\delta^{-1}\,,
 \end{align}
 where we use $t_1\le n-\delta n$, and we justify that the integral and sum are close by observing that for each $i$ in the summation, if $(i-1)\eps n\le x\le i\eps n$ then we have
 \[\tfrac{1}{n-t_0-i\eps n}\le\tfrac{1}{n-t_0-x}\le\tfrac{1}{n-t_0-(i-1)\eps n}\le(1+\alpha)\tfrac{1}{n-t_0-i\eps n}\,,\]
 where the final inequality uses $n-t_0-i\eps n\le n-t_1\le \delta n$ and the choice of $\eps$. By choice of $\eps$, this gives part~\ref{e:f1}. Furthermore, \eqref{eq:histensnotsmall}, and the fact $t_1\le n-\delta n$, imply that $\PROB(\histens_k)\ge n^{-5}$. Since the $\histens_i$ form a decreasing sequence of events the same bound holds for each $\histens_i$.
   
\medskip
   For part~\ref{e:f2}, we use the identical approach, replacing Lemma~\ref{lem:probint}\ref{e:probintv} with~Lemma~\ref{lem:probint}\ref{e:probintuv}. Since the difference between these equations is a factor of $2$, we obtain twice all the terms other than the term $\log\PROB(\histens_0)$ in the above equation, and hence the second statement of the claim.
 \end{proof}

Next, we estimate the probability that the edge $uv\in E(H)$ is used by \RandomEmbedding\ when embedding $G$ to $H$. The idea is the following. In order for $uv$ to be used, there must be some $xy\in G$ such that $x$ is embedded to $u$ and $y$ to $v$, or vice versa. These events are disjoint, and so it suffices to estimate the probability of each separately and sum them. Without loss of generality, we can assume $x$ is embedded before $y$. We need to calculate the probability that $x$ is embedded to $u$ and $y$ to $v$. In other words, we need that all left-neighbours of $x$ are embedded to neighbours of $u$, all left-neighbours of $y$ are embedded to vertices of $v$, other vertices are not embedded to $\{u,v\}$, and when we come to embed $x$ and $y$ we actually do embed them to $u$ and $v$. The point of phrasing it like this is that, provided the diet condition holds, we can estimate accurately all the (conditional) probabilities of embedding individual vertices in $N(x)\cup N(y)\cup\{x,y\}$ to neighbourhoods or to $u$ or $v$, while Lemma~\ref{lem:probvtx} gives accurate estimates for the probability of any other vertex being embedded to $u$ or $v$. Putting this together yields the desired accurate estimate for the probability that we have $x\AlgMap u$ and $y\AlgMap v$.

\begin{lem}\label{lem:probedge}
 Given~$D\in \mathbb{N}$, and $\gamma>0$, let constants $\delta,\eps,C,\alpha_0,\alpha_{2n}$ be as in Setting~\ref{set:graphs}. Then the following holds for any $\alpha_0\le\alpha\le\alpha_{2n}$ and all sufficiently large $n$. Suppose
that $G$ is a graph on $[n]$ such that $\LEFTDEG(x)\le D$ for each $x\in V(G)$, and $H$ is an $(\alpha,2D+3)$-quasirandom graph
with $n$ vertices and $p\binom{n}{2}$ edges, with $p\ge\gamma$. Let $uv$ be an edge of $H$. When \RandomEmbedding{} is run to embed $G[{\scriptstyle [n-\delta n]}]$ into $H$, the probability that an edge of $G$ is embedded to $uv$ is
\[\big(1\pm500C\alpha\delta^{-1}\big)^{4D+2}p^{-1}n^{-2}\cdot 2e(G)\,.\]
\addtocounter{counterProbedge}{\value{thm}}
\end{lem}
\begin{proof}
 We first calculate the probability that a given pair $(x,y)$, such that $xy$ is an edge of $G$, is embedded to $(u,v)$, in that order.  Without loss of generality, suppose that $x<y$. Let $z_1,\dots,z_k$ be the vertices $\LNBH(x)\cup\LNBH(y)\setminus\{x,y\}$ in increasing order. Let $j\in\{0,\ldots,k\}$ be such that $z_j<x<z_{j+1}$ (where the case $j=0$ and $j=k$ corresponds to the situations when all $z_i$'s are to the right or to the left of $x$, respectively; in these cases some notation below has to be modified in a straightforward way).
Define time intervals using $z_1,\ldots,z_j,x,z_{j+1},\ldots,z_k,y$ as separators: $I_0=[1,z_1-1]$, $I_1=[z_1+1,z_2-1]$, \ldots, $I_j=[z_j+1,x-1]$, $I_{j+1}=[x+1,z_{j+1}-1]$, \ldots, $I_{k+1}=[z_k+1,y-1]$.

We now define a nested collection of events, the first being the trivial (always satisfied) event and the last being the event $\{x\AlgMap u,y\AlgMap v\}$, whose probability we wish to estimate. These events are simply that we have not yet (by given increasing times in \RandomEmbedding) made it impossible to have $\{x\AlgMap u,y\AlgMap v\}$. We will see that we can estimate accurately the probability of each successive event, conditioned on its predecessor.

Let $\histens'_{-1}$ be the trivial (always satisfied) event. If $\histens'_{i-1}$ is defined, we let $\histens_i$ be the event that $\histens'_{i-1}$ holds intersected with the event that 
\begin{enumerate}[label=(A\arabic*)]
 \item\label{e:A1} (if $i\le j$:) no vertex of $G$ in the interval $I_i$ is mapped to $u$ or $v$, or
\item\label{e:A2} (if $i> j$:) no vertex of $G$ in the interval $I_i$ is mapped to $v$.
\end{enumerate}
In other words, $\histens_i$ is the event that we have not covered $u$ or $v$ in the interval $I_i$. It turns out that we do not need to know anything else about the embeddings in the interval $I_i$.

If $\histens_i$ is defined, we let $\histens'_i$ be that event that $\histens_i$ holds and that
 \begin{enumerate}[label=(B\arabic*)]
 \item\label{e:Ch1} (if $i<j$:) 
 \begin{enumerate}[label=(\roman*)]
\item\label{e:Ch1A} (subcase $z_{i+1}\in \LNBH(x)\setminus\LNBH(y)$:)
 we have the event $z_{i+1}\AlgMap \NBH_H(u)\setminus\{v\}$,
\item\label{e:Ch1B} (subcase $z_{i+1}\in \LNBH(y)\setminus\LNBH(x)$:)
  we have the event $z_{i+1}\AlgMap \NBH_H(v)\setminus\{u\}$,
\item\label{e:Ch1C} (subcase $z_{i+1}\in \LNBH(x)\cap \LNBH(y)$:)
  we have the event $z_{i+1}\AlgMap \NBH_H(u)\cap\NBH_H(v)$,
 \end{enumerate} 
 \item\label{e:Ch2} (if $i=j$:)
 we have the event $x\AlgMap u$,
 \item\label{e:Ch3} (if $j<i\le k:$)   
 we have the event $z_{i}\AlgMap \NBH_H(v)\setminus\{u\}$ (unlike the range $i<j$, there are no subcases here, as necessarily $z_i\in \LNBH(y)\setminus\LNBH(x)$),
 \item\label{e:Ch4} (if $i=k+1:$) we have the event $y\AlgMap v$.
 \end{enumerate}
Again, in order for $\{x\AlgMap u,y\AlgMap v\}$ to occur we obviously need that a neighbour of $x$ is embedded to a neighbour of $u$ and so on, hence the above conditions.

By definition, we have $\histens'_{k+1}=\{x\AlgMap u,y\AlgMap v\}$. Since we have $\histens'_i\subseteq\histens_i\subseteq\histens'_{i-1}$ for each $i$ and $\histens'_{-1}$ is the sure event, we see
\begin{align}
\PROB\left(x\AlgMap u,y\AlgMap v\right)=\prod_{i=0}^{k+1}
\frac{\PROB(\histens_i)}{\PROB(\histens'_{i-1})}\cdot\frac{\PROB(\histens'_i)}{\PROB(\histens_i)}=\prod_{i=0}^{k+1}
\PROB\left(\histens_i\:|\:\histens'_{i-1}\right)
\PROB\left(\histens'_i\:|\:\histens_i\right)
\label{eq:nicelygrouped}
\;.
\end{align}
Thus, we need to estimate the factors in~\eqref{eq:nicelygrouped}. This is done in the two claims below. In each claim we assume $\PROB(\histens'_i),\PROB(\histens_i)>n^{-4}$. This assumption is justified, using an implicit induction, since the smallest of all the events we consider is $\histens'_{k+1}$, whose probability according to the following~\eqref{eq:probxy} is bigger than $n^{-4}$.

\begin{claimProbedge}\label{cl:ng1}
We have
\begin{equation*}
\prod_{i=0}^{k+1} \PROB\left(\histens_i\:|\:\histens'_{i-1}\right)=(1\pm 200C\alpha\delta^{-1})^{2k+2}\cdot \frac{(n-x)(n-y)}{n^2}\;.
\end{equation*}
\end{claimProbedge}
\begin{claimproof}
By definition of~\ref{e:A1}, for each $i=0,\ldots, j$, we have 
\begin{equation}\label{eq:cancelme}
\PROB\left(\histens_i\:|\:\histens'_{i-1}\right)=(1\pm 200C\alpha\delta^{-1})\cdot \frac{(n-1-\max(I_{i}))^2}{(n-\min(I_i)+1)^2}
\end{equation}
by Lemma~\ref{lem:probvtx}\ref{e:f2}, with $\histens=\histens'_{i-1}$.   Note that looking at two consecutive indices $i$ and $i+1$ in~\eqref{eq:cancelme} we have cancellation of the former nominator and the latter denominator, $n-1-\max(I_{i})=n-\min(I_{i+1})+1$. Thus,
\begin{equation}\label{eq:cancelmex}
\prod_{i=0}^{j} \PROB\left(\histens_i\:|\:\histens'_{i-1}\right)=(1\pm 200C\alpha\delta^{-1})^{2j+2}\cdot \frac{(n-x)^2}{n^2}\;.
\end{equation}
To express $\prod_{i=j+1}^{k+1} \PROB\left(\histens_i\:|\:\histens'_{i-1}\right)$, by definition of~\ref{e:A2} we have to repeat the above replacing Lemma~\ref{lem:probvtx}\ref{e:f2} by Lemma~\ref{lem:probvtx}\ref{e:f1}. We get that
\begin{equation}\label{eq:cancelmexx}
\prod_{i=j+1}^{k+1} \PROB\left(\histens_i\:|\:\histens'_{i-1}\right)=(1\pm 200C\alpha\delta^{-1})^{2(k-j)+2}\cdot \frac{n-y}{n-x}\;.
\end{equation}
Putting~\eqref{eq:cancelmex} and~\eqref{eq:cancelmexx} together, we get the statement of the claim.
\end{claimproof}

\begin{claimProbedge}\label{cl:ng2}
We have
$$
\prod_{i=0}^{k+1}\PROB\left(\histens'_i\:|\:\histens_i\right)=(1\pm 100C\alpha)^{2D}\cdot \frac{1}{p(n+1-x)(n+1-y)}\;.
$$
\end{claimProbedge}
\begin{claimproof}
Suppose that we have embedded up to vertex $\max(I_i)$, and that $\histens_i$ holds. The probability of the event $\histens'_i$ depends on which of the cases in~\ref{e:Ch1}--\ref{e:Ch3} applies. When $\histens'_i$ is defined using~\ref{e:Ch1}\ref{e:Ch1A} then the probability $\PROB(\histens'_i|\histens_i)$ is equal to $\PROB(\{z_{i+1}\AlgMap \NBH_H(u)\setminus\{v\}\}|\histens_i)$. Let $X:=N_H\big(\psi(\LNBH_G(z_{i+1}))\big)\setminus\im\psi_{z_{i+1}-1}$ be the set of vertices in $H$ to which we could embed $z_{i+1}$, given the embedding of all vertices before $z_{i+1}$. Suppose that the $(C\alpha,2D+3)$-diet condition holds for $(H,\im\psi_{z_{i+1}-1})$. Then we have 
\begin{align*}
\PROB\left(z_{i+1}\AlgMap \NBH_H(u)\setminus\{v\}|\histens_i\right)&=\frac{|(\NBH_H(u)\setminus\{v\})\cap X|}{|X|}=\frac{|\NBH_H(u)\cap X|\pm 1}{|X|}
\\
&=
\frac{(1\pm C\alpha)p^{1+\LEFTDEG(z_{i+1})}(n-(z_{i+1}-1))\pm 1}
{(1\pm C\alpha)p^{\LEFTDEG(z_{i+1})}(n-(z_{i+1}-1))}
=(1\pm 4C\alpha)p\;,
\end{align*}
where the last line uses the $(C\alpha,2D+3)$-diet condition for $(H,\im\psi_{z_{i+1}-1})$ twice, in the denominator with the set $\psi(\LNBH(z_{i+1}))$ and in the numerator with the set $\{u\}\cup\psi(\LNBH(z_{i+1}))$. Recall that we assume the event $\histens_i$, and so we have $u\not\in \im \psi_{z_{i+1}-1}$. Therefore, the set $\{u\}\cup \psi(\LNBH_G(z_{i+1}))$ has indeed size $1+\LEFTDEG(z_{i+1})$.

Likewise, when $\histens'_i$ is defined using~\ref{e:Ch1}\ref{e:Ch1B}, using~\ref{e:Ch1}\ref{e:Ch1C}, or using~\ref{e:Ch3} then $\PROB(\histens'_i|\histens_i)$ is the probability of $\{z_{i+1}\AlgMap \NBH_H(v)\setminus\{u\}\}$, of $\{z_{i+1}\AlgMap \NBH_H(u,v)\}$, or of $\{z_i\AlgMap \NBH_H(v)\setminus\{u\}$, respectively. If the $(C\alpha,2D+3)$-diet condition holds for $(H,\im\psi_{z_{i+1}-1})$, this probability is equal to $(1\pm 4C\alpha)p$, $(1\pm 4C\alpha)p^2$, or $(1\pm 4C\alpha)p$, respectively. 

Let us now deal with the terms $\PROB\left(\histens'_j\:|\:\histens_j\right)$ and $\PROB\left(\histens'_{k+1}\:|\:\histens_{k+1}\right)$ which correspond to~\ref{e:Ch2} and~\ref{e:Ch4}, respectively. Suppose first that $\histens_j$ holds. In particular, $\LNBH(x)$ is embedded to $\NBH_H(u)$.
Suppose first that the $(C\alpha,2D+3)$-diet condition for $(H,\im\psi_{x-1})$ holds. With this, conditioning on the embedding up to time $x-1$, the probability of embedding $x$ to $u$ is $(1\pm 2C\alpha)p^{-\LEFTDEG(x)}\tfrac{1}{n+1-x}$. Similarly, if the $(C\alpha,2D+3)$-diet condition for $(H,\im\psi_{y-1})$ holds, the probability of embedding $y$ to $v$, provided $\LNBH(y)$ is embedded to $\NBH_H(v)$, and conditioning on the embedding up to time $y-1$, is $(1\pm 2C\alpha)p^{-\LEFTDEG(y)}\tfrac{1}{n+1-y}$.

Thus, letting $\cF$ be the event that the $(C\alpha,2D+3)$-diet condition fails at least once for $(H,\im\psi_t)$, where $t$ runs between $1$ and $y$, we have
\begin{multline*}
\prod_{i=0}^{k+1}\PROB\left(\histens'_i\:|\:\histens_i\right)
=\Big(\big((1\pm 4C\alpha)p\big)^{\ell_1}\cdot\big((1\pm 4C\alpha)p^2\big)^{\ell_2} \\
\cdot
(1\pm 2C\alpha)p^{-\LEFTDEG(x)}\tfrac{1}{n+1-x}
\cdot
(1\pm 2C\alpha)p^{-\LEFTDEG(y)}\tfrac{1}{n+1-y}
\Big) \pm \PROB(\cF)\;,
\end{multline*}
where we write $\ell_1$ for the number of times~\ref{e:Ch1}\ref{e:Ch1A}, \ref{e:Ch1}\ref{e:Ch1B}, or \ref{e:Ch3} applies, and $\ell_2$ for the number of times~\ref{e:Ch1}\ref{e:Ch1C} applies. We have $\ell_1+2\ell_2=\LEFTDEG(x)+\LEFTDEG(y)-1$. Indeed, $\ell_1$ and $\ell_2$ count the left neighbours of $x$ and $y$, but $x$, which is a left neighbour of $y$, is omitted. Finally, $\PROB(\cF)\le 2n^{-9}$ by Lemma~\ref{lem:diet}. Thus we obtain
\begin{equation*}
\prod_{i=0}^{k+1}\PROB\left(\histens'_i\:|\:\histens_i\right)
=(1\pm 4C\alpha)^{\ell_1+\ell_2+2}p^{-1} 
\cdot
\tfrac{1}{n+1-x}
\cdot
\tfrac{1}{n+1-y}
 \pm 2n^{-9}\;,
\end{equation*}
which gives the claim since $\ell_1+\ell_2+2\le 2D+1$.
\end{claimproof}

Plugging Claims~\ref{cl:ng1} and~\ref{cl:ng2} into~\eqref{eq:nicelygrouped}, we get 
\begin{equation}\label{eq:probxy}
\PROB\left(x\AlgMap u,y\AlgMap v\right)=(1\pm 500C\alpha\delta^{-1})^{4D+2}\cdot p^{-1}n^{-2}\;.
\end{equation}

 We now sum over the choices of $(x,y)$ such that $xy\in E(G)$. There are $2e(G)$ such choices, so we conclude that the probability that some edge of $G$ is embedded by \RandomEmbedding\ to $uv$ is
 \[\big(1\pm500C\alpha\delta^{-1}\big)^{4D+2}p^{-1}n^{-2}\cdot 2e(G)\]
 as desired.
\end{proof}

We can now estimate the probability that, again for fixed $G$ and $H$, at least one edge in a given star in $H$ is used by \RandomEmbedding.
\begin{lem}\label{lem:probstar}
 Given $D\in\mathbb{N}$ and $\gamma>0$, let the constants $\delta,\eps,\alpha_0,\alpha_{2n},C$ be as in Setting~\ref{set:graphs}. Then the following holds for any $\alpha_0\le\alpha\le\alpha_{2n}$ and all sufficiently large $n$. Suppose that $G$ is a graph on $[n]$ such that $\LEFTDEG(x)\le D$ for each $x\in V(G)$, with at least $n/4$ edges and maximum degree $\Delta(G)\le n/\log n$, and $H$ is an $(\alpha,2D+3)$-quasirandom graph with $n$ vertices and $p\binom{n}{2}$ edges, where $p\ge\gamma$. Let $u_1,\dots,u_k,v$ be vertices of $H$ for some $k\le 2D+3$, and suppose $u_iv$ is an edge of $H$ for each $i$. When \RandomEmbedding\ is run to embed $G[{\scriptstyle [n-\delta n]}]$ into $H$, the probability that there is at least one $u_iv$ to which some edge of $G$ is embedded is
 \[\big(1\pm1000C\alpha\delta^{-1}\big)^{4D+2}p^{-1}n^{-2}\cdot 2ke(G)\,.\]
\end{lem}
\begin{proof} Given $u_1,\dots,u_k,v$ and $G$ and $H$, let $S$ be the event that there is at least one $u_iv$ to which some edge of $G$ is embedded.

 The expected number of edges $u_iv$ embedded to by \RandomEmbedding\ is, by Lemma~\ref{lem:probedge} and linearity of expectation,
 \[E:=\big(1\pm500C\alpha\delta^{-1}\big)^{4D+2}p^{-1}n^{-2}\cdot 2kE(G)\,,\]
 and by inclusion-exclusion, we have
 \[E-\sum_{1\le i<i'\le k}\PROB\big(u_iv\text{ and }u_{i'}v\text{ are embedded to by \RandomEmbedding}\big)\le\PROB(S)\le E\,.\]
 We thus simply have to show that the above sum, which has $\binom{k}{2}\le\binom{2D+3}{2}$ terms, is small. We will show that the probability of \RandomEmbedding\ embedding to any two fixed edges $uv,u'v$ is small. This probability is equal to the sum over triples $x,x',y\in V(G)$ such that $xy,x'y\in E(G)$ of the probability that $x\AlgMap u$, $x'\AlgMap u'$ and $y\AlgMap v$. For any given $y\in V(G)$ there are at most $\deg_G(y)^2$ choices of $(x,x')$, so by Lemma~\ref{lem:squarebound}, there are at most $2Dn\Delta(G)$ such triples. It is now enough to make the estimate for one such triple. Assuming the $(C\alpha,2D+3)$-diet condition holds throughout \RandomEmbedding, we embed each of $x$, $x'$ and $y$ uniformly at random into a set of size at least $\tfrac12p^D\delta n\ge\tfrac12\gamma^D\delta n$, so the probability of the event $x\AlgMap u,x'\AlgMap u',y\AlgMap v$ is at most $8\gamma^{-3D}\delta^{-3}n^{-3}$. Finally, the probability of the $(C\alpha,2D+3)$-diet condition failing for some $(H,\im\psi_i)$ is by Lemma~\ref{lem:dietcon} at most $2n^{-9}$. Putting this together, we have
 \[\PROB(S)=\big(1\pm500C\alpha\delta^{-1}\big)^{4D+2}p^{-1}n^{-2}\cdot 2ke(G)\pm\tbinom{2D+3}{2}\cdot 2Dn\Delta(G)\cdot 8\gamma^{-3D}\delta^{-3}n^{-3}\pm 2n^{-9}\,.\]
 Because $e(G)\ge n/4$ the first term in the above is $\Theta(n^{-1})$, while since $\Delta(G)\le n/\log n$ the other two terms are of asymptotically smaller order. Since $n$ is sufficiently large, this gives the desired result.
\end{proof}

In Lemma~\ref{lem:probstar} we estimated the probability of using an edge in a star with a given centre and a given set $R$ of ends. In particular, looking at all stars in $H$ whose ends are $R$, we get an estimate of the expected number of them from which an edge is used in the embedding. In the following lemma we prove an upper bound on the second moment of this random variable.

\begin{lem}\label{lem:starsecond}
 Let $D\in\mathbb{N}$ and let $\gamma>0$. Let $\delta,\eps,c,C,\alpha_0,\alpha_{2n}$ be as in Setting~\ref{set:graphs}. Then the following holds for any $\alpha_0\le\alpha\le\alpha_{2n}$ and all sufficiently large $n$. Suppose that $G$ is a graph on $[n]$ such that $\LEFTDEG(x)\le D$ for each $x\in V(G)$, with at least $n/4$ edges and maximum degree $\Delta(G)\le c n/\log n$, and $H$ is an $(\alpha,2D+3)$-quasirandom graph with $n$ vertices and $p\binom{n}{2}$ edges, where $p\ge\gamma$. Given $R\subset V(H)$ with $|R|\le 2D+3$ and any subset $T$ of $\NBH_H(R)$, let $X$ count the number of vertices $v\in T$ such that an edge from $v$ to $R$ is used by \RandomEmbedding\ when embedding $G$ to $H$. Then we have
  \[\EXP(X^2)\le 2^{30}D^4\Delta(G)\gamma^{-4D}\delta^{-4}\,.\]
\end{lem}
 \begin{proof}
 We can write $X=\sum_{v\in T}W_v$, where $W_v$ is the indicator random variable of the event that some edge from $R$ to $v$ is used in embedding $G$. We have
 \[\EXP(X^2)=\sum_{(v,v')\in T^2}\EXP(W_vW_{v'})=\EXP(X)+2\sum_{\{v,v'\}\subset T}\EXP(W_vW_{v'})\,.\]
 Since $e(G)\le Dn$, by Lemma~\ref{lem:probstar}, applied with $\{u_1,\dots,u_k\}=R$ and for each $v\in T$, we have
 \[\EXP(X)\le \big(1+1000C\alpha\delta^{-1}\big)^{4D+2}p^{-1}n^{-2}\cdot 2|R|\cdot Dn\cdot |T|\le  4\gamma^{-1}D(2D+3)\,,\]
 where we use $|R|\le 2D+3$ and $|T|\le n$.
 Thus the main task is thus to estimate $\EXP(W_vW_{v'})$ for $v\neq v'$. Now $W_vW_{v'}$ is equal to~1 if and only if there is an edge of $G$ embedded to some edge between $R$ and $v$, and another to an edge between $R$ and $v'$. So, in order to refine our strategy, for $v\in T$ and $u\in R$, let $Y_{v,u}$ be the indicator random variable of the event that the edge $uv$ is used in embedding $G$. For each $\{v,v'\}\subset T$ we have
\begin{equation}
\label{summeifyoulike}
\EXP(W_vW_{v'})= \sum_{u,u'\in R,u\neq u'}\EXP(Y_{v,u}Y_{v',u'})+\sum_{u\in R}\EXP(Y_{v,u}Y_{v',u})\;.
\end{equation}
 First, we focus on the first term of the right-hand side of~\eqref{summeifyoulike}. That is, we need to find an upper bound for the probability that two given disjoint edges $xy$ and $x'y'$ of $G$ are embedded to respectively $uv$ and $u'v'$ for some fixed $u,u'\in R$ and fixed $v,v'$. As \RandomEmbedding\ runs, either for some $t$ we observe that the $(C\alpha,2D+3)$-diet condition fails for $(H,\im\psi_t)$, or it is successful and at each time $t$, the vertex $t$ is embedded uniformly at random into a set of size at least $\tfrac12\gamma^D\delta n$. The probability of the former occurring is at most $2n^{-9}$ by Lemma~\ref{lem:diet}, while in the latter case the probability of embedding $x,y,x',y'$ to $u,v,u',v'$ in that order is at most $16\gamma^{-4D}\delta^{-4}n^{-4}$. Putting these together the probability of $xy,x'y'$ being embedded to $uv,u'v'$ in that order is at most $32\gamma^{-4D}\delta^{-4}n^{-4}$. Summing over the at most $8\binom{e(G)}{2}\le 8\binom{Dn}{2}$ choices of edges $xy,x'y'$ and their orderings, we get
 $$\EXP(Y_{v,u}Y_{v',u'})\le 8\binom{Dn}{2}\cdot 32\gamma^{-4D}\delta^{-4}n^{-4}\;.$$
There are exactly $|R|^2-|R|\le(2D+3)^2$ choices of distinct vertices $u,u'\in R$. Hence
\begin{equation}\label{zne:1}
\sum_{u,u'\in R,u\neq u'}\EXP(Y_{v,u}Y_{v',u'})\le (2D+3)^2\cdot 8\binom{Dn}{2}\cdot 32\gamma^{-4D}\delta^{-4}n^{-4}\le 2^{15}D^4\gamma^{-4D}\delta^{-4}n^{-2}\,.
\end{equation}

 Next, we focus on the second term of the right-hand side of~\eqref{summeifyoulike}. That is, we now find an upper bound for the probability that \RandomEmbedding\ uses both $uv$ and $uv'$ for some $u\in R$. The only way this can happen is that for some $x,y,y'\in V(G)$ with $xy,xy'\in E(G)$, the vertex $x$ is embedded to $u$ and $y,y'$ to $v,v'$. Again, by Lemma~\ref{lem:diet}, the probability that a fixed such triple $x,y,y'$ are embedded to $u,v,v'$ is at most $2n^{-9}+8\gamma^{-3D}\delta^{-3}n^{-3}$. By Lemma~\ref{lem:squarebound} there are at most $2Dn\Delta(G)$ such triples. Hence, we get
 $$\EXP(Y_{v,u}Y_{v',u})\le 2Dn\Delta(G)\cdot(2n^{-9}+8\gamma^{-3D}\delta^{-3}n^{-3})\le 2Dn\Delta(G)\cdot 16\gamma^{-3D}\delta^{-3}n^{-3}\;.$$
There are exactly $|R|\le 2D+3$ choices of $u$, so the probability that \RandomEmbedding\ uses both $uv$ and $uv'$ for some $u\in R$ is at most
\begin{equation}\label{eq:zne2}
 \sum_{u\in R}\EXP(Y_{v,u}Y_{v',u})\le (2D+3)\cdot 2Dn\Delta(G)\cdot 16\gamma^{-3D}\delta^{-3}n^{-3}\le 2^{10}D^2\Delta(G)\gamma^{-3D}\delta^{-3}n^{-2}\,
 .
\end{equation}
 
 We can now plug in~\eqref{zne:1} and~\eqref{eq:zne2} into~\eqref{summeifyoulike},
 \[\EXP(W_vW_{v'})\le 2^{20}D^4\Delta(G)\gamma^{-4D}\delta^{-4}n^{-2}\,.\]
 Summing over the at most $n^2$ choices of $v,v'\in T$, we obtain the desired bound.
 \end{proof}

We are now in a position to prove Lemma~\ref{lem:aggregate}. 
\begin{proof}[Proof of Lemma \ref{lem:aggregate}]
We define $\hat{p}$ by $e(H^*_0)=\hat{p}\binom{n}{2}$. By assumption we have $\hat{p}=(1\pm\eta)\gamma$. 

 Our aim is to show that with high probability, for any given $s$, either \PackingProcess\ fails before completing stage $s$ or the pair $(H_s,H^*_0)$ is $(\alpha_s,2D+3)$-coquasirandom. Let $S$ be a set of at most $2D+3$ vertices in $V(H^*_0)$, and let $R\subseteq S$. Recall that for $(H_s,H^*_0)$ to be $(\alpha_s,2D+3)$-coquasirandom means that $\NBH_{H_{s}}(R)\cap \NBH_{H^*_0}(S\setminus R)$ has about the size one would expect if both graphs were random. For each $1\le i\le s$, let
 \[Y_i=\big|\NBH_{H_{i-1}}(R)\cap \NBH_{H^*_0}(S\setminus R)\setminus \NBH_{H_i}(R)\big|\,.\]
 In other words, $Y_i$ is the number of vertices which are removed to form $\NBH_{H_{i}}(R)\cap \NBH_{H^*_0}(S\setminus R)$ when we embed $G_i[{\scriptstyle [n-\delta n]}]$ to $H_{i-1}$. To prove coquasirandomness of $(H_s,H^*_0)$, what we want is for $\sum_{i=1}^sY_i$ to be sufficiently concentrated to take a union bound over choices of $R$ and $S$. For this purpose we aim to apply Lemma~\ref{lem:freedman} with $\cE$ being the event that after each stage $i=0,\ldots,s-1$ the pair $(H_i,H^*_0)$ is $(\alpha_i,2D+3)$-coquasirandom. The probability space in which we work is the set of all possible histories of \RandomEmbedding, and the sequence of partitions required by Lemma~\ref{lem:freedman} is given by the histories up to increasing times $1\le i\le s$ of \RandomEmbedding. We thus have to estimate $\EXP(Y_s|H_{s-1})$ and $\Var(Y_s|H_{s-1})$ only in the case $(H_{s-1},H^*_0)$ is $(\alpha_{s-1},2D+3)$-coquasirandom.
 
 So suppose that $(H_{s-1},H^*_0)$ is $(\alpha_{s-1},2D+3)$-coquasirandom. Let $p_s$ be such that $p_s\binom{n}{2}=e(H_s)=e(H_0)-\sum_{i=1}^s e(G_i[{\scriptstyle [n-\delta n]}])$. Then by Lemma~\ref{lem:probstar} and linearity of expectation, we have
 \begin{align}
 \nonumber
  \EXP(Y_s|H_{s-1})&=(1\pm\alpha_{s-1})p_{s-1}^{|R|}\hat{p}^{|S\setminus R|}n\cdot \big(1\pm1000C\alpha_{s-1}\delta^{-1}\big)^{4D+2}p^{-1}_{s-1}n^{-2}\cdot 2|R|e(G_s[{\scriptstyle [n-\delta n]}])\\
 \label{Iaahm}
  &=\big(2|R|\pm 10^6CD^2\delta^{-1}\alpha_{s-1}\big)p_{s-1}^{|R|-1}\hat{p}^{|S\setminus R|}e(G_s[{\scriptstyle [n-\delta n]}])/n\,.
 \end{align}
 We now need to estimate the sum $\sum_{i=1}^s\EXP(Y_i|H_{i-1})$, on the assumption that each $(H_{i-1},H^*_0)$ is $(\alpha_{i-1},2D+3)$-coquasirandom. We first estimate the sum of the main terms of~\eqref{Iaahm}:
 \begin{align}
 \nonumber
  &\sum_{i=1}^s2|R|p_{i-1}^{|R|-1}\hat{p}^{|S\setminus R|}e(G_i[{\scriptstyle [n-\delta n]}])/n\\
  \label{eq:continueme}
\JUSTIFY{we have $e(G_i[{\scriptstyle [n-\delta n]}])=(p_{i-1}-p_i)\binom{n}{2}$}  =&\sum_{i=1}^s|R|p_{i-1}^{|R|-1}(p_{i-1}-p_i)\hat{p}^{|S\setminus R|}(n-1)\;.
\end{align}
Note that for every $x,h\in[0,1]$ and $a\in\mathbb{N}$, we have $(x+h)^a-x^a=ah(x+h)^{a-1}\pm 2^ah^2$. We use this with $x:=p_i$, $h:=p_{i-1}-p_i$, and $a:=|R|$, and observe that $\binom{n}{2}(p_{i-1}-p_i)=e(G_i[{\scriptstyle [n-\delta n]}])\le Dn$. We continue~\eqref{eq:continueme} as follows:
\begin{align}
\nonumber
\sum_{i=1}^s2|R|p_{i-1}^{|R|-1}\hat{p}^{|S\setminus R|}e(G_i[{\scriptstyle [n-\delta n]}])/n &=(n-1)\hat{p}^{|S\setminus R|}\sum_{i=1}^s\Big(\big(p_{i-1}^{|R|}-p_i^{|R|}\big)\pm 16D^22^{|R|}/n^2\Big)\\
\nonumber
  &=(n-1)\hat{p}^{|S\setminus R|}\big(p_0^{|R|}-p_s^{|R|}\big)\pm 64D^22^{|R|} \\
\label{eq:NCBHQ}
  &=\big(p_0^{|R|}-p_s^{|R|}\big)\hat{p}^{|S\setminus R|}n\pm 100D^22^{2D+3}\,.
 \end{align}
 Next, we bound the sum of the error terms of~\eqref{Iaahm}:
 \begin{align}
 \nonumber
  &\sum_{i=1}^s10^6CD^2\delta^{-1}\alpha_{i-1}p_{s-1}^{|R|-1}\hat{p}^{|S\setminus R|}e(G_i[{\scriptstyle [n-\delta n]}])/n\\
  \nonumber
\JUSTIFY{we have $e(G_s)\le Dn$} 
  &\le \int_{-\infty}^s10^7CD^3\delta^{-1}\hat{p}^{|S\setminus R|}\alpha_x\,\mathrm{d}x\\
\label{eq:NCBLQ}
\JUSTIFY{by~\eqref{eq:defconsts}}  &\le \hat{p}^{|S\setminus R|}\gamma^{4D+6}\alpha_s n/4\,.
 \end{align}
Plugging~\eqref{eq:NCBHQ} and~\eqref{eq:NCBLQ} into~\eqref{Iaahm}, we get
 \[\sum_{i=1}^s\EXP(Y_i|H_{i-1})=\big(p_0^{|R|}-p_s^{|R|}\big)\hat{p}^{|S\setminus R|}n\pm \hat{p}^{|S\setminus R|}\gamma^{4D+6}\alpha_s n/2\,,\]
 provided that $H_{i-1}$ is $(\alpha_{i-1},2D+3)$-quasirandom for each $1\le i\le s$.
 
 Let us write $\Delta:=cn/\log n$.
 
 We wish to estimate $\Var(Y_i|H_{i-1})$. Trivially, we have $\Var(Y_s|H_{s-1})\le\EXP(Y_i^2|H_{i-1})$. By Lemma~\ref{lem:starsecond}, \begin{equation*}
\EXP(Y_s^2|H_{s-1})\le 2^{30}D^4\Delta(G_s)\gamma^{-4D}\delta^{-4}\le 2^{30}D^4\Delta\gamma^{-4D}\delta^{-4} \;.
\end{equation*}
Summing this up, we obtain
\begin{equation*}
\sum_{i=1}^{s}\EXP(Y_s^2|H_{s-1})\le  2^{31}D^4\Delta\gamma^{-4D}\delta^{-4}n=:\sigma^2 \;.
\end{equation*}

 Furthermore, the range of each $Y_i$ is at most $|S|\Delta(G_i)\le |S|\Delta$. We apply Lemma~\ref{lem:freedman} with $\sigma^2$ as above, $\rho=\eps n$ and $\cE$ the event that the pair $(H_i,H^*_0)$ is $(\alpha_i,2D+3)$-coquasirandom for each $0\le i\le s-1$. We obtain that the probability that
 \[\sum_{i=1}^sY_i\neq \big(p_0^{|R|}-p_s^{|R|}\big)\hat{p}^{|S\setminus R|}n\pm (\alpha_s n/2+\eps n)=\big(p_0^{|R|}-p_s^{|R|}\big)\hat{p}^{|S\setminus R|}n\pm \tfrac34\hat{p}^{|S\setminus R|}\gamma^{4D+6}\alpha_s n \]
 is at most
 \[2\exp\Big(\frac{-\eps^2n^2}{2^{31}D^4\Delta\gamma^{-4D}\delta^{-4}n+2(2D+3)\Delta\eps n}\Big)<n^{-2D-30}\,,\]
 where the last inequality is by choice of $c$.
 
 Taking the union bound over all choices of $R\subset S$ and $S$ of size at most $2D+3$, and applying Lemma~\ref{lem:dietcon}, we see that the following event has probability at most $3n^{-9}$. The pair $(H_i,H^*_0)$ is $(\alpha_i,2D+3)$-coquasirandom for each $0\le i\le s-1$, but either \RandomEmbedding\ fails to embed $G_s$ or $(H_s,H^*_0)$ is not $(\alpha_s,2D+3)$-coquasirandom. Taking now the union bound over all choices of $1\le s\le s^*$, and recalling that $(H_0,H^*_0)$ is by assumption $\big(\tfrac14\alpha_0,2D+3\big)$-coquasirandom, we conclude that the probability that for any $1\le s\le s^*$, \RandomEmbedding\ fails to embed $G_s$ or the pair $(H_s,H^*_0)$ fails to be $(\alpha_s,2D+3)$-coquasirandom is at most $1.5n^{-8}$. This completes the proof.
\end{proof}

\section{Completing the embedding}\label{sec:completion}
Recall that we complete the embedding of each graph $G_s$ by embedding the final $\delta n$ vertices using only edges of $H^*_{s-1}$. From Setting~\ref{set:graphs}, these unembedded of $G_s$ vertices form an independent set and each of them has degree $d_s$. Lemma~\ref{lem:degHstar} states that it is very likely, provided \PackingProcess\ does not fail and provided $(H_s,H^*_0)$ is coquasirandom for each $s$, that only a few edges of $H^*_0$ are used at any given vertex to form $H^*_s$, and hence $(H_s,H^*_s)$ is also coquasirandom. Complementing this, Lemma~\ref{lem:completion} states that this coquasirandomness guarantees that completing the embedding is possible. We prove these two lemmas in this section.

To prove Lemma~\ref{lem:degHstar}, we give an upper bound for the expected number of edges used at $v$ in each stage, and apply Lemma~\ref{lem:freedman} to show that the actual outcome is with high probability not much larger than this upper bound. For each $x\in V(G_s)$, we define the \emph{completion degree} of $x$, written $\deg^*(x)$, to be the degree of $x$ in the bipartite graph $G_s\big[[n-\delta n],[n]\setminus[n-\delta n]\big]$. Then the number of edges of $H^*_0$ at $v$ used in stage $s$ is $\deg^*(x)$ where $x$ is the vertex of $G_s$ embedded to $v$.
Note that since $\sum_{x=n-\delta n+1}^{n}\deg^{*}(x)=\delta n d_s$, the hand-shaking lemma tells us that
\begin{equation}\label{eq:handshake}
\sum_{x=1}^{n-\delta n}\deg^{*}(x)=\delta n d_s\;.
\end{equation}

We note that the number of edges of $H^*_{s-1}$ used in stage $s$ at any given vertex $v$ does not depend upon how the embedding of $G_s$ is completed, but only on how \RandomEmbedding\ embeds the first $n-\delta n$ vertices, so the proof of Lemma~\ref{lem:degHstar} will only need to analyse \RandomEmbedding. Indeed, if some vertex $x \in V(G_s)$, $x\le n-\delta n$ is mapped onto $v$, then this number is $\deg^*(x)$. If on the other hand, $v$ is not in the image of $G_s\big[[n-\delta n]\big]$ then $v$ will be used is the completion phase. In this case, the number of edges used at $v$ will be $d_s$ irrespective of which particular vertex $v$ will host.

\begin{proof}[Proof of Lemma~\ref{lem:degHstar}]
 Fix $v\in V(H^*_0)$. For each $s\in[s^*]$, let $Y_s$ be the number of edges of $H^*_0$ at $v$ used in stage $s$. We have
 \begin{equation}
 \label{eq:YU}
 Y_s=\sum_{x\in V(G_s)}\deg^*(x)\mathbbm{1}_{x\AlgMap v}=
 \sum_{x=1}^{n-\delta n}\deg^*(x)\mathbbm{1}_{x\AlgMap v}
 +\sum_{x=n-\delta n+1}^{n}\deg^*(x)\mathbbm{1}_{x\AlgMap v}\;.
 \end{equation} 
 
 We define $\cE$ to be the event that \PackingProcess\ succeeds and $(H_{s-1},H^*_0)$ is $(\alpha_{s-1},2D+3)$-coquasirandom for each $1\le s\le s^*$. In other words, $\cE$ is the complement of the first two events in the statement of Lemma~\ref{lem:degHstar}, so to prove Lemma~\ref{lem:degHstar} we want to show that the probability of $\cE$ occurring and the third event not occurring is very small.
 
 Suppose that $\hist_{s-1}$ is an arbitrary history of \PackingProcess{} up to and including stage $s-1$ for which $(H_{s-1},H^*_0)$ is $(\alpha_{s-1},2D+3)$-coquasirandom. We begin by estimating $\EXP(Y_s|\hist_{s-1})$.

 To estimate the desired expectation, we first aim to show
 \begin{align}
 \label{eq:degHstar:probx1} \PROB(x\AlgMap v|\hist_{s-1})&\le 5\gamma^{-D}n^{-1}&\text{if $1\le x\le n-\delta n$, and\quad}\\
 \label{eq:degHstar:probx2}
\PROB(\nexists x\in [1,n-\delta n]:x\AlgMap v|\hist_{s-1})& \le 2\delta.
 \end{align}

  In order to establish~\eqref{eq:degHstar:probx1} and~\eqref{eq:degHstar:probx2}, we need the following consequence of Lemma~\ref{lem:probvtx}. Conditioning on $\hist_{s-1}$, for each $1\le t\le n-\delta n$, the probability that \RandomEmbedding\ does not embed any of the first $t$ vertices of $G_s$ to $v$ is at most $2\tfrac{n-1-t}{n}< 2\tfrac{n-t}{n}$. This readily establishes~\eqref{eq:degHstar:probx2}.

Furthermore, under the same conditioning, by Lemma~\ref{lem:diet}, for each $1\le t\le n-\delta n$, with probability at least $1-2n^{-9}$, we have $\big|\CANDSET^{t-1}_{G_s\AlgMap H_{s-1}}(t)\big|\ge\tfrac12\gamma^D(n+1-t)$. Now, for each $1\le t\le n-\delta n$, the probability that \RandomEmbedding, conditioning on $\hist_{s-1}$, embeds $t$ to $v$ is the probability that no vertex is embedded to $v$ at time $t-1$ times the probability of picking $v$ when choosing uniformly from the candidate set of $t$. This is at most
$$2n^{-9}+2\frac{n+1-t}{n}\cdot\frac{1}{\big|\CANDSET^{t-1}_{G_s\AlgMap H_{s-1}}(t)\big|}\le 2n^{-9}+\frac{2}{\tfrac12\gamma^Dn}\;.$$ This establishes~\eqref{eq:degHstar:probx1}.

Now, we are going to substitute~\eqref{eq:degHstar:probx1} and~\eqref{eq:degHstar:probx2} into~\eqref{eq:YU}. To this end, recall that for each $x\in[n-\delta n+1,n]$ we have $\deg^*(x)=d_s$. It follows that
 \begin{align}\nonumber
  \EXP(Y_s|\hist_{s-1})&\le 5\gamma^{-D}n^{-1}\sum_{1\le x\le n-\delta n}\deg^*(x)+d_s \sum_{x=n-\delta n+1}^{n}\PROB(x\AlgMap v|\hist_{s-1})\;
  \\
\label{eq:IloveRef}
 \JUSTIFY{by \eqref{eq:handshake}, \eqref{eq:degHstar:probx2}}&\le 5\gamma^{-D}n^{-1}\cdot\delta nd_s+d_s\cdot 2\delta\le 7\gamma^{-D}D\delta\,.
 \end{align}
 Next, we obtain a similar upper bound for the second moment. Since only one vertex gets embedded to $v$, we have 
 \begin{align*}
  \EXP(Y_s^2|\hist_{s-1})&=\sum_{x\in V(G_s)}\deg^*(x)^2\PROB(x\AlgMap v|\hist_{s-1})\\
  &\le \Delta(G_s) \cdot \sum_{x\in V(G_s)}\deg^*(x)\PROB(x\AlgMap v|\hist_{s-1})=\Delta(G_s)\cdot\EXP(Y_s|\hist_{s-1})\\
  &\leByRef{eq:IloveRef} 7\gamma^{-D}D\delta\cdot\Delta(G_s)\,.
 \end{align*}
 Since $0\le Y_s\le\Delta(G_s)\le \Delta$ holds for each $s$, and since $s^*\le 2n$, we can apply Lemma~\ref{lem:freedman}, with $\rho=\delta n$ and with $\cE$ as defined above, to give
 \[\PROB\left(\cE\text{ and }\sum_{i=1}^{s^*}Y_s>50\gamma^{-D}D\delta n\right)\le \exp\big(-\tfrac{\delta^2n^2}{28\gamma^{-D}D\delta\cdot\Delta  n+2\Delta \delta n}\big)<n^{-100}\,,\]
 where the final inequality is since $\Delta= cn/\log n$ and by choice of $c$. Taking the union bound over all choices of $v$, we see that the probability that $\cE$ occurs and yet more than $50\gamma^{-D}D\delta n$ edges of $H^*_0$ are deleted at any vertex in the running of \PackingProcess\ is at most $n^{-99}$. Because the degree of each vertex in $H^*_s$ is monotone decreasing as $s$ increases, in particular this implies that the probability that there exists $1\le s\le s^*$ such that \PackingProcess\ completes stage $s$, and $(H_i,H^*_0)$ is $(\alpha_i,2D+3)$-coquasirandom for each $i<s$, yet more than $50\gamma^{-D}D\delta n$ edges of $H^*_0$ are deleted at any vertex of $H^*_s$, is at most $n^{-99}$.
 
 It remains to argue that since few edges are deleted at each vertex of $H^*_0$ to form $H^*_s$, the pair $(H_s,H^*_s)$ is coquasirandom. Suppose now that $\Delta\big(H^*_0-H^*_s\big)\le 50\gamma^{-D}D\delta n$ for some $s$, and that $(H_s,H^*_0)$ is $(\alpha_s,2D+3)$-coquasirandom. Then for any $R\subseteq S\subset V(H_s)$ with $|S|\le 2D+3$, we have
\[\big|\NBH_{H_s}(R)\cap\NBH_{H^*_0}(S\setminus R)\big|=(1\pm\alpha_s)p^{|R|}\gamma^{|S\setminus R|}n\]
 and hence
\begin{align*}
 \big|\NBH_{H_s}(R)\cap\NBH_{H^*_s}(S\setminus R)\big|&=(1\pm\alpha_s)p^{|R|}\gamma^{|S\setminus R|}n\pm (2D+3)\cdot 50\gamma^{-D}D\delta n\\
 &=\big(1\pm\eta\big)p^{|R|}\gamma^{|S\setminus R|}n
\end{align*}
where the final line is by choice of $\delta$ in~\eqref{eq:defconsts} and since $p\ge\gamma$, so that $(H_s,H^*_s)$ is $(\eta,2D+3)$-coquasirandom, as desired.
\end{proof}

Recall that Lemma~\ref{lem:completion} states that it is likely that the partial embedding $\phi_s$ of each $G_s$ provided by \RandomEmbedding\ can be extended to an embedding $\phi_s^*$ of $G_s$, with the completion edges used for the extension lying in $H^*$. Since the neighbours of each of the last $\delta n$ vertices of $G_s$ are embedded by $\phi_s$, the set of candidate vertices
\[\CANDSET_s^*(x):=\big\{v\in V(H^*_{s-1})\setminus\im\phi_s:\phi_s(y)\in \NBH_{H^*_{s-1}}(v)\text{ for each }y\in \NBH_{G_s}(x)\big\}\]
for each $x$ of these last $\delta n$ vertices in $V(H^*_{s-1})\setminus\im\phi_s$ are already fixed, and the desired $\phi_s^*$ exists if and only if there is a system of distinct representatives for the $\CANDSET^*_s(x)$ as $x$ ranges over the last $\delta n$ vertices of $G_s$. Recall that Lemma~\ref{lem:diet} states in particular that $(H^*,\im\phi_s)$ is likely to satisfy the $(2\eta,2D+3)$-diet condition, which implies both that $\CANDSET^*_s(x)$ is of size roughly $p^{d_s}\delta n$ for each of these last $x$, and also that the collection of sets is well-distributed (in a sense we will make precise later). We will see that this is almost enough to verify Hall's condition for the existence of a system of distinct representatives, but we need in addition to know that every vertex of $H^*_{s-1}-\im\phi_s$ is in sufficiently many of these candidate sets. The following lemma states that this typically is the case.

\begin{lem}\label{lem:cover*}
 Let $D\in\mathbb{N}$ and let $\gamma>0$. Let $\eta,\delta,\eps,c$ and $\alpha_x$ be as in Setting~\ref{set:graphs}. Suppose that~$G$ is a graph on vertex set $[n]$, with $\LEFTDEG(x)\le D$ for each $x\in V(G)$, with maximum degree at most $cn/\log n$ and whose last~$\delta n$ vertices all have degree $d$, where $0\le d\le D$, and form an independent set. Suppose that~$H$ is an $(\alpha_{s^*},2D+3)$-quasirandom $n$-vertex graph and that~$H^*$ is a graph on $V(H)$ with $(1\pm\eta)\gamma\binom{n}{2}$ edges such that $(H, H^*)$ forms an $(\eta, 2D+3)$-coquasirandom pair. When \RandomEmbedding\ is run to embed $G[{\scriptstyle [n-\delta n]}]$ into~$H$, with probability at least~$1-3n^{-9}$ we have for all \corrig{$v\in V(H^*)\setminus\im\psi_{n-\delta n}$ that }  %$v\in V(H^*)$
\[\Big|\big\{x\in V(G)\,:\, n-\delta n< x\le n,\;  \psi_{n-\delta n}
(\LNBH(x))\subseteq \NBH_{H^*}(v)
\big\}\Big|=(1\pm 10 D\eta )\gamma^d\delta n\;.\]
\end{lem}

The proof of this lemma is similar to the proof of Lemma~\ref{lem:cover}. 
\begin{proof}
Fix $v\in V(H^*)$ and let~$I$ be the last $\delta n$ vertices of~$G$, which by assumption form an independent set. If at any time during the run of \RandomEmbedding{} we embed a vertex to $v$, then there is nothing to prove, so we will suppose that this does not occur. Denote by $\LNBH_k(x)$ the first $k$ neighbours of $\LNBH(x)$. Let $\mathcal Y_k$ be the event that the vertices $\LNBH_k(x)$ are all embedded to $\NBH_{H^*}(v)$ for about as many $x\in I$ as one would expect, more formally that
\begin{equation}\label{eq:cover*:set}
 \Big|\big\{x\in I\::\: \psi_{n-\delta n}(\LNBH_k(x))\subseteq\NBH_{H^*}(v)\big\}\Big|=(1\pm 10 k\eta)\gamma^k\delta n\,.
\end{equation}
Let $\cB$ be the event that the $(2\eta, 2D+3)$-codiet condition fails at some time~$t\le n-\delta n$. Let \[Z_{k,t}:=\Big|\big\{x\in I\::\:\psi_{n-\delta n}(\LNBH_{k-1}(x))\subseteq \NBH_{H^*}(v)\text{ and } t\mbox{ is the $k$th vertex of }\LNBH(x)\big\}\Big|\;.\]
In other words, when we embed the vertex $t$, if it is embedded to $\NBH_{H^*}(v)$ it will add $Z_{k,t}$ more vertices to the set in~\eqref{eq:cover*:set}. Let $Y_{k,t}:=Z_{k,t}\cdot\mathbbm{1}_{\psi_{n-\delta n}(t)\in\NBH_{H^*}(v)}$. 

We want to show that if $\mathcal{Y}_{k-1}$ occurs, then $\mathcal{Y}_k$ is very likely to occur. We will then show this implies the lemma.
Observe that $\mathcal Y_k$ is the event that $\sum_{t=1}^{n-\delta n}Y_{k,t}=(1\pm 10 k\eta)\gamma^k\delta n$. Furthermore, $\mathcal Y_{k-1}$ implies that $\sum_{t=1}^{n-\delta n}Z_{k,t}=(1\pm 10 (k-1) \eta)\gamma^{k-1}\delta n$. We would like to calculate $\sum_{t=1}^{n-\delta n}\EXP(Y_{k,t}\big|\hist_{t-1})$, where $\hist_{t-1}$ denotes the embedding history of \RandomEmbedding\ up to and including embedding $t-1$. Given a time $t$, if $t$ is the $k$th vertex of $\LNBH(x)$, then at time $t-1$ the first $k-1$ vertices of $\LNBH(x)$ have already been embedded, so $Z_{k,t}$ is determined. Thus we have 
\begin{equation*}
\EXP(Y_{k,t}\big|\hist_{t-1})= \PROB\Big(\psi_t(t)\in \NBH_{H^*}(v)\big| \hist_{t-1}\Big)\cdot Z_{k,t}\;.
\end{equation*}
Suppose that at time $t-1$ we have not seen a witness that $\mathcal{B}$ fails. Then, using the $(2\eta,2D+3)$-codiet condition once with $S=\
\psi_{t-1}\big(\LNBH(t)\big)\cup \{v\}$ \corrig{(since $v$ is not in $\im\psi_{t-1}$, we indeed have $v\not\in \psi_{t-1}\big(\LNBH(t)\big)$)} and $R=\psi_{t-1}\big(\LNBH(t)\big)\subseteq S$ and once with $S=R=\psi_{t-1}\big(\LNBH(t)\big)$, we obtain
\[\PROB\Big(\psi_t(t)\in \NBH_{H^*}(v)\big| \hist_{t-1}\Big)=\frac{(1\pm 2\eta)(1\pm\eta)\gamma p^{|\LNBH(t)|}(n-t+1) }{(1\pm 2\eta) p^{|\LNBH(t)|}(n-t+1)}=(1\pm6\eta)\gamma\;.\]
Therefore, if $\overline{\cB}$ and 
$\mathcal Y_{k-1}$ hold,  we have
\[\sum_{t=1}^{n-\delta n}\EXP(Y_{k,t}\big| \hist_{t-1})=(1\pm 10(k-1)\eta)(1\pm 6\eta)\gamma^k\delta n\,.\]

Applying Lemma~\ref{lem:concentration} with $\rho=\eta\gamma^k\delta n$, we deduce that the probability that $\mathcal Y_k$ fails is very small. Indeed, the probability that $\overline{\cB}$ holds but
$\sum_{t=1}^{n-\delta n}Y_{k,t}\neq (1\pm 10 k\eta)\gamma^k\delta n$ is at most $2\exp\big(-\frac{ \eta^2\gamma^{2k}\delta^2n^2\log n}{2Dcn^2}\big)\le n^{-20}$, where we use that $Y_{k,t}\le \deg(t)$ and observe that Lemma~\ref{lem:squarebound} gives $\sum_{t=1}^{n-\delta n}\deg(t)^2\le 2D\Delta(G)n\le 2Dcn^2/\log n$.

As $\mathcal Y_0$ holds trivially with probability one, by a union bound over the choices of $k$ and $v$ we obtain that the probability that $\overline{\cB}$ holds but there is some $1\le k\le d$ for which $\mathcal Y_k$ fails is at most $2dn^{-19}$. Finally, Lemma~\ref{lem:diet} states that $\cB$ holds with probability at most $2n^{-9}$, giving the lemma statement by the union bound.
\end{proof}

We are now in a position to prove the completion lemma, Lemma~\ref{lem:completion}.

\begin{proof}[Proof of Lemma~\ref{lem:completion}]
 Suppose $H$ is an $n$-vertex $(\alpha_{s^*},2D+3)$-quasirandom graph, and $(H,H^*)$ is $(\eta,2D+3)$-coquasirandom, with $e(H)=p\binom{n}{2}$ and $e(H^*)=(1\pm\eta)\gamma\binom{n}{2}$. Let $G$ be a graph on $[n]$ with $\LEFTDEG(x)\le D$ for each $x\in[n]$ and such that the last $\delta n$ vertices of $G$ form an independent set all of whose vertices have degree $d$. When \RandomEmbedding\ is run to produce a partial embedding $\phi$ of $G$ into $H$, by Lemma~\ref{lem:diet} with probability at least $1-2n^{-9}$ the algorithm succeeds and the triple $(H,H^*,\im\phi)$ satisfies the $(2\eta,2D+3)$-diet condition. By Lemma~\ref{lem:cover*}, with probability at least $1-3n^{-9}$ in addition we have, for every vertex $v$ of $V(H^*)\setminus\im\phi$,
\begin{equation}\label{eq:hall:cover}
  \Big|\big\{x\in V(G)\,:\, n-\delta n< x\le n,\;  \phi
(\LNBH(x))\subseteq \NBH_{H^*}(v)
\big\}\Big|=(1\pm 10 D\eta )\gamma^d\delta n\;.
\end{equation}

Suppose that both good events occur, which happens with probability at least $1-5n^{-9}$. We will now show that (deterministically) this implies the existence of a system of distinct representatives for the candidate sets $\big\{\CANDSET^*(x):n-\delta n+1\le x\le n\big\}$, which trivially gives an embedding $\phi^*$ of $G$ into $H\cup H^*$ such that all edges in $[n-\delta n]$ are embedded to $H$ and the rest to $H^*$, as desired.

We prove the existence of a system of distinct representatives by verifying Hall's condition. To that end, let $X$ be a subset of $\{n-\delta n+1,\dots,n\}$. We need to show
\begin{equation}\label{eq:Hall}
 \left|\bigcup_{x\in X}\CANDSET^*(x)\right|\ge|X|\,.
\end{equation}
We separate three cases. The two easy cases are $|X|\le\tfrac12\gamma^D\delta n$ and $|X|\ge\delta n-\tfrac12\gamma^D\delta n$. For the former, if $X=\emptyset$ the statement is trivial. If not, pick any $x\in X$. We have 
\begin{equation}\label{eq:cover:CANDSET}
 \big|\CANDSET^*(x)\big|\ge (1-2\eta)(1-\eta)^d\gamma^d\delta n\ge\tfrac12\gamma^D\delta n
\end{equation}
since $\NBH_G(x)$ is a set of $d\le D$ vertices and $(H^*,\im\phi)$ satisfies the $(\eta,2D+3)$-diet condition, which in particular verifies~\eqref{eq:Hall}. For the latter, by~\eqref{eq:hall:cover} and choice of $\eta$, every vertex of $V(H^*)\setminus\im\phi$ is in more than $\tfrac34\gamma^d\delta n$ of the sets $\CANDSET^*(x)$ for $x\in\{n-\delta n+1,\dots,n\}$. In particular, every vertex $v\in V(H^*)\setminus\im\phi$ is in $\CANDSET^*(x)$ for some $x\in X$, giving~\eqref{eq:Hall}.

The final, harder, case is $\tfrac12\gamma^D\delta n<|X|<\delta n-\tfrac12\gamma^D\delta n$. Given $X$ in this size range, let $X'$ be a maximal subset of $X$ with the property $\NBH_G(x)\cap\NBH_G(x')=\emptyset$ for each $x,x'\in X'$. Since each vertex of $X'$ has $d\le D$ neighbours, the set $Y=\bigcup_{x\in X'}\NBH_G(x)$ has size at most $D|X'|$. By maximality of $X'$, every vertex in $X$ is adjacent to some vertex of $Y$. Since no vertex of $Y$ has degree more than $\Delta(G)\le cn/\log n$, we conclude
\[\tfrac12\gamma^D\delta n<|X|\le\Delta(G)|Y|\le \Delta(G) D|X'|\le cnD|X'|/\log n\,,\]
and hence $|X'|\ge \log n$ by choice of $c$ in~\eqref{eq:defconsts}. We will now argue that $Z:=\bigcup_{x\in X'}\CANDSET^*(x)$ satisfies $|Z|\ge\big(1-\tfrac12\gamma^D\big)\delta n$, which implies~\eqref{eq:Hall}.

Suppose for a contradiction that $|Z|<\big(1-\tfrac12\gamma^D\big)\delta n$. By definition, we have $\CANDSET^*(x)\subset Z$ for each $x\in X'$. We now aim to estimate the number $N$ of triples $(x,x',z)$ with $x,x'\in X$ distinct and $z\in Z$ satisfying $z\in\CANDSET^*(x)\cap\CANDSET^*(x')$. For each $z$, let $d_z=\big|\{x\in X':z\in\CANDSET^*(x)\}\big|$. Using Jensen's inequality (since $\binom{\cdot}{2}$ is convex), we have
\begin{align*}
 N=\sum_{z\in Z}\binom{d_z}{2}&\ge |Z|\cdot\binom{|Z|^{-1}\sum_{z\in Z}d_z}{2}\\
 \JUSTIFY{by~\eqref{eq:cover:CANDSET}}&\ge |Z|\cdot\binom{|Z|^{-1}|X'|(1-2\eta)(1-\eta)^d\gamma^d\delta n}{2}\\
 &=\tfrac12|X'|(1-2D\eta)\gamma^d\delta n\big(|Z|^{-1}|X'|(1-2D\eta)\gamma^d\delta n-1\big)\\
 &\ge\tfrac12(1-2D\eta)^3|X'|^2|Z|^{-1}\gamma^{2d}\delta n\\
 &\ge\tfrac12(1-2D\eta)^3|X'|^2\big(1-\tfrac12\gamma^D\big)^{-1}\gamma^{2d}n\,,
\end{align*}
where the penultimate inequality holds since $|Z|<\delta n$ and $|X'|\ge\log n$ is sufficiently large, and the final inequality uses our assumed upper bound on $|Z|$.
On the other hand, since $\NBH_G(x)$ and $\NBH_G(x')$ are disjoint, we have
\[N=\sum_{x,x'\in X'}\big|\CANDSET^*(x)\cap\CANDSET^*(x')\big|\le\tbinom{|X'|}{2}(1+2\eta)(1+\eta)^{2d}\gamma^{2d}\delta n\le\tfrac12|X'|^2(1+4D\eta)\gamma^{2d}\delta n\]
using the $(2\eta,2D+3)$-diet condition which $(H^*,\im\phi)$ satisfies. We conclude
\[\tfrac12(1-2D\eta)^3|X'|^2\big(1-\tfrac12\gamma^D\big)^{-1}\gamma^{2d}n\le\tfrac12|X'|^2(1+4D\eta)\gamma^{2d}\delta n \]
which is false since by choice of $\eta$ in~\eqref{eq:defconsts} we have $(1-2D\eta)^3(1+4D\eta)^{-1}>1-\tfrac12\gamma^D$. Thus~\eqref{eq:Hall} holds for all $X$, so the desired $\phi^*$ exists.
\end{proof}

\section{Concluding remarks}
\label{sec:concluding}

\subsection{Constants in Theorem~\ref{thm:MAINunbounded}}
Given~$\gamma$ and $D$ in Theorem~\ref{thm:MAINunbounded}, the constant $c$ is
set in Setting~\ref{set:graphs}. All the dependencies in~\eqref{eq:defconsts}
are polynomial, except for the exponentials used to define $C$ and $\alpha_x$.
As a result, $c$ depends roughly doubly-exponentially on $D$ and $\gamma$, more
precisely $c\approx \exp(-\exp(D^{5+o(1)}\cdot \gamma^{-24D-10+o(1)}))$ (where
$o(1)\rightarrow 0$ as $D,1/\gamma\rightarrow \infty$). This of course puts an
implicit requirement on $n_0$, as instances of the result for which the maximum
degree bound $\frac{cn}{\log n}$ are less than~$1$ are vacuous. 

% The only other
% bound on $n_0$ comes from the analysis of our randomised algorithm; in
% particular, we need that the lower bound of $1-(n^{-50} + s^* \cdot 5n^{-9})$ on
% the probability of the algorithm succeeding, which appears in the proof of
% Theorem~\ref{thm:maintechGENERAL}, is positive. For that, $n_0\ge 2$ suffices.\marginpar{We don't ever claim that some log eventually beats some constant or such a thing..? If not I suggest we don't mention the lower bound $2$. -PA}

By way of brief comparison with other recent packing results, we believe most of
the results we cited earlier obtain broadly similar or better constant
dependencies to our results (though these bounds are generally not given
explicitly and we did not check carefully), unless the Regularity Lemma is used.

\subsection{Limits of the method}\label{ssec:limits}
As Ferber and Samotij~\cite{FeSa:PackingTrees} point out, a randomised strategy such as the one we use here will not succeed in packing graphs with many vertices of degree $\omega\big(\tfrac{n}{\log n}\big)$, because it is likely to put these vertices unevenly into the host graph and after packing only half the guest graphs one vertex will probably have degree substantially less than the average. If the remaining graphs are for example Hamilton cycles, this vertex will become a bottleneck which causes the strategy to fail. One might try to pick vertices non-uniformly in order to correct such imbalances as they form, but analysing such a strategy would be challenging and it is not clear that it would work: common neighbourhoods of several vertices will also occasionally be far from the expected size.

Although it might well be that we can obtain near-perfect packings of graphs with degeneracy much bigger than $\log n$ into $K_n$, any strategy like the one we use here will certainly not succeed in doing so. The reason is simply that strategies like ours work by maintaining quasirandomness, and hence work equally well starting with a dense random graph rather than the complete graph. Take $H$ to be a clique of order $3\log_2 n$. Then a well-known calculation shows that $\mathbb G\big(n,\tfrac12\big)$ typically does not even contain one copy of $H$.

We have not tried to analyse our approach more carefully in order to work with sparse random or quasirandom graphs. We are confident that (with substantially more work, and using ideas from~\cite{SparseBU}) one could prove a near-perfect packing result for typical $\mathbb G\big(n,p\big)$, where $p>n^{-\eps}$ for some $\eps>0$ depending on the degeneracy bound $D$. But we suspect that our approach would not then allow for maximum degrees of the guest graphs as large as $\Omega(pn/\log n)$, even if we asked only to pack almost-spanning graphs, and certainly we cannot take $\eps$ as big as $\tfrac{1}{2D+3}$, since at this point $\mathbb G(n,p)$ itself is typically not $(\tfrac12,2D+3)$-quasirandom. In particular, our approach cannot challenge the tree packing results of~\cite{FeSa:PackingTrees} in sparse random graphs.

\subsection{Perfect packings}
It is easy to check that the graph of uncovered edges in the packing of Theorem~\ref{thm:maintechGENERAL} is $(2\eta,2D+3)$-quasirandom, and $\eta$ can be chosen arbitrarily small by increasing $D$ if necessary. In particular, this means that the result of Joos, Kim, K\"uhn and Osthus~\cite{JoKiKuOs:Packing} applies to this leftover. Thus we can extend the result of~\cite{JoKiKuOs:Packing} on the Tree Packing Conjecture to allow many trees where the maximum degree is bounded only by $\tfrac{cn}{\log n}$, provided that it is bounded by $D$ in the remainder. This is however a rather peculiar condition. 

\section{Acknowledgements}
We thank Pavel Hladk\'y for looking after our little children while we were proving these results. We thank the anonymous referees for their very helpful comments.

\medskip{}

The contents of this publication reflects only the authors' views
and not necessarily the views of the European Commission of the European
Union.
\bibliographystyle{plain}
\bibliography{bibl}

\begin{thebibliography}{10}

\bibitem{AAGH:AlmostAllTrees}
A.~Adamaszek, P.~Allen, C.~Grosu, and J.~Hladk{\'y}.
\newblock Almost all trees are almost graceful.
\newblock arXiv:1608.01577.

\bibitem{SparseBU}
P.~Allen, J.~B{\"o}ttcher, H.~H{\`a}n, Y.~Kohayakawa, and Y.~Person.
\newblock Blow-up lemmas for sparse graphs.
\newblock arXiv:1612.00622.

\bibitem{Balogh2013}
J.~Balogh and C.~Palmer.
\newblock {On the {T}ree {P}acking {C}onjecture}.
\newblock {\em SIAM J. Discrete Math.}, 27(4):1995--2006, 2013.

\bibitem{barber2015fractional}
B.~Barber, D.~K\"uhn, A.~Lo, R.~Montgomery, and D.~Osthus.
\newblock Fractional clique decompositions of dense graphs and hypergraphs.
\newblock {\em J. Combin. Theory Ser. B}, 127:148--186, 2017.

\bibitem{Bollobas1983}
B.~Bollob{\'a}s.
\newblock {Some remarks on packing trees}.
\newblock {\em Discrete Math.}, 46(2):203--204, 1983.

\bibitem{Boettcherb}
J.~B{\"o}ttcher, J.~Hladk{\'y}, D.~Piguet, and A.~Taraz.
\newblock An approximate version of the tree packing conjecture.
\newblock {\em Israel J. Math.}, 211(1):391--446, 2016.

\bibitem{Chung1989}
F.~R.~K. Chung, R.~L. Graham, and R.~M. Wilson.
\newblock {Quasi-random graphs}.
\newblock {\em Combinatorica}, 9(4):345--362, 1989.

\bibitem{CsKueLoOsTr:HamDec}
B.~Csaba, D.~K\"uhn, A.~Lo, D.~Osthus, and A.~Treglown.
\newblock Proof of the $1$-factorization and {H}amilton decomposition
  conjectures.
\newblock {\em Mem. Amer. Math. Soc.}, 244(1154):170, 2016.

\bibitem{DorTar}
D.~Dor and M.~Tarsi.
\newblock Graph decomposition is {NP}-complete: a complete proof of {H}olyer's
  conjecture.
\newblock {\em SIAM J. Comput.}, 26(4):1166--1187, 1997.

\bibitem{FeLeMou:PackingSpanning}
A.~Ferber, C.~Lee, and F.~Mousset.
\newblock Packing spanning graphs from separable families.
\newblock {\em Israel J. Math.}, 219(2):959--982, 2017.

\bibitem{FeSa:PackingTrees}
A.~Ferber and W.~Samotij.
\newblock Packing trees of unbounded degrees in random graphs.
\newblock arXiv:1607.07342, to appear in J. Lond. Math. Soc.

\bibitem{Freedman}
D.~A. Freedman.
\newblock On tail probabilities for martingales.
\newblock {\em Ann. Probability}, 3:100--118, 1975.

\bibitem{glock2016decomposition}
S.~Glock, D.~K{\"u}hn, A.~Lo, R.~Montgomery, and D.~Osthus.
\newblock On the decomposition threshold of a given graph.
\newblock arXiv:1603.04724, to appear in J. Combin. Theory Ser. B.

\bibitem{GKLO:Designs}
S.~Glock, D.~K{\"u}hn, A.~Lo, and D.~Osthus.
\newblock The existence of designs via iterative absorption.
\newblock arXiv:1611.06827.

\bibitem{GKLO:Fdesigns}
S.~Glock, D.~K{\"u}hn, A.~Lo, and D.~Osthus.
\newblock {Hypergraph $F$-designs for arbitrary $F$}.
\newblock arXiv:1706.01800.

\bibitem{Gyarfas1978}
A.~Gy{\'a}rf{\'a}s and J.~Lehel.
\newblock {Packing trees of different order into {$K_n$}}.
\newblock In {\em {Combinatorics (Proc. Fifth Hungarian Colloq., Keszthely,
  1976)}}, volume~18 of {\em {Colloq. Math. Soc. J{\'a}nos Bolyai}}, pages
  463--469. North-Holland, Amsterdam, 1978.

\bibitem{Hobbs1981}
A.~M. Hobbs.
\newblock {Packing trees}.
\newblock In {\em {Proceedings of the {T}welfth {S}outheastern {C}onference on
  {C}ombinatorics, {G}raph {T}heory and {C}omputing, {V}ol. {II} ({B}aton
  {R}ouge, {L}a., 1981)}}, volume~33, pages 63--73, 1981.

\bibitem{JoKiKuOs:Packing}
F.~Joos, J.~Kim, D.~K{\"u}hn, and D.~Osthus.
\newblock Optimal packings of bounded degree trees.
\newblock arXiv:1606.03953, to appear in J. Eur. Math. Soc.

\bibitem{Kee:ExistenceOfDesigns}
P.~Keevash.
\newblock The existence of designs.
\newblock arXiv:1401.3665.

\bibitem{KiKuOsTy:Packing}
Jaehoon Kim, Daniela K\"{u}hn, Deryk Osthus, and Mykhaylo Tyomkyn.
\newblock A blow-up lemma for approximate decompositions.
\newblock {\em Trans. Amer. Math. Soc.}, 371(7):4655--4742, 2019.

\bibitem{Kirkman}
T.~P. Kirkman.
\newblock On a problem in combinations.
\newblock {\em Cambridge and Dublin Math. J.}, 2:191--204, 1847.

\bibitem{MeRoSch:PackingMinorClosed}
S.~Messuti, V.~R\"odl, and M.~Schacht.
\newblock Packing minor-closed families of graphs into complete graphs.
\newblock {\em J. Combin. Theory Ser. B}, 119:245--265, 2016.

\bibitem{Pluecker}
J.~Pl\"ucker.
\newblock {\em System der analytischen {G}eometrie, auf neue
  {B}etrachtungsweisen gegr\"undet, und insbesondere eine ausf\"uhrliche
  {T}heorie der {C}urven dritter {O}rdnung enthaltend}.
\newblock Duncker und Humboldt, Berlin, 1835.

\bibitem{RCW}
D.~K. Ray-Chaudhuri and R.~M. Wilson.
\newblock Solution of {K}irkman's schoolgirl problem.
\newblock In {\em Combinatorics ({P}roc. {S}ympos. {P}ure {M}ath., {V}ol.
  {XIX}, {U}niv. {C}alifornia, {L}os {A}ngeles, {C}alif., 1968)}, pages
  187--203. Amer. Math. Soc., Providence, R.I., 1971.

\bibitem{Ringel1963}
G.~Ringel.
\newblock {Problem 25}.
\newblock In {\em {Theory of Graphs and its Applications (Proc. Int. Symp.
  Smolenice 1963)}}. Czech. Acad. Sci., Prague, 1963.

\bibitem{Steiner}
J.~Steiner.
\newblock Combinatorische aufgabe.
\newblock {\em Journal f{\"u}r die reine und angewandte Mathematik},
  45:181--182, 1853.

\bibitem{Yus:Sur}
R.~Yuster.
\newblock Combinatorial and computational aspects of graph packing and graph
  decomposition.
\newblock {\em Comp. Sci. Rev.}, 1(1):12--24, 2007.

\end{thebibliography}

\end{document}